\documentclass[twoside]{article}
\usepackage{RR}
\usepackage[T1]{fontenc}
\usepackage[utf8]{inputenc}
\usepackage[square,sort,comma,numbers]{natbib}
\usepackage{graphicx} 
\usepackage[colorlinks=true, allcolors=blue]{hyperref}
\usepackage[a4paper]{geometry} 
\usepackage{amsfonts}
\usepackage{amsmath}
\usepackage{amssymb}
\usepackage{amscd}
\usepackage{amsthm}
\usepackage{mathrsfs}
\usepackage{stmaryrd}
\usepackage{enumitem}
\usepackage{fancybox}
\usepackage{datetime}
\usepackage{algorithm}
\usepackage{algorithmic}
\usepackage{tikz}
\usetikzlibrary{shapes,arrows,cd}
\input xy
\xyoption{all}
\def\A{\mathcal{A}}
\def\B{\mathcal{B}}
\def\C{\mathbb{C}}

\def\E{\mathbb{E}}
\def\I{\mathbb{I}}
\def\IC{\mathcal{I}}
\def\K{\mathbb{K}}

\def\L{\mathcal{L}}
\def\M{\mathcal{M}}
\def\N{\mathbb{N}}
\def\NC{\mathcal{N}}
\def\O{\mathbb{O}}
\def\OC{\mathcal{O}}
\def\P{\mathbb{P}}
\def\PC{\mathcal{P}}
\def\R{\mathbb{R}}

\def\S{\mathbb{S}}
\def\U{\mathbb{U}}
\def\X{\mathcal{X}}
\def\a{\textsc{a}}
\def\b{\textsc{b}}
\def\c{\textsc{c}}
\def\ee{\textsc{e}}
\def\f{\textsc{f}}

\def\h{\textsc{h}}
\def\i{\textsc{i}}
\def\j{\textsc{j}}
\def\m{\textsc{m}}
\def\n{\textsc{n}}
\def\p{\textsc{p}}
\def\q{\textsc{q}}
\def\r{\textsc{r}}
\def\t{\textsc{t}}
\def\corr{\mathrm{corr}\:}
\def\cov{\mathrm{cov}\:}
\def\diag{\mathrm{diag}\:}
\def\GL{\mathbb{G}\mathbb{L}}
\def\im{\mathrm{im}\:}
\def\ker{\mathrm{ker}\:}
\def\Log{\mathrm{Log}\:}
\def\ones{\mathbf{1}}
\def\rank{\mathrm{rank}\:}
\def\SO{\mathbb{S}\mathbb{O}}
\def\Sp{\mathrm{Sp}\:}
\def\sp{\mathrm{sp}\:}
\def\span{\mathrm{span}\:}
\def\st{\mathrm{s.t.}}
\def\std{\mathrm{std}\:}
\def\Tr{\mathrm{Tr}\:}
\def\var{\mathrm{var}\:}

\def\zeros{\mathbf{0}}
\def\CCA{\texttt{CCA }}
\def\CoA{\texttt{CoA }}

\def\MCoA{\texttt{MCoA }}
\def\MDA{\texttt{MDA }}
\def\MDS{\texttt{MDS }}
\def\PCA{\texttt{PCA }}
\def\PCAda{\texttt{PCAda }}
\def\PCAiv{\texttt{PCAiv }}
\def\PCAmet{\texttt{PCAmet }}
\def\PCAsc{\texttt{PCAsc }}

\def\numpy{\textsl{numpy}}
\def\numpy_{\textsl{numpy }}
\def\nB{\vspace*{5mm}\noindent}
\def\nD{\vspace*{5mm}\noindent$\diamond\quad$}
\def\nS{\vspace*{5mm}\noindent$\bullet\quad$}
\def\nT#1{\vspace*{5mm}\noindent$\bullet\quad$\textbf{#1:}}
\def\notes{\paragraph{Notes and references:}}

\usepackage{todonotes}

\usepackage{amsthm}

\numberwithin{equation}{subsection}

\RRNo{9488}
\RRdate{May 2023 (2nd revision)}

\RRauthor{%
  Alain Franc%
  \thanks{Pleiade team and INRAE, Biogeco, University of Bordeaux,
    69, route d'Arcachon, 33610, Cestas}%
  \thanks{Correspondence: \texttt{alain.franc@inrae.fr}}
}

\authorhead{Franc}
\RRtitle{%
  Méthodes linéaires de Réduction de Dimension
}
\RRetitle{%
  Linear Dimensionality Reduction
  }
\titlehead{Linear Dimensionality Reduction}

\RRresume{%
  Ce document brosse un panorama des méthodes linéaires de l'Analyse de données multivariées. Il s'agit d'un domaine ancien et classique, bien établi depuis les années 60, et redevenu d'actualité en tant qu'étape clé dans l'apprentissage statistique. On peut considérer ces méthodes comme faisant partie d'une approche algébrique de l'apprentissage statistique ou bien comme une réduction de dimension avec une tonalité plus géométrique. Ces deux approches sont étroitement liées : il est plus facile d'apprendre des patterns des données dans des espaces à faible dimension que dans des espaces à grande dimension. Nous montrons comment une apparente diversité de méthodes et outils peut se ramener à une seule méthode : l'Analyse en Composantes Principales, avec la  \texttt{SVD} (Singular Value Decomposition), de telle sorte que les efforts d'optimisation des codes pour l'analyse de jeux de données massives peut se focaliser sur cette méthode centrale partagée, au bénéfice de toutes les méthodes. Une extension à l'étude de plusieurs tableaux est présentée (Analyse canonique).
}
\RRabstract{%
 These notes  are an overview of some classical linear methods in Multivariate Data Analysis. This is a good old domain, well established since the 60's, and refreshed timely as a key step in statistical learning. It can be presented as part of statistical learning, or as dimensionality reduction with a geometric flavor. Both approaches are tightly linked: it is easier to learn patterns from data in low dimensional spaces than in high-dimensional spaces. It is shown how a diversity of methods and tools boil down to a single core methods, \PCA with \texttt{SVD}, such that the efforts to optimize codes for analyzing massive data sets like distributed memory and task-based programming or to improve the efficiency of the algorithms like Randomised SVD  can focus on this shared core method, and benefit to all methods. 
}
\RRmotcle{%
  Réduction de dimension,
  Analyse de données multivariées,
  Apprentissage statistique,
  Analyse en Composantes Principales,
  Analyse Factorielle des Correspondances,
  Analyse avec variables instrumentales,
  Analyse canonique
}
\RRkeyword{%
  Dimensionality reduction,
  Multivariate Data Analysis,
  Statistical Learning,
  Principal Components Analysis,
  Correspondence Analysis,
  Analysis with Instrumental Variables,
  Canonical Analysis
}

\RRprojet{Pleiade}
\RCBordeaux

\begin{document}
\makeRR

%
\section*{Motivation}
%

The writing of this document is motivated by the convergence of several observations.

\nD The diversity of methods described as Data Analysis in the 1970s, nowadays attached to Statistical Learning, can be organised in a circular way where one is a variation of the other via certain choices. This global approach is essentially algebraic, based on matrix calculation, and has been developed by a French school, in parallel with more statistical approaches in Anglo-Saxon countries.

\nD The main methods are PCA (Principal Component Analysis) associating items and features, CoA (Correspondence Analysis) of contingency tables, Canonical Analysis for analysis of two tables, Multiple Correspondence Analysis for analysis of two tables or more, MDS (Multidimensional Scaling) for building a point cloud from a table of distances, etc...  These methods had fallen into disuse as descriptive methods, but have experienced a revival of interest in exploring structures in massive data (which is sometimes called ``pattern discovery'' and is related to unsupervised learning), mainly because statistical learning methods are sometimes inoperative in spaces of very large dimension. A pretreatment as a projection on a space of lower dimension can make them operational, provided the projection respects the structure of the dataset.

\nD The starting point for each of these methods is a data array, which can be a cross-tabulation between objects and variables, a contingency table, a distance or dissimilarity matrix, etc. ... A key observation is that each of the methods is organised according to the triptych:
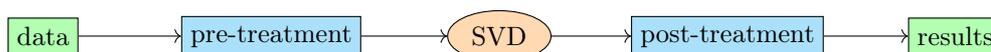
\begin{center}
        \begin{tikzpicture}
            \node[rectangle,draw,fill=green!30](dat) at (0,0){data};            \node[rectangle,draw,fill=cyan!30](pre) at (3,0){pre-treatment};
            \node[ellipse,draw,fill=orange!30](meth) at (6,0){SVD};
            \node[rectangle,draw,fill=cyan!30](post) at (9,0){post-treatment};
            \node[rectangle,draw,fill=green!30](res) at (12,0){results};
    \draw[->] (dat)--(pre);
   \draw[->] (pre)--(meth);
   \draw[->] (meth)--(post);
    \draw[->] (post)--(res);
        \end{tikzpicture}
    \end{center}
where each method can be read as a diabolo, where a pre-processing of the data constructs a matrix, which is itself decomposed via a Singular Value Decomposition (SVD), the outputs of which are post-processed to produce the desired result. The SVD step is generally the limiting step for scaling up, i.e. processing massive data from very large arrays: the complexity is cubic with the size of the arrays.

\nB Significant progress has been made recently in scaling up the SVD by combining three elements:
\begin{itemize}
    \item an evolution of the SVD computation algorithm by including the Gaussian Random Projection (rSVD, see \cite{Bingham2001,Halko2011}),
    \item a distributed memory implementation of the basic matrix operations,
    \item the use of a task-based programming paradigm for the assembly of these steps.
\end{itemize}

\nB Random projection if often presented itself as a dimensionality reduction method, by projecting a point cloud on a space of lower dimension while respecting the pairwise distances (see e.g. \cite{Bingham2001}). In these notes, it is used as a tool for lowering the complexity of the calculation of the SVD of a large matrix, as in \cite{Halko2011}. The integration of the rSVD in the MDS algorithm has been made in \cite{Paradis2018,blanchard:hal-01685711}. Numerical implementation for very large matrices (namely $10^6 \times 10^6$) with distributed-memory and task-based programming has been made in \cite{agullo:hal-03773985}. A motivation for these notes is to show how such an approach developed for MDS can be operational for any linear dimensionality reduction method that relies on a SVD.

\paragraph{Acknowledgements:} These notes have been written as a methodological companion to two ADT (Action de Développement Technologique) built on collaborations between INRIA and Inrae at Bordeaux: Gordon (2019-2020) and Diodon (2021-2022), the aim of which was to provide librairies for running linear dimension reduction on massive data sets \cite{agullo:hal-03773985} with rSVD, possibly (in C++) distributed memory and task based programming. Most of the pseudocodes presented in these notes, including the integration of Randomized SVD, have been implemented in
\begin{itemize}[label=$\diamond$]
\item a C++ library written by INRIA, called \texttt{cppdiodon}, publicly available at \\ \url{https://diodon.gitlabpages.inria.fr/cppdiodon/index.html}, 
\item a python library written by Inrae, called \texttt{pydiodon}, publicly available at \\ \url{https://gitlab.inria.fr/diodon/pydiodon}.
\end{itemize}
I am particularly grateful to Emmanuel Agullo, Pierre Blanchard, Olivier Coulaud, Jean-Marc Frigerio, Romain Péressoni, Florent Pruvost for many discussions throughout these projects and before, especially on linear algebra, and their encouragements to make explicit this methodological companion which enabled the continuation of the Gordon ADT through the Diodon ADT.
I am particularly indebted to Francis Cailliez, Daniel Chessel, Jean-Baptiste Denis, Yves Escoufier, Jean-Dominique Lebreton and Robert Sabatier, who taught me multivariate data analysis many decades ago.

\clearpage

\tableofcontents

\clearpage

%
\section{Introduction}\label{sec:semantics}
%

Let us start with an example: supervised learning with Support Vector Machine (SVM). Imagine a training set of $n$ observations, each observation being a pair $(x,y)$, with $x \in \R^p$ and $y \in \{-1,1\}$. One wishes to predict $y$ as an outcome of a new observation $x$, not in the training set, and where $y$ is unknown. Avoiding technicalities, and keeping the presentation short for this introduction, this can be done in adequate situations when there exists a linear discriminating function $f(x) = \beta + \sum_iw_ix_i$, and that $y =1$ if $f(x) >0$ and $y= -1$ if $f(x) < 0$. A separating hyperplane is an hyperplane such that all points $x \in \R^p$ of pairs $(x,1)$ (say, blue points) are on one side and all points $x$ of pairs $(x,-1)$ (say red points) are on the other side of the hyperplane. Such a function $f$ is not unique, and here support vectors come into the game. The margin is defined as the minimal distance between the points $x_i$ of the training set and the separating hyperplane $f(x)=0$. SVM is finding a linear function $f$ with maximum margin. This is equivalent to finding two parallel separating hyperplanes with maximal mutual distance. If the dimension $p$ of the space where the observations $x$ are given is large, this can lead to high computation load. However, there is a lemma\footnote{Which deserves to be called a theorem, but is is referred to classically as a lemma} called Johnson-Lindenstrauss lemma (J.-L.) telling that there exists a space $H$ of smaller dimension $d \ll p$ (of $\OC(\mathrm{Log} \: p)$) such that all pairwise distances are preserved up to a high accuracy while projecting on $H$. Therefore, SVM can be fit on the projection of the training set on a space of much lower dimension. It appears that J-L lemma is at the heart of very efficient methods to perform Singular Value Decomposition of very large matrices, and SVD is at the heart of most of, if not all, linear dimension reduction methods in multivariate data analysis. This establishes a tight link between machine learning and multivariate data analysis.\\  
\\
Multivariate Data Analysis (MDA) could be read, and is presented here, as an algebraic construction in linear algebra, around Singular Value Decomposition. This is deliberate, as we wish to focus on recent progress in High Performance Computing to handle very large data sets, and this requires a sound basis in linear algebra. Let us keep in mind that data distinguish statistics from probability: statistics are about inference of some models from data. MDA is about inferring some patterns from data, like correlation structure between items and features, or indviduals and variables. As such, it is now part of Data Mining, Statistical Learning and Machine Learning. This can be summarized with the subtitle of \cite{Hastie2009}: data mining, inference and prediction. Tibshirani has provided a dictionary between statistics and machine learning\footnote{See \url{http://www-stat.stanford.edu/~tibs/stat315a/glossary.pdf}}, referred to in \cite{Murphy2012}, partly given here:\\
\\
\begin{center}
\ovalbox{
 \begin{tabular}{ll}
  \textbf{Machine learning} & \textbf{Statistics} \\
  \hline
  weights & parameters \\
  learning & fitting \\
  supervised learning & regression/classification \\
  unsupervised learning & density estimation, clustering \\
 \end{tabular}
 }
\end{center}

\nB Beyond this dictionary, there is genuine innovation in machine learning, which is about inferring meaningful pattern in data in a very elaborate way, far beyond data analysis (see e.g. \cite{Murphy2012,Shalev2014}). A grail is to mimic the way a brain learns from experience. This is beyond the scope of those notes.

\nB MDA is a tool for learning patterns about correlations between features and items, which is one type among many others of data structure. Supervised or unsupervised learning may rely, at one step or another, on techniques inherited from multivariate data analysis, in a same way that multivariate data analysis relies at one step or another on tools inherited from linear algebra. This can be formalized by the following succession of steps:\\
\\
\tikzstyle{block} = [rectangle, draw, fill=blue!20, text width=6em, text centered, rounded corners, minimum height=4em]
\tikzstyle{line} = [draw, -latex'] 
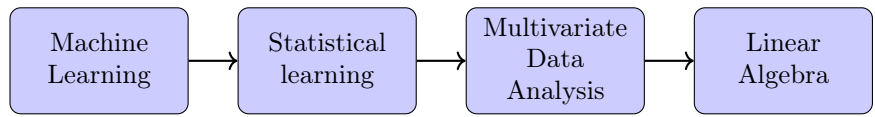
\begin{center}
\begin{tikzpicture}
\node [block] (ml) {Machine Learning};
\node [block, right of=ml, node distance = 3cm] (sl) {Statistical learning}; 
\node [block, right of=sl, node distance = 3cm] (mda) {Multivariate Data Analysis}; 
\node [block, right of=mda, node distance = 3cm] (la) {Linear Algebra}; 
\draw[->, thick] (ml) -- (sl) ;
\draw[->, thick] (sl) -- (mda) ;
\draw[->, thick] (mda) -- (la) ;
\end{tikzpicture}
\end{center}

\vspace*{3mm}\noindent We have this inheritance in mind for these notes, with an objective of implementation of calculations for very large data sets.\\
\\
There is a second reason for emphasizing the role of MDA in Machine Learning. Murphy recognized two approaches in Machine Learning \cite[sec. 1.1.2]{Murphy2012}, which are consistent with the synthetic table of Tibshirani:
\begin{itemize}
 \item a predictive or supervised approach, where a response variable is predicted from some features, from a set of observations where the response is known, called the training set
 \item a descriptive or unsupervised approach, where the objective is to find some interesting patterns, and is referred to as ``pattern discovery''.
\end{itemize}
Many of these techniques work well in low dimension, i.e. when the number of features to work with is small. MDA plays a key role in such a framework to produce a low dimension representation of an array connecting items and features, as close as possible to the original array. This is called Dimensionality Reduction (see e.g. \cite{Lee2007}), and one of the key (linear) technique therefore is Principal Component Analysis (PCA, see section \ref{sec:pca}). So, one can say that MDA paves the way for elaborate machine learning processes.

\notes There exists many excellent surveys for learning patterns from data.
See e.g. \cite{Cristianini2000,Murphy2012,Shalev2014}. For a concise introduction to SVM, see \cite{Cristianini2000} or \cite{Smola2004}.

%
\section{Multivariate Data Analysis}\label{sec:mda}
%

MDA is at the crossroads of three domains:
\begin{itemize}[label=$\rightarrow$]
 \item It is about finding structures in data presented as arrays. This is algebra. 
 \item It is possible to attach to an array a point cloud in a Euclidean space, and study the shape of the cloud. This is geometry.
 \item Data are modeled as realizations of random variables. This is statistics.
\end{itemize}
Hence, \MDA is at the crossroads between algebra, geometry and statistics.

\nB As arrays of data are matrices, \MDA heavily relies on Linear Algebra. Producing a matrix $A$ is one of the most classical mathematical formalization of some information gathered on a set of items. Let us consider a set of $n$ items each characterized by $p$ variables. The rows $i \in \llbracket 1,n \rrbracket$ are the items, and the columns $j \in \llbracket 1,p \rrbracket$ are some variables, often referred to as features in machine learning. The value of the feature $j$ for the item $i$ is the coefficient $\alpha_{ij}$ of the matrix. One objective is to describe how items and features are related, which is usually addressed by a low rank approximation of matrix $A$.
 
\nB A point cloud is a geometric object associated to such a feature matrix. Provided the variables are quantitative, i.e. numbers in $\R$, a set of $n$ points in $\R^p$ is built from $A$, with one point $a_i = (\alpha_{i1}, \ldots, \alpha_{ip}) \in \R^p$ for item $i$. This point cloud as a geometric object is denoted $\A$. In many real cases, the points are located in a low dimensional manifold. Finding such a manifold is called \emph{dimensionality reduction}. Efficient and well understood techniques exist when the manifold is linear (or affine): the best approximation of $\A$ by a projection in a space of dimension $r$ can be solved by finding the best approximation of $A$ by a matrix of rank $r$.

\nB Row $i$ of matrix $A$ or point $a_i \in \A$ can be modeled as the realization of a set of random variables $(X_1, \ldots, X_p)$. Questions of interest are the study of the dependence structure of the $X_j$, given by the variance-covariance matrix of observed data, or exhibiting low dimension latent variables $z$, such that observations $x$ can be modeled as $\P(x\mid z)$.

\nB Many of these techniques have been progressively selected as core techniques along several decades since the beginning of 20th century. Many of these methods are presented and studied within the algebraic framework of linear algebra. For example, PCA can be built as a consequence of the Singular Value Decomposition of $A$: $A = U\Sigma V^\t$. These methods are constantly co-evolving with numerical linear algebra, and have benefited from\footnote{at least ...} two revolutions:
\begin{description}
 \item[computing revolution:] the development of computing infrastructures especially in the 70's has led to the mushrooming of scientific libraries first for mainframes, and later for a range of machines from laptops to computing clusters 
 \item[massive data revolution:] many sensors produce now myriads of bytes of data, like telescopes, satellites, sequencers, etc., raising the challenge to overcome the walls of time and memory while implementing  those methods on massive data sets, leading to matrices of very large dimensions (like $10^5$ to $10^6$ rows or columns).
\end{description}
The progress in these domains are due simultaneously to the derivation of new algorithms (like the random projection method for computing the Singular Value Decomposition of a large dense matrix, see \cite{Halko2011}), development of new paradigms for implementing some algorithms (like Message Passing Interface for distributed-memory parallelization), and progress in technology of computing infrastructures (like Graphic Processing Units).\\
\\
The high diversity of existing method can be organized into a small set of iconic methods, knowing that such a classification is far from being either unique or universally adopted. However, most textbooks presenting these methods progressively reach an agreement on the poles around which to organize their variety and similarities/dissimilarities. We will not discuss this here, but select among many possibilities a small set of methods, organized along a variety of questions they answer. These methods are presented in the following table.\\ 
\\
\begin{center}
\ovalbox{
\begin{tabular}{cll}
\textbf{Acronym} & \textbf{Method} & \textbf{Section} \\
&&\\
\PCA & Principal Component Analysis & \ref{sec:pca}\\
\PCAsc & Scaled-Centered \PCA & \ref{sec:pca:classical}\\
\PCAda & \PCA with double averaging & \ref{sec:pca:classical} \\
\PCAiv& PCA with instrumental variables & \ref{sec:pcaiv}\\
\PCAmet & PCA with metrics & \ref{sec:pcamet} \\
\CoA & Correspondence Analysis & \ref{sec:coa} \\
\CCA & Canonical Correlation Analysis & \ref{sec:cca} \\
\MCoA & Multiple Correspondence Analysis & \ref{sec:mcoa} \\
\MDS & Multidimensional Scaling & \ref{sec:mds} \\
\end{tabular}
}
\end{center}

\vspace*{5mm} 
\noindent The organization followed here is to analyze each of these methods as a pipeline\\
\\
\tikzstyle{block} = [rectangle, draw, fill=blue!20, text width=6em, text centered, rounded corners, minimum height=4em]
\tikzstyle{line} = [draw, -latex'] 
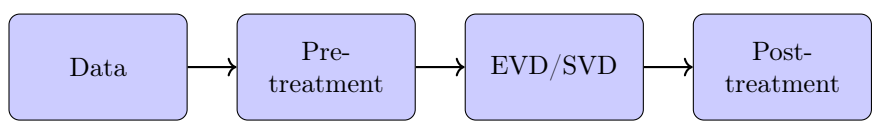
\begin{center}
\begin{tikzpicture}
\node [block] (data) {Data};
\node [block, right of=data, node distance = 3cm] (pre) {Pre-treatment}; 
\node [block, right of=pre, node distance = 3cm] (svd) {EVD/SVD}; 
\node [block, right of=svd, node distance = 3cm] (post) {Post-treatment}; 
\draw[->, thick] (data) -- (pre) ;
\draw[->, thick] (pre) -- (svd) ;
\draw[->, thick] (svd) -- (post) ;
\end{tikzpicture}
\end{center}

\notes Multivariate Data Analysis is a classical domain in data analysis, still relevant, and underestimated (although \cite{Jolliffe2002} and a current research with ``Principal Component Analysis'' on Google Scholar provided more than 2 millions of hits). It has been developed along several lines, all over 20th century. Most of seminal papers have been published before 1935. Two trends coexist: a statistics oriented trend, developed by an Anglo-Saxon school, in UK, US, India and Scandinavia, and a French school, more algebraic and geometrical, developed in the 60's under leadership of J.-P. Benzecri. A comprehensive and very informative paper for comparing both approaches with historical insight is \cite{Tenenhaus1985}. The seminal paper on PCA by Pearson in 1901 however is geometrical, and Hotelling presentation 30 years later is algebraic. Several recent and excellent textbooks exist for a global presentation of MDA or dimensionality reduction, like \cite{Anderson1958,Mardia1979,Chatfield1980,Lee2007,Izenman2008,Wang2012}. Each has a special flavour: \cite{Anderson1958} is the seminal book on MDA (Anderson was in Stanford) and has been a bedside book of statisticians for decades; \cite{Mardia1979} is the most comprehensive, with all demonstrations of results presented as theorems, \cite{Chatfield1980} establishes a link with statistics, and is easier to read, \cite{Lee2007} goes beyond linear methods, focusing on nonlinear methods, \cite{Izenman2008} is comprehensive too, focusing on a diversity of examples, \cite{Wang2012} addressed explicitly the new challenge raised by massive data and work in high dimensional spaces. The seminal books for French school, more geometrical, are \cite{Benzecri1973a,Benzecri1973b} and \cite{Cailliez1979} who introduced the duality diagram as a unifying framework (for presentation of the duality diagram, see as well \cite{DrayDufour2007,CruzHolmes2011}). See as well \cite{Lebart1977,LMF82,LMP2000,Escofier1990} among others. The very nice paper \cite{Pages1979} gives an historical sketch of the French school as well as a comparison with Anglo-Saxon school.

%
\section*{Notations}
%

The following notations are adopted throughout these notes:\\
\\
\begin{center}
 \ovalbox{
 \begin{tabular}{cl}
 $\|.\|$ & Frobenius norm (unless otherwise stated) \\
 $\|.\|_\mathrm{sp}$ & Spectral norm \\
 $\llbracket a,b \rrbracket$ & the set of integers $i$ with $a \leq i \leq b$ \\
 $a_{ij}$ & coefficient in row $i$ and column $j$ of matrix $A$ \\
 $a_{i*}$ & row $i$ of matrix $A$ $\quad( \in \R^p)$\\
 $a_{*j}$ & column $j$ of matrix $A$ $\quad (\in \R^n)$ \\
  $A$ & a matrix in $\R^{n \times p}$ \\
  $A_r$ & a matrix of rank $r$ \\
  $A \geq 0$ & a non-negative matrix \\
  $\A$ & a point cloud of $n$ points in $\R^p$ \\ 
  $E_r$ & a subspace of $\R^p$ of dimension $r$ \\
  $\I_n$ & the identity matrix in $\R^{n \times n}$ \\
  $\L(E,F)$ & the space of linear functions from $E$ to $F$ \\
  $r$ & prescribed rank for best low rank approximation \\
  $\R^{n \times p}$ & the space of matrices with $n$ rows and $p$ columns \\  
 \end{tabular}
 }
\end{center}

\nB Notations more specific to a section are given at the beginning of the corresponding section. We have tried to be as consistent as possible between sections, but the diversity of situations and matrices involved makes such a challenge sometimes ... challenging.

%
\clearpage
\section{Principal Component Analysis (\texttt{PCA})}\label{sec:pca}
%

\subsection*{Notations}
\vspace*{5mm}

\begin{center}
\ovalbox{
\begin{tabular}{cll}
Symbol & space & meaning \\
\hline
$\alpha$ & $\N$ & index of eigenvalues of $C$\\
$a_{i*}$ & $\R^p$ & row $i$ of $A$, point $i$ in $\A$ \\
$a_{i,:}$ & $\R^p$ & row $i$ of $A$, point $i$ in $\A$ \\
$\widetilde{a_{i*}}$ & $\R^p$ & projection of $a_{i*}$ on $E_r$ \\
$A$ & $\R^{n \times p}$ & matrix to analyse \\
$\A$ & & point cloud associated with $A$ \\
$\A_r$ & &  projection of $\A$ on $E$ \\
$A_r$ & $\R^{n\times p}$ & best rank $r$ matrix for approximating $A$\\
$C$ & $\R^{p \times p}$ &  variance-covariance matrix \\
$E_r$ & $\subset \R^p$ &  best $r-$dimensional subspace \\
$\lambda$ & $\R$ &  eigenvalue of $C$ \\
$\Lambda$ & $\R^{p \times p}$ & diagonal matrix of eigenvalues of $C$ \\
$r$ & $\N$ &  rank \\
$\sigma$ & $\R$ & singular value of $A$ \\
$\Sigma$ & $\R^{p \times p}$ &  diagonal matrix of singular values of $A$ \\
$U$ & $\R^{n \times p}$ & left singular vectors of $U$, colulmnwise \\
$v$ & $\R^p$ &  eigenvector of $C$ (principal axis)\\
$V$ & $\R^{p \times p}$ & matrix of eigenvalues of $C$, columnwise
 \\
 & & right singular vectors of $A$, columnwise \\
 $Y$ & $\R^{n \times p}$ &  matrix of principal components \\
\end{tabular}
}
\end{center}

\subsection*{Matrices involved (PCA in a nutshell)}

\vspace*{5mm}

\begin{center}
\ovalbox{
\begin{tabular}{ccl}
   Matrix & dimensions & what it is \\
   \hline 
   $A$ & $n \times p$ & matrix to be analyzed \\
   $C$ & $p \times p$ & variance - covariance matrix of $A$\\
   $\Lambda$ & $p \times p$ & diagonal matrix of eigenvalues of $C$ \\
   $\Sigma$ & $p \times p$ & diagonal matrix of singular values of $A$ \\
   $U$ & $n \times p$ & left singular vectors of $A$ \\
   $V$ & $p\times p$ & principal axis (columnwise) \\
   & & right singular vectors of $A$ \\
   $Y$ & $n \times p$ & principal components \\
\end{tabular}
}
\end{center}

\noindent with
\begin{center}
\ovalbox{
\begin{tabular}{ccl}
    $C$ &=& $A^\t A$ \\
    $CV$ &=& $V\Lambda$ \\
    $A$ &=& $U\Sigma V^\t$ \\
    $\Lambda$ &=& $\Sigma^2$ \\
    $Y$ &=& $AV$ \\
    $Y$ &=& $U \Sigma$ \\
    $U^\t U$ &=& $\I_p$ \\
    $V^\t V$ &=& $\I_p$ 
\end{tabular}
}
\end{center}

\nB PCA is a Dimensionality Reduction technique which can be read in three different ways:
\begin{description}
 \item[geometric:] a point cloud $\A$ of $n$ points in $\R^p$ being given, as well as an integer $r < p$, find an affine subspace $E \subset \R^p$ of dimension $r$ such that the projection $\A_\ee$ of $\A$ on $E$ is as close as possible to $\A$.
 \item[algebraic:] a $n \times p$ matrix $A$ being given, as well as an integer $r < p$, find a matrix $A_r$ of rank $r$ such that $\|A-A_r\|$ is minimum with Frobenius $(\ell^2)$ norm. 
 \item[statistical:] a set of $p$ random variables being observed on $n$ items independently, find $r$ independent linear combination of these variables with maximum variance (the first one is with maximum variance, the second one is uncorrelated with the first one and with maximum variance, and so on). 
\end{description}
Here, we adopt the algebraic viewpoint, but start with giving some links with the geometric viewpoint, which is important for visualization of point clouds. Geometric approach is about dimension reduction, and algebraic approach is about best low rank approximation. Low rank approximation has many applications in numerical linear algebra. There are many links between algebraic and statistical approach too, which deserve to be further studied.

\notes The historical development of PCA is well known and well documented. According to \cite[chap. 3]{Basilevsky1994}, its origin can be traced in the work of Bravais in 1846 (where the notion of principal axis emerged, in ``in the form of rotating an ellipse to
`axes principaux' in order to achieve independence in a multivariate normal
distribution''). Classically, its origin is attributed to Pearson in \cite{Pearson1901} under the guise of both a statistical and a geometrical derivation, as an extension of the linear regression. The aim is statistical, but the idea behind is clearly geometric. Pearson's aim was to escape from the non symmetry of dependent and independent variables in linear regression (the regression line changes if the status dependent - independent are reversed), by giving equal status for variations on both types of variables. He was led to find the best fit of, say, a system of points in the plane by a line. One guise of his approach is statistical, in the sense that it uses the notions of mean, variance, standard deviation of a sample, but the notion of statistical model as it is understood nowadays was not available at this time. So, Pearson's approach can be qualified as geometrical. The term Factor Analysis has been introduced by Thurstone in 1931 (Thurstone, L. - 1931 - Multiple Factor Analysis, \emph{Psychological Review}, \textbf{38:}406-427). His purpose was to find a general method for finding factors which could explain correlations, following an idea published by Spearman in 1904. The presence of an underlying model in factor analysis has been the cause of numerous and fierce discussions (see \cite{Jolliffe2002}). Hotelling gave an algebraic framework in 1933 (in Hotelling, H. - 1933 - Analysis of a complex of statistical variables into principal components. \emph{J. Educ. Psychol.}, \textbf{24:}498-520), where the term PCA first appeared. Following the work of Pearson, he showed that the principal axis are the eigenvectors of the covariance matrix of the sample. PCA is based on SVD of a matrix. The link between PCA and SVD is classically attributed to a theorem published by Eckart \& Young in 1936 \cite{Eckart1936}. Classical textbooks in anglo-saxon litterature dedicated to PCA are \cite{Jackson1991}, and \cite{Jolliffe2002}. The triple nature of PCA (algebraic, geometrical, statistical) can be found in \cite{Vidal2016}, where chapter 2 is a thorough presentation of PCA. PCA is developed in every textbook in MDA, like \cite{Mardia1979,Chatfield1980,Izenman2008}. \cite{Wold1987} is a survey of PCA tools with an introduction on main milestones in the development of the method. Classical textbooks in French literature are \cite{Cailliez1979,LMF82}. A recent survey of \PCA and its recent developments can be found in \cite{Jolliffe2016}.\\
\\
PCA is probably one of the most used tool in multivariate statistics.  It is known under various names according to the field it is applied to, as mentioned  in \url{https://en.wikipedia.org/wiki/Principal_component_analysis}: Principal Component Analysis, Karhunen–Loève transform (KLT) in signal processing,  proper orthogonal decomposition (POD) in mechanical engineering, empirical orthogonal functions (EOF) in meteorological science, among others.\\
\\
If the theory can be derived for any pair $(p,r)$ of dimensions with $1 \leq r < p$, it is most interesting when $p$ is large and $r \ll p$.

%
\subsection{Setting the problem}\label{sec:pca:set}
%

Let us recall that the norm here and in the following sections is the Frobenius or $\ell^2$ norm unless otherwise stated:
\begin{equation}
 \|x\| = \left(\sum_ix_i^2\right)^{1/2}
\end{equation}
for a vector and 
\begin{equation}
 \|A\| = \left(\sum_{i,j}a_{ij}^2\right)^{1/2}
\end{equation}
for a matrix.

\nB The problem can be set as follows:

\nT{Algebraic approach}\\
\begin{center}
 \shadowbox{
  \begin{tabular}{ll}
      Given & $A \in \R^{n \times p}$ \\
      & $0 < r < p$ \\
      \\
      Find & $A_r \in \R^{n \times p}$\\ 
      with & $\rank A_r=r$ \\
      \\
      such that & $ \|A-A_r\|\quad \mbox{minimal}$ \\
  \end{tabular}
 }
\end{center}

\nT{Geometric approach}\\
\\
In this approach, a cloud denoted $\A$ of $n$ points in $\R^p$, labeled $(a_i)_i$ with $1 \leq i \leq n$, is associated to a matrix $A \in \R^{n \times p}$, where point $a_i$ is row $i$ of $A$: $a_i = a_{i,:}$. The point cloud associated to $A_r$ is denoted $\A_r$\\
\begin{center}
 \shadowbox{
  \begin{tabular}{ll}
      Given a point cloud & $\A = (a_1, \ldots, a_n)$ \\
      & $a_i \in \R^p$ \\
      & $0 < r < p$ \\
      \\
      Find & a subspace $E_r \subset \R^p$\\ 
      with & $\dim E_r=r$ \\
      \\
      such that & $ d(\A, \A_r) \quad \mbox{minimal}$ \\
      \\
      where & $d(\A,\A_r) = \sum_i\|a_i - \widetilde{a_i}\|^2$ \\
      and & $\widetilde{a_i}$ is the projection of $a_i$ on $E_r$
  \end{tabular}
 }
\end{center}

\nS Let us note that a point cloud $\A$ is equivalent to a matrix $A$ with the row $i$ of $A$ being the point $a_i \in \A$. Knowing that, both approaches are equivalent. To see that, let us consider an orthonormal basis $(v_1, \ldots, v_r)$ of $E_r$, and let us complete it to have an orthonormal basis of $\R^p$ as $(v_1, \ldots, v_r, v_{r+1}, \ldots, v_p)$.  If the projection is exact, i.e. if $\forall \: i, \: a_i=\widetilde{a_i}$, the columns $r+1$ to $p$ of $A$ in basis $V$ are zero, and $A$ is of rank $r$. The converse is true as well.

%
\subsection{Solving the problem}\label{sec:pca:solve}
%

The solution to this problem is well known (see any textbook mentioned in the introduction of this section). Finding the subspace $E_r$ is finding an orthonormal basis for it. Let us fix a subspace $E_r \subset \R^p$ of dimension $r$. If $\widetilde{a_i}$ is the projection of $a_i$ on $E_r$, we have, by Pythagore theorem
\[
 \forall \: i \in \llbracket 1,n \rrbracket, \quad \|a_i\|^2 = \|\widetilde{a_i}\|^2 + \|a_i - \widetilde{a_i}\|^2
\]
Then, setting $\sum_i \|a_i - \widetilde{a_i}\|^2$ minimum is equivalent to setting $\sum_i\|\widetilde{a_i}\|^2$ maximum. We then can set PCA as\\
\\
\begin{equation}
 \left|
  \begin{tabular}{ll}
      Given a point cloud & $\A = (a_1, \ldots, a_n)$ \\
      & $a_i \in \R^p$ \\
      & $0 < r < p$ \\
      \\
      Find & a subspace $E_r \subset \R^p$\\ 
      with & $\dim E_r=r$ \\
      \\
      such that & $\sum_i\|\widetilde{a_i}\|^2 \quad \mbox{maximal}$ \\
      \\
      where & $\widetilde{a_i}$ is the projection of $a_i$ on $E_r$
  \end{tabular}
 \right. 
\end{equation}

\vspace*{5mm}\noindent Let us consider the simple case $r=1$ and $0 \in E_r$. If $u$ with $\|u\|=1$ is a basis of $E_1$, we have $\widetilde{a_i}= \langle a_i,u\rangle u$, and the optimization problem can be stated as\\
\\
\begin{equation}
 \left|
  \begin{tabular}{ll}
      Given a point cloud & $\A = (a_1, \ldots, a_n)$ \\
      & $a_i \in \R^p$ \\
      \\
      Find & a vector $u \in \R^p$\\ 
      with & $\|u\|=1$ \\
      \\
      such that & $\|Au\| \quad \mbox{maximal}$ \\
  \end{tabular}
 \right. 
\end{equation}
Indeed, $\widetilde{a_i}= \langle a_i,u\rangle \, u$. Then $\sum_i\|\widetilde{a_i}\|^2 = \sum_i\langle a_i,u\rangle^2$. The elements $\langle a_i,u\rangle$ are the coordinates of the vector $y=Au$. The solution of such a problem is classical: $u$ is the eigenvector of $A^\t A$ associated with the largest eigenvalue
\begin{equation}
 A^\t Av = \lambda v, \qquad \lambda = \max \{ \lambda \in \Sp A^\t A\}
\end{equation}
The Eckart-Young theorem \cite{Eckart1936} extends this result to $r > 1$ and states that an orthonormal basis of $E_r$ is the set $(v_1, \ldots,v_r)$ with
\begin{equation}
 A^\t Av_j = \lambda_jv_j\qquad \mbox{with} \quad \lambda_1 \geq \ldots \geq \lambda_r \geq \lambda_{r+1} \geq \ldots \geq \lambda_p \geq 0
\end{equation}
This is a direct application of the variational properties of the Rayleigh quotients (see \cite[sec. 4.2]{Horn2012}). This can be written
\begin{equation}
 A^\t AV = V\Lambda
\end{equation}
where $V$ is the $p \times p$ matrix with column $j$ being $v_j$ and $\Lambda$ is the diagonal matrix with terms $(\lambda_i)_i$ in the diagonal (in decreasing order). The coordinates of the point cloud $\A$ in new basis $V$ are given by
\begin{equation}
 Y = AV
\end{equation}
Hence a first algorithm, for a function called \textsc{pca\_evd}$()$:\\
\\
\begin{algorithm}[H]
\begin{algorithmic}[1]
\STATE \textbf{input} $A \in \R^{n \times p}$
\STATE \textbf{compute} $C=A^\t A$
\STATE \textbf{compute} $(\lambda_\alpha, v_\alpha)$ such that $Cv_\alpha = \lambda_\alpha v_\alpha$, or $CV=V\Lambda$
\STATE \textbf{compute} $Y=AV$
\RETURN $Y,\Lambda, V$
\end{algorithmic}
\caption{PCA of a matrix with EVD: \textsc{pca\_evd}$(A)$}
\label{alg:pca}
\end{algorithm}

\nS The vectors in new basis are called \emph{principal axis}, and the coordinates along the principal axis are called \emph{principal components}.

\nS Here is a summary of the results:\\
\\
\begin{center}
 \ovalbox{
 \begin{tabular}{c|ll|l}
  $C$ & $\in \R^{p \times p}$ & matrix of correlations of columns of $A$ & $C = A^\t A$ \\
  $V$ & $\in \R^{p \times p}$ & new basis & $Cv = \lambda v$ \\ 
  $\Lambda$ & $\in \R^p$ & eigenvalues of $C$ in decreasing order& \\
$Y$ & $\in \R^{n \times p}$ & matrix of coordinates in new basis & $Y=AV$\\
 \end{tabular}
 }
\end{center}

%
\subsection{Link with SVD}\label{sec:pca:svd}
%

Let $(U,\Sigma, V)$ be the SVD of $A$.
\begin{equation}
 A=U\Sigma V^\t
\end{equation}
Then
\begin{equation}
 \begin{aligned}
   C &= A^\t A \\
   &= (V \Sigma U^\t)(U\Sigma V^\t) \\
   &= V \Sigma^2V^\t \qquad \mbox{as}\quad U^\t U=\I_n
 \end{aligned}
\end{equation}
and
\begin{equation}
 CV = V\Sigma^2 \qquad \mbox{as} \quad V^\t V = \I_p
\end{equation}
Hence, $V$ as new basis for PCA of $A$ is the matrix $V$ in SVD of $A$, and $\Lambda=\Sigma^2$. We have
\begin{equation}
 \begin{aligned}
   Y &= AV \\
   &= U\Sigma V^\t V \\
   &= U\Sigma
 \end{aligned}
\end{equation}
This yields a second algorithm for PCA:\\
\\
\begin{algorithm}[H]
\begin{algorithmic}[1]
\STATE \textbf{input} $A \in \R^{n \times p}$
\STATE \textbf{compute} $U, \Sigma, V = \textsc{svd}(A)$
\STATE \textbf{compute} $\Lambda=\Sigma^2$
\STATE \textbf{compute} $Y=U\Sigma$
\RETURN $Y,\Lambda,V$
\end{algorithmic}
\caption{PCA of a matrix with SVD: \textsc{pca\_svd}$(A)$)}
\end{algorithm}

\paragraph{About the choice between SVD and EVD:} Adopting the SVD viewpoint has some numerical advantages:
\begin{itemize}[label=$\rightarrow$]
 \item there exists efficient and stable algorithms for computing an SVD, 
 \item it avoids a matrix $\times$ matrix computation $(C=A^\t A)$ which can be costly when $n$ and $p$ are large,
 \item it leads to easier generalization with instrumental variables or metrics on row or column space (see section \ref{alg:pcamet}),
 \item If the dimensions $n$ and $p$ are (very) large, the SVD can be computed with random projection (see \cite{Halko2011}). Such a calculation is presented below.
\end{itemize}
However, matrix $C=A^\t A$ is the variance covariance matrix of the distribution of the variables in the statistical approach, and must not be ignored.

%
\subsection{Randomized SVD}
%

Let $A \in \R^{n \times p}$ with $n \geq p$. The complexity (number of operations) of the SVD of $A$ is in $\mathcal{O}(n^2p)$. SVD becomes untractable for large values of $n$ and $p$, say $10^4$. Fortunately, there are some very efficient heuristics, with bounds on errors, to compute the first singular values and vectors, based on randomized algorithms.  The idea behind is the following. 

\nB If $Q \in \R^{n \times k}$ is columnwise orthonormal, the projection of the columns of $A$ on the vector space spanned by the columns of $Q$ is
\[
 \widetilde{A}= QQ^\t A
\]
Let us denote
\[
 B = Q^\t A, \qquad B \in \R^{k \times p}
\]
Then, $\widetilde{A}= QB$ is a rank $k$ approximation of $A$. The SVD of $B$ $(B=U_\b\Sigma V)$ is in $\mathcal{O}(p^2k)$ instead of $\mathcal{O}(n^2p)$, and we have 
\begin{equation*}
 A \approx QB = Q(U_\b\Sigma V)= (QU_\b)\Sigma V=U\Sigma V
\end{equation*}
So 
\begin{equation}
 A \approx U\Sigma V \qquad \mbox{with} \quad U= QU_\b 
\end{equation}
Next step is to show that (what we will not develop here, see appendix \ref{sec:random_projection}), when $n$ and $p$ are large, $A \approx QQ^\t A$ with high quality for any random matrix $Q$. This comes from deep theorems in geometry of Banach spaces, and from Johnson-Lindenstrauss lemma which states that, for any $\epsilon>0$, and any dimension $n$, there exists a dimension $k$ such that for any cloud $X$ of $n$ points, there exists an embedding
\[
 \begin{CD}
  f \: : \: \R^n @>>> \R^k
 \end{CD}
\]
such that for any $x,y \in X$
\begin{equation}
 (1-\epsilon)\|x-y\|^2 \leq \|f(x)-f(y)\|^2 \leq (1+\epsilon)\|x-y\|^2
\end{equation}
A demonstration of the existence relies on showing that such an embedding exists with probability one. The dimension $k$ must comply with
\begin{equation}
 k \geq \frac{8\; \Log n}{\epsilon^2}
\end{equation}
The bad news is that $\epsilon^2$ is at the denominator (so, $k$ is large when $\epsilon$ is small), but the good news is that $k$ grows with $\Log n$ and not $n$. This becomes efficient when $n$ is large. So, $Q$ as an orthonormal basis is built as the $QR-$ decomposition of $Y=A\Omega$ where $\Omega$ is a random matrix ($Y$ is in the span of $A$). There are several ways to chose $\Omega$, and here we restrict ourselves to the Gaussian random projection, i.e. $\Omega$ is a random Gaussian matrix with
\[
 \Omega[i,j] \sim \mathcal{N}(0,1)
\]
Usually, for a good accuracy at rank $k$, it is advised to select $\Omega$ as $n \times k'$ with $k'=k+s$ where $s$ is called the oversampling. Usually, taking $s=5$ is said to be sufficient. The reader is encouraged to read \cite{Halko2011} for further details and explanations on randomized algorithms in matrix computations (what is presented here is the tip of the iceberg). The algorithm runs as follows:\\
\begin{algorithm}[H]
\begin{algorithmic}[1]
\STATE \textbf{input} $A \in \R^{n \times p}$, $k$ as prescribed rank
\STATE \textbf{build} $\Omega \in \R^{p \times k}$, random $(\Omega[i,j] \sim \mathcal{N}(0,1))$
\STATE \textbf{compute} $Y=A\Omega$
\STATE \textbf{compute} the $QR-$decomposition of $Y$: $Y=QR$ 
\STATE \textbf{build} $B=Q^\t A$
\STATE \textbf{run} the SVD of $B$: $B=U_\b\Sigma V^\t$, $\quad$ or $(U_\b,\Sigma,V) = \textsc{svd}(B)$ 
\STATE \textbf{compute} $U=QU_\b$ 
\RETURN $U, \Sigma,V$
\end{algorithmic}
\caption{SVD of a matrix with Gaussian Random Projection \textsc{svd\_grp}$(A,k)$)}
\label{alg:svd_grp}
\end{algorithm}

\nS Here are the dimensions of the involved matrices:\\
\begin{center}
 \begin{tabular}{c|c|l}
   \hline  
   Matrix & dimensions & computation \\
   \hline 
   $A$ & $n \times p$ & data \\
   $\Omega$ & $p \times k$ & Gaussian random matrix\\
   $Y$ & $n \times k$ & $Y=A\Omega$\\
   $Q$ & $n \times k$ & $Y=QR$\\
   $B$ & $k \times p$ & $B=Q^\t A$\\
   $U_\b$ & $k \times k$ & $B=U_\b \Sigma V^\t$\\
   $\Sigma$ & $k \times k$ & idem \\
   $V$ & $k \times p$ & idem \\
   $U$ & $n \times k$ & $U=Q U_\b$\\
   \hline
 \end{tabular}
\end{center}

\nS In section \ref{sec:mds}, we will see that MDS relies on the SVD or EVD of the Gram matrix of a point cloud (if $\X = (x_i)_i$ is a point cloud with $x_i \in \R^p$ and $1 \leq i \leq n$, the Gram matrix $G$ of $\X$ is the symmetric definite positive matrix of elements $g_{ij} = \langle x_i,x_j \rangle$). Then, $G=XX^\t$. MDS is solving the reverse problem: finding $X$ knowing $G$. If $G=U\Sigma U^\t$ is the SVD of $G$, then a solution to MDS is $X=U\Sigma^{1/2}$. The SVD  can be done by Random Projection by selecting $\Omega$, computing $Y=G\Omega$, $Q$ such that $Y=QR$, defining $G'=QQ^\t G \approx G$ and doing the SVD of $G'$. If $B=Q^\t G$, the SVD of $B$ is $B=U_\b \Sigma V^\t$ and $G'=QB = U\Sigma V^\t$ with $U=QU_\b$. However, one cannot write $G'=X'X^{'\t}$ because $G'$ is not symmetric, and $V \neq U$. To solve this, one can do a double projection on $G$, rowwise and columnwise, and define $G''=(QQ^\t)G(QQ^\t)$, with still $G'' \approx G$. One defines $C=Q^\t G Q$, which is symmetric, do its SVD by $C=U_\c \Sigma U_\c^\t$, hence $G'' = U\Sigma U^\t$ with $U=QU_\c$. One then can write $G \approx X''X^{''\t}$ with $X'' = U\Sigma^{1/2}$. This can be derived as follows:\\
\begin{algorithm}[H]
\begin{algorithmic}[1]
\STATE \textbf{input} $G \in \R^{n \times n}$, $k$ as prescribed rank
\STATE \textbf{build} $\Omega \in \R^{n \times k}$, random $(\Omega[i,j] \sim \mathcal{N}(0,1))$
\STATE \textbf{compute} $Y=G\Omega$
\STATE \textbf{compute} the $QR-$decomposition of $Y$: $Y=QR$ 
\STATE \textbf{build} $C=Q^\t C Q$
\STATE \textbf{run} the SVD of $C$: $C=U_\c\Sigma U_\c^\t$, $\quad$ or $(U_\c,\Sigma,U_\c) = \textsc{svd}(C)$ 
\STATE \textbf{compute} $U=QU_\c$ 
\RETURN $U, \Sigma$
\end{algorithmic}
\caption{SVD of a Gram matrix with Gaussian Random Projection \textsc{svd\_grp\_gram}$(G,k)$)}
\label{alg:svd_grp_gram}
\end{algorithm}

\nS Here are the dimensions of the involved matrices:\\
\begin{center}
 \begin{tabular}{c|c|l}
   \hline  
   Matrix & dimensions & computation \\
   \hline 
   $G$ & $n \times n$ & Gram matrix \\
   $\Omega$ & $n \times k$ & Gaussian random matrix\\
   $Y$ & $n \times k$ & $Y=G\Omega$\\
   $Q$ & $n \times k$ & $Y=QR$\\
   $C$ & $k \times k$ & $B=Q^\t G Q$\\
   $U_\c$ & $k \times k$ & $C=U_\c \Sigma U_\c^\t$\\
   $\Sigma$ & $k \times k$ & idem \\
   $U$ & $n \times k$ & $U=Q U_\c$\\
   \hline
 \end{tabular}
\end{center}

\noindent One observes that the complexity of the calculation of the SVD of $C$ is in $\OC(k^3)$, whereas the calculation of $Y=G\Omega$ is in $\OC(n^2k)$, and more ``expensive''.

\notes The property behind Random Projection (RP) is counter-intuitive: ``RP, while reducing dimensionality, approximatively preserves pairwise distances with high probability'' (\cite[p. 2]{Vempala2004}). This is a consequence of Johnson-Lindenstrauss lemma (\cite{Johnson1984}, see appendix \ref{sec:random_projection} for details). \cite{Vempala2004} presents a variety of domains of application of random projection. This is highly counter-intuitive, because it can undermine the very notion of PCA. Indeed, a point cloud and a rank being given, PCA is about finding the best subspace of dimension $k$ as far as quality of projection is concerned. RP answer is that any randomly selected subspace of dimension $k$ will make the job! Actually, this is not true for very small $k$, like first two or three axis useful for visualizing the shape of the point cloud. However, when $k$ is significantly larger, this is true. It is possible to show (see appendix \ref{sec:random_projection}) that the best choice is that $k$ is in $\OC(\frac{\Log n}{\epsilon^2})$ where $\epsilon$ is the relative accuracy in distances preservation. The presence of $\epsilon^2$ at the denominator forces $n$ to be very large for $\Log n < n \epsilon^2$. This is however true for any point cloud, including random ones. In applications, point clouds have a pattern or a structure, and RP works for much smaller values of $n$. See for example \cite{Bingham2001} for a comparison in accuracy and computing time with some other dimensionality reduction tools like PCA and an overview of various methods for selecting the matrix $\Omega$ (Gaussian matrix is not the only possible choice). The use of randomized algorithms to compute efficiently the SVD of a very large matrix is fully developed in \cite{Halko2011}. This section explains how SVD with Gaussian Random Projection works. To understand why it works, see appendix \ref{sec:random_projection}.

%
\subsection{Core algorithm for PCA}
%

\nS Wrapping all this together leads to a core algorithm for PCA, where the user can select which method to implement, presented hereafter in pseudocode:\\
\\
\begin{algorithm}[H]
\begin{algorithmic}[1]
\STATE \textbf{input} $A \in \R^{n \times p}$, $k \in \N^* \cup \{-1\}$,  $\texttt{meth} \in \{\textsc{evd}, \textsc{svd}, \textsc{grp}\}$
\IF{\texttt{meth}==\textsc{evd}}
\STATE \textbf{compute} $C=A^\t A$
\STATE \textbf{compute} $(\lambda_\alpha, v_\alpha)$ such that $Cv_\alpha = \lambda_\alpha v_\alpha$, or $CV=V\Lambda$
\STATE \textbf{compute} $Y=AV$
\IF{$k>0$}
\STATE $Y=Y\texttt{[:,0:k]}; \quad \Lambda=\Lambda\texttt{[0:k]}; \quad  V=V\texttt{[:,0:k]}$
\ENDIF
\ENDIF
\IF{\texttt{meth}==\textsc{svd}}
\STATE \textbf{compute} $U, \Sigma, V = \textsc{svd}(A)$
\STATE \textbf{compute} $\Lambda=\Sigma^2$
\STATE \textbf{compute} $Y=U\Sigma$
\IF{$k>0$}
\STATE $Y=Y\texttt{[:,0:k]}; \quad \Lambda=\Lambda\texttt{[0:k]}; \quad  V=V\texttt{[:,0:k]}$
\ENDIF
\ENDIF
\IF{\texttt{meth}==\textsc{grp}}
\STATE \textbf{compute} $U, \Sigma, V = \textsc{svd\_grp}(A,k)$
\STATE \textbf{compute} $\Lambda=\Sigma^2$
\STATE \textbf{compute} $Y=U\Sigma$
\ENDIF
\RETURN $Y,\Lambda,V$
\end{algorithmic}
\caption{PCA of a matrix: \textsc{pca\_core}$(A, k=-1,\texttt{meth=}\textsc{svd})$}
\label{alg:pca:core}
\end{algorithm}

\nT{Comments} Here are some comments:
\begin{description}
 \item[Why \textsc{pca\_core}?] The reason for the name is the following: PCA is a method which very seldom runs dimension reduction or approximation by a low rank matrix directly on the data matrix. Most of the times, there is a preatreatment (see section \ref{sec:pca:classical}), like centering and scaling columnwise, and dimension reduction is performed after pre-treatment. The name PCA designates classically the whole analysis: 
 \begin{enumerate}
  \item a pretreatment
 \[
  \begin{CD}
   A @>\mathrm{pretreatment}>> A'
  \end{CD}
  \]
  \item the treatment itself
  \[
    Y,\Lambda,V=\textsc{pca\_core(A')}
  \]
  which is denoted $\textsc{pca\_core}()$.
 \end{enumerate}
 \item[Choice of the method:] Three methods are described here: eigenvalues of the correlation matrix (\textsc{evd}), SVD of the data matrix (\textsc{svd}), and SVD with Gaussian Random Projection (\textsc{grp}). Here is an advice for selecting the right method 
 \begin{itemize}[label=$\rightarrow$]
  \item If the size of the matrix (number $p$ of columns) is small to medium (say $\approx 10^3$), then \textsc{evd} or \textsc{svd} can be used
  \item If the size of the matrix is large to very large $(10^4 < p < 10^6)$, gaussian random projection must be used.
 \end{itemize}
 \item[Prescribed rank:] $k$ is the prescribed rank. Its default value is $k=-1$, which means that all the eigenvalues, or singular values, and components and axis will be computed. This is relevant for methods \textsc{evd} or \textsc{svd} only. A rank $k >0$ must be prescribed for method \textsc{grp}. If a rank $k > 0$ is prescribed, the first $k$ eigenvalues or singular values, components and axis only will be computed.

\end{description}

%
\subsection{Interpretation and plotting}\label{sec:pca:plot}
%

In this section, the geometrical viewpoint is adopted. Interpretation of the results of a \PCA is about quantifying the fraction of inertia of the point cloud (i.e. variance of the associated variables, or norm of the associated matrix) which is preserved by projection either on one axis or on a space spanned by $r$ first axis. When the point cloud is centered,  \PCA is finding a rotation in $\R^p$ such that these quantities are maximal.

\nS Let $A \in \R^{n \times p}$ be a matrix, $\A$ its associated point cloud in $\R^p$, and $(Y,\Lambda,V)$ the \PCA of $A$. Then, the column $y_j$ of $Y$ with $j \in \llbracket 1,p\rrbracket$ is the vector in $\R^n$ of the coordinates of the points of $\A$ on principal axis $j$.

\begin{proof}
Indeed, let us have a point cloud $\A$ made of $n$ points $a_i \in \R^p$ with $a_i$ being row $i$ of $A$. \PCA of $A$ is performing an SVD of $A$ as
\begin{equation}
 A = U\Sigma V^\t
\end{equation}
and yields a new orthogonal basis $(v_1, \ldots,v_p)$ of $\R^p$, where $v_j$ is the column $j$ of $V$. Let us recall that
\begin{equation}
 Y = U\Sigma
\end{equation}
Then
\begin{equation}
 Y = U \Sigma = U \Sigma (V^\t V) = (U \Sigma V^\t) V =  AV
\end{equation}
which means that the rows of $Y$ are the coordinates of the points of $\A$ in new basis $V$.
\end{proof}

\nS If $y \in \R^n$ and $v \in \R^p$, let us recall the notation $\otimes$ for tensor product
\[
 y \otimes v \equiv yv^\t 
\]
(i.e. $(y \otimes x)_{ij}=y_iv_j$). Let us denote by $y_j$ the column $j$ of $Y$, and $v_j$ the column $j$ of $V$. Then, 
\begin{equation}
 A = \sum_{j=1}^p y_j \otimes v_j
\end{equation}
from which
\begin{equation}
 \|A\|^2 = \sum_j \|y_j\|^2
\end{equation}
\begin{proof}
Indeed,
\[
 \begin{aligned}
  \|A\|^2 &= \left\|\sum_{j=1}^p y_j \otimes v_j\right\|^2 \\
  &= \left\langle \sum_{j=1}^p y_j \otimes v_j \; , \; \sum_{k=1}^p y_k \otimes v_k\right\rangle  \\
  &= \sum_{j,k} \langle y_j \otimes v_j \; , \; y_k \otimes v_k \rangle \\
  &= \sum_{j,k} \langle y_j \, , \, y_k \rangle. \; \langle v_j \, , \, v_k \rangle   \\
  &= \sum_j \langle y_j \, , \, y_j \rangle \\
  &= \sum_j\|y_j\|^2
 \end{aligned}
\]
\end{proof}
\noindent (another way to see this is to observe that $Y$ is deduced from $A$ by a rotation, which is an isometry, then $\|A\|=\|Y\|$). As $Y=U\Sigma$ with $U^\t U=\I_p$, we have
\begin{equation}
 \|y_j\| = \sigma_j
\end{equation}
Hence, the norm of $A$ can be partitioned as 
\begin{equation}
 \|A\|^2 = \sum_j \sigma_j^2
\end{equation}

\nB Let us recall that
\[
 \lambda_j = \sigma_j^2, \qquad A^\t Av_j = \lambda_jv_j
\]
Then
\begin{equation}
 \|A\|^2 = \sum_{i=1}^p\lambda_i
\end{equation}
and the quality of the representation of $A$ by its projection on the axis spanned by $v_j$ is
\begin{equation}
 \varrho_j = \frac{\lambda_j}{\sum_i\lambda_i}
\end{equation}
The quality of representation of the point cloud (i.e. of array $A$) by its projection $A_r$ on the subspace spanned by vectors $(v_1,\ldots,v_r)$ is
\begin{equation}
  \begin{aligned}
    \rho_r &= \sum_{j=1}^r \varrho_j \\
    &= \frac{\sum_{j=1}^r\lambda_j}{\sum_i\lambda_i}
  \end{aligned}
\end{equation}

\nS The quality of representation of item $i$ on axis $ j \in \{1,p\}$ is
\begin{equation}
 \psi(i,j) = \frac{y_{ij}^2}{\sum_{\ell=1}^p y_{i\ell}^2}
\end{equation}
and the quality of projection of item $i$ on the subspace spanned by vectors $(v_1,\ldots,v_r)$ is
\begin{equation}
 \begin{aligned}
    \theta(i,r) &= \sum_{j=1}^r\psi(i,j) \\
    &= \frac{\sum_{j=1}^ry_{ij}^2}{\sum_{\ell=1}^p y_{i\ell}^2}
 \end{aligned}
\end{equation}

\nB We have 
\begin{equation}
 \left\{
  \begin{aligned}
    \varrho_j &= \sum_{i=1}^n\psi(i,j) \\
    \rho_r &= \sum_{i=1}^n\theta(i,r)
  \end{aligned}
 \right.
\end{equation}

\nB This can be summarized as\\
\\
\begin{center}
 \ovalbox{
  \begin{tabular}{lccl}
    Quality of representation of & Notation & Calculation & Observation\\
    \hline
    item $i$ on axis $j$ & $\psi(i,j)$ & $\displaystyle \frac{y_{ij}^2}{\sum_{\ell=1}^p y_{i\ell}^2}$ & \\  
    &&&\\
    item $i$ on subspace $E_r$ & $\theta(i,r)$ & $\sum_{j=1}^r \psi(i,j)$ & \\
    &&&\\
    point cloud on axis $j$ & $\varrho_j$ & $\displaystyle \frac{\lambda_j}{\sum_{i=1}^p\lambda_i}$ & $= \sum_{i=1}^n\psi(i,j)$\\
    &&&\\
    point cloud on subspace $E_r$ & $\rho_r$ & $\displaystyle \sum_{j=1}^r \varrho_j$ & $= \sum_{i=1}^n\theta(i,r)$
  \end{tabular}
 }
\end{center}

\nT{Prescribed rank or accuracy} In the algebraic framework, \PCA is about the best low rank approximation of a matrix. This can be set in two guises:
\begin{itemize}[label=$\rightarrow$]
 \item select a rank $r$, and deduce the quality $\rho$ of the approximation, 
 \item select a quality of an approximation, and deduce the rank at which it should be done.
\end{itemize}
The former is \PCA at \emph{prescribed rank}, whereas the latter is \PCA at \emph{prescribed accuracy}. The key tool for implementing the one or the other is the curve $r \mapsto \rho(r)$.

%
\subsection{Classical analysis}\label{sec:pca:classical}
%

Here, we denote $A \geq 0$ if all coefficients in $A$ are nonnegative. Let us recall that $a_i \in \R^p$, $a_{i*}$ denotes the row $i$ of $A \in \R^{n \times p}$, and $a_{*j} \in \R^n$ its column $j$.

\nB  The mean and standard deviation of a distribution are often the best summary of it. If \PCA yields the best rank one approximation, it is likely that first axis and components mirror this best summary and bring few information on the inner structure of matrix $A$. Hence a standard procedure is to center and scale a dataset, and run the \PCA on the scaled and centered dataset to focus on the inner structure (e.g. correlations between columns). 

\nT{Centering} Centering a matrix $A$ is translating the attached point cloud $\A$ to its barycenter:
\begin{equation}
 \begin{CD}
   \mbox{in}\: \R^p, \qquad a_i @>\mathrm{centering}>> \overline{a_i} = a_i-g, 
 \end{CD}
\end{equation}
where $g \in \R^p$ is the barycenter of the point cloud, i. e. 
\begin{equation}
 g_j = \frac{1}{n}\sum_i \, a_{ij}
\end{equation}
It is easy to check that $\sum_i\overline{a_i}= \sum_i\left(a_i-\frac{1}{n}\sum_ia_i\right)=\sum_ia_i-\sum_ia_i=0$. Matrix $\overline{A}$ is centered columnwise:
\begin{equation}
 \forall \: j, \quad \sum_i \;  \overline{a}_{ij}=0
\end{equation}

\nT{Scaling} Scaling a matrix columnwise is dividing each column vector $a_{*j}$ by its norm, aka its standard deviation if it is centered:
\begin{equation}
 \begin{CD}
  a_{*j} @>\mathrm{scaling}>> \displaystyle \frac{a_{*j}}{\|a_{*j}\|}
 \end{CD}
\end{equation}
The centered-scaled matrix $A'$ is defined by
\begin{equation}
 \begin{CD}
   a_{*j} @>>> \displaystyle \frac{\overline{a}_{*j}}{\|\overline{a}_{*j}\|} \qquad \mbox{with} \quad \overline{a}_{*j} = a_{*j}-g_j\ones_n
 \end{CD}
\end{equation}

\nS Scaled-centered \PCA of a matrix $A$ is defined as:\\
\begin{algorithm}[H]
\begin{algorithmic}[1]
\STATE \textbf{input} $A \in \R^{n \times p}$
\STATE \textbf{compute} the barycenter of $A$: $g = \frac{1}{n}\sum_ia_{i*}$
\STATE \textbf{center} $A$ : $A \longrightarrow \overline{A}$, with $\forall \: i, \quad a_i \longrightarrow \overline{a}_i = a_i-g$
\STATE \textbf{scale} $\overline{A}$: $\overline{A} \longrightarrow A'$, with $\forall \: j, \quad \overline{a}_{*j} \longrightarrow a'_{*j} = \frac{\overline{a}_{*j}}{\|\overline{a}_{*j}\|}$ 
\STATE \textbf{do} $Y,\Lambda,V =\textsc{pca\_core}(A')$
\RETURN $Y, \Lambda, V$
\end{algorithmic}
\caption{$\textsc{pca-sc}(A)$}
\label{alg:scpca}
\end{algorithm}

\nS There are a few elementary results for scaled-centered PCA. By definition, the coefficients $c_{j\ell}$ of $C=A^{'\t}A'$ are the correlations between centered scaled variables $a'_{*j}$ and $a'_{*\ell}$. Hence, we have
\begin{equation}
 -1 \leq c_{j\ell} \leq 1
\end{equation}
We have as well
\begin{equation}
 \forall \: j, \qquad c_{jj} = 1
\end{equation}
Hence
\begin{equation}
 \sum_j \lambda_j = \Tr C = p
\end{equation}
Hence, the quality of approximation at rank $r$ of $A'$ is
\begin{equation}
 \rho_r = \frac{1}{p}\sum_{j \leq r}\, \lambda_j
\end{equation}

\nT{Double averaging or bicentering} Let us have a matrix $A \geq 0$ which is for example an array of counts. A classical example is a contingency table (contingency tables can be analysed with Correspondecne Analysis, see section \ref{sec:coa}). The structure of $A$ is dominated by the property $A \geq 0$. If a \PCA of $A$ is run, this will be the main (trivial) information given by axis 1. This trivial information can be filtered out by setting the model \begin{equation}
 a_{ij} = \underbrace{m}_{\mathrm{global}\:\mathrm{mean}} +   \underbrace{x_i}_{\mathrm{effect}\:\mathrm{of}\:\mathrm{row}\: i} + \underbrace{y_j}_{\mathrm{effect}\:\mathrm{of}\:\mathrm{column}\: j}  + \underbrace{r_{ij}}_{residuals}
\end{equation}
with
\begin{equation}
 \left\{
   \begin{array}{rcl}
    \displaystyle \sum_i x_i &=& 0 \\
    \displaystyle  \sum_j y_j &=& 0 \\
    \displaystyle \forall \: j, \quad \sum_ir_{ij} &=& 0 \\
    \displaystyle  \forall \: i, \quad \sum_jr_{ij} &=& 0 
   \end{array}
 \right.
\end{equation}
Then, we have
\begin{equation}
 \left\{ 
    \begin{array}{lcl}
      m &=& \displaystyle \frac{1}{np}\sum_{i,j}a_{ij} \\
      x_i &=& \displaystyle \left(\frac{1}{p}\sum_ja_{ij}\right)-m \\
      y_j &=& \displaystyle \left(\frac{1}{n}\sum_ia_{ij}\right)-m
    \end{array}
 \right.
\end{equation}
\begin{proof}
We have
\begin{equation*}
 \sum_{i,j}a_{ij} = np\;m
\end{equation*}
so 
\begin{equation*}
m = \frac{1}{np}\sum_{ij}a_{ij} 
\end{equation*}
Then
\begin{equation*}
   \sum_ia_{ij} = n\, m + ny_j
\end{equation*}
and 
\begin{equation}
 y_j = \left(\frac{1}{n}\sum_ia_{ij}\right) - m
\end{equation}
Similarly
\begin{equation*}
 \sum_ja_{ij} = p\, m + px_i
\end{equation*}
and 
\begin{equation}
 x_i = \left(\frac{1}{p}\sum_ja_{ij}\right) -m
\end{equation}
\end{proof}
\noindent So, we have
\begin{equation}
 \begin{array}{lcl}
  a_{ij} &=& m + x_i + y_i + r_{ij} \\
  &=& \displaystyle m + \left(\frac{1}{p}\sum_ja_{ij}\right) -m + \left(\frac{1}{n}\sum_ia_{ij}\right) - m + r_{ij} \\
  &=& \displaystyle - \frac{1}{np}\sum_{i,j}a_{ij} + \left(\frac{1}{p}\sum_ja_{ij}\right) + \left(\frac{1}{n}\sum_ia_{ij}\right) + r_{ij}
 \end{array}
\end{equation}
which is denoted as well
\begin{equation}
 r_{ij} = a_{ij} - a_{i.} - a_{.j} + a_{..} 
\end{equation}

\nS \PCA with double averaging is
\begin{enumerate}
 \item computing the global mean $m$, each effect $x_i$ and $y_j$ and the matrix of residuals $R$
 \item run the \PCA of $R$, which is already centered, without scaling.
\end{enumerate}

%

%
\section{Complements on \PCA}\label{sec:comp_pca}
%

\subsection{Preliminaries}

Even if there exists a rigorous definition in relation with the theory of the measure, let us assume that a random variable is a probability distribution on the outcome of a random experiment. Typical examples are dice rolling or coin flipping. In dice rolling, $X$ is the observed value on the upwards face of the dice when it stops rolling. The possible outcomes are $\{1,2,3,4,5,6\}$ and the random variable $X$ is defined by
\begin{equation}
    p_1= \P(X=1), \quad p_2 = \P(X=2), \quad \ldots 
\end{equation}
If the dice rolls $n$ times, the outcome is a tuple in $\{1,\ldots,6\}^n$, referred to as $n$ realisations if $X$. They are called i.i.d. for independent identically distributed. 

\nS Let $X$ be a random variable with values in $\R$. Let $x=(x_i)_i$ with $1 \leq i \leq n$ be a vector of $n$ realizations of $X$ (it is called a sample). Then, the sample mean is
\begin{equation}
    \E(x) = \frac{1}{n}\sum_{i=1}^n x_i
\end{equation}
It must not be confounded with the population mean $\E(X)$. For sake of clarity, we denote
\begin{center}
  \begin{tabular}{|cl}
    $\E(X)$ & the population mean \\
    $\E(x)$ or $\overline{x}$ & the sample mean
  \end{tabular}
\end{center}
A sample is a subset of a possibly infinite population. Statistics are about inferring properties for the population knowing a sample. For example, if an infinite population follows a Gaussian law of mean $\mu$ and variance $\sigma^2$, the population mean is $\mu$ whereas the sample mean is $\overline{x}=(1/n)\sum_i x_i$. The sample mean is an unbiased estimator of the mean of the population.\\
\\
The population variance is defined as
\begin{equation}
      \var X = \E\left((X-\E(X))^2\right) 
\end{equation}
The sample variance is defined as
\begin{equation}
      \var x = \frac{1}{n-1} \sum_{i=1}^n (x_i-\overline{x})^2
\end{equation}
It is an unbiased estimator of the population variance. \\
\\
The standard deviation is the square root of the variance
\begin{equation}
    \std x = \sqrt{\var x} = \sqrt{\frac{1}{n-1} \sum_{i=1}^n (x_i-\overline{x})^2}
\end{equation}
If $\E(x)=0$, $\var x = \frac{1}{n-1}\sum_i x_i^2 = \|x\|^2/(n-1)$, which establishes a direct link between the geometric approach relying on Frobenius norm and statistical approach relying on variance. Let $E_r$ be a $r-$dimensional subspace in $\R^p$. Each point $x_i \in \R^p$ has a projection $\widetilde{x_i}$ on $E_r$. Pythagore theorem 
\[
\|x_i\|^2 = \|x_i - \widetilde{x_i}\|^2 + \|\widetilde{x_i}\|^2
\]
can be read as a decomposition of the variance of $X=(x_1, \ldots,x_n)$. The projection is optimal when the discrepancy $\|x_i - \widetilde{x_i}\|^2$ is minimal, or $\|\widetilde{x_i}\|^2$ is maximal, i.e. $\var \widetilde{x_i}$ maximal. This leads to statistical approach of PCA.\\
\\
The covariance of two random variables $X$ and $Y$ is
\begin{equation}
       \cov (X,Y) = \E((X-\E(X))(Y-\E(Y))) 
\end{equation}
and
\begin{equation}
       \cov (x,y) = \frac{1}{n-1}\sum_{i=1}^n \, (x_i-\overline{x})(y_i-\overline{y})
\end{equation}
It is a symmetric, semi-definite, positive bilinear form. It is definite, hence an inner product, for those samples with $\overline{x}=0$: the map $(x,y) \longrightarrow \cov (x,y)$ is an inner product, and $x \longrightarrow \var x$ is a norm, proportional to the Frobenius norm. Two random variables are independent if $\cov (X,Y)=0$. In such a case
\begin{equation}
    \var (X+Y) = \var X + \var Y
\end{equation}
This is Pythagore theorem. 

%
\subsection{Statistical approach}\label{sec:comp_pca:stat}
%

The statistical approach to PCA is a complex topic with many guises. Some rely on empirical distributions, without inference, and some rely on inference, with or without latent variables. Some are presented here.


\nS Statistical approach has been set first with data sets being realisation of random variables. Let us have a $p-$variate  random variable (rv) $X=(X_1, \ldots,X_p)$ with zero mean: $\E(X)=(0, \ldots,0)$, observed on $n$ independent items . A typical example is a $p-$variate Gaussian distribution with zero mean and $\Sigma$ as variance-covariance matrix. The $n$ observations for variable $X_j$ are the column $j$ of a matrix $X$. Each variable $X_j$ is centered: $\E(X_j)=0$. \PCA at rank $r$ is about finding $r$ independent (non correlated) linear combinations $\widetilde{X_k} = Xu_k = \sum_{j=1}^p u_{kj}X_j$, $u_k \in \R^p$, $1 \leq k \leq r$, of the $X_j$ which have maximum variance under the constrain $\|u_k\|=1$. Here, means, variances, covariances are sample mean, sample variances and sample covariances. Statistical PCA relies on the observation that the sample variance of a realisation $x \in \R^p$ of a centered $p-$variate rv is proportional to the square of its Frobenius norm: if $\E(x)=0$, $\var x \propto  \|x\|^2=\langle x,x\rangle$. 

\nS Let us denote as in \cite{Anderson1958} the variance-covariance of $X$ as $\Sigma$, which is standard in statistics\footnote{$\Sigma$ is the standard notation for the diagonal matrix of singular values of a given matrix as well. Confusion can easily be avoided from context.}: $\Sigma=X^\t X$ as $X$ is centered. There is no need for the random variable $X$ to be Gaussian, although some more can be said if it is Gaussian. Let $u \in \R^p$ and let us consider the r.v. $Xu = \sum_ju_jX_j$. It is centered as $\E(Xu) = \sum_ju_j\E(X_j)=0$, and its variance is
\[
 \begin{aligned}
   \var Xu &= \langle Xu \, , \, Xu \rangle \qquad \mbox{as } Xu \mbox{ is centered} \\
   &= \langle u , X^\t Xu \rangle  \\
   &= \langle u , \Sigma u \rangle
 \end{aligned}
\]
Maximizing $\var Xu$ with the constraint $\|u\|=1$ yields that $u$ is the eigenvector of $\Sigma$ associated to its largest eigenvalue, denoted $\lambda$. Then, as $\Sigma u = \lambda u$, $\var Xu = \langle u,\Sigma u\rangle = \lambda \langle u,u\rangle = \lambda$.

\nS The general result is \cite[th. 11.2.1]{Anderson1958}: Let $X$ be a random variable on $\R^p$, with $\E(X)=0$. Let us denote $\var X= \Sigma$. Then, there exists a rotation $U \in \O(\R^p)$ defining 
\[
 V=XU
\]
such that the covariance matrix of $V$ is diagonal and the $k$th component of $V$ has maximum variance among all normalized linear combinations uncorrelated with $V_1,\ldots,V_{k-1}$.

\notes A standard textbook for statistical appproach to \PCA is \cite{Anderson1958}. This section is adapted from \cite[chap. 11]{Anderson1958}. Another key reference is \cite{Rao1964}.

\subsection{Factor Analysis and Probabilistic PCA}

Let us have a $p-$variate Gaussian variable, with $n$ iid realisations, considered as observations. Correlations between observed variables are expressed by the variance-covariance matrix $\Sigma$. Factor analysis (FA) models a situation where there exists $r$ independent unobserved variables, with $r < p$, the response to which explains the correlation between the observed variables. It is a latent variable statistical model. Unobserved variables are called \emph{hidden variables} or \emph{latent variables}. Let $z$ denote the latent variable, and $x$ the observed variable (the data). Let us assume that random variable $Z$ follows a certain law, denoted 
\[
p(z_i \mid \theta)
\]
(it will be specified later). Then, the model for $x_i$ in FA can be written
\[
p(x_i \mid z_i, \theta')
\]
and solving FA is about inferring parameters in $\theta,\theta'$ and recovering the latent variables $z_i$. In Bayesian perspective, $p(z \mid \theta)$ is called the \emph{prior}. 

\nS This makes FA rather complex, and this complexity has nourished decades of debates about Factor Analysis and PCA. Therefore, many authors (including \cite[chap. 7]{Jolliffe2002}) have insisted on the difference between PCA and FA. In fact, FA is a statistical approach of PCA where principal axis (and their realizations as components) are latent variables. 

\nS 
Let us select as prior for the latent variables a multivariate Gaussian probability
\begin{equation}
    p(z_i) = \NC(z_i \mid \mu_0, \Sigma_0)
\end{equation}
(here, $\theta = (\mu_0, \Sigma_0)$). Then, we select for $p(x \mid z)$ a Gaussian law as well, denoted
\begin{equation}
    p(x_i \mid z_i, \theta) = \NC(Wz_i+\mu, \Psi)
\end{equation}
($\theta' = (W, \mu, \Psi)$) where the mean has been selected as a linear function of the hidden input ($x_i \in \R^p$, $z_i \in \R^r$ and $W \in \R^{p \times r}$), and $\Psi$ to be diagonal. The special case where $\Psi = \sigma^2 \I_p$ is known as \emph{probabilistic PCA}. If $\sigma \rightarrow 0$, probabilistic PCA converges to PCA which is non-probabilistic. We then have the following reductions\\
\[
\begin{CD}
\mbox{FA:} \; \Psi @>\Psi=\sigma^2\I >> \mbox{probabilistic PCA} @>\sigma=0>> \mbox{PCA}
\end{CD}
\]

\nS Marginalizing over the latent variable leads to
\begin{equation}
  \begin{array}{lcl}
        p(x_i \mid \theta) &= & \displaystyle \int_{z_i} p(x_i \mid z_i)\, p(z_i \mid \theta) \, dz_i\\
        &=& \NC(x_i \mid W\mu_0 + \mu \, , \, \Psi + W\Sigma_0 W^\t)
  \end{array}
\end{equation}
This shows that t is possible to chose $\mu_0=0$ and $\Sigma_0=\I_r$ without loss of generality:
\begin{equation}
    p(z_i \mid \theta) = \NC(0, \I_r)
\end{equation}
which leads to
\begin{equation}
    p(x_i \mid z_i,\theta') = \NC(Wz_i + \mu \, , \, WW^\t + \Psi)
\end{equation}
and
\begin{equation}
    p(x_i \mid \theta) = \NC(x_i \mid \mu, WW^\t + \Psi).
\end{equation}
The link with \PCA can be read in the observation that the variance-covariance matrix of the data $x_i$ is given by a low-rank matrix $WW^\t + \Psi$, as $W \in \R^{p \times r}$ and $\rank WW^\t=r$. The variance-covariance structure of the observed variables is split into a variance structure for each variable given by $\Psi$ and the covariance structure given by $W$.

\nS Next step is to estimate the parameters $(W,\mu,\sigma)$ or $(W,\mu,\Psi)$ of the model knowing the observations, not presented here. This is not so easy, and is presented in \cite[sect. 12.2.1]{Bishop2006}. It can be done with EM procedure (see \citealp[sec. 12.1.5]{Murphy2012}.

\notes Factor analysis has been developed all over the 20th century in parallel with PCA, sometimes with some controversies. One probable reason for those controversies is that sometimes it is said to inherit from the work of Pearson in 1901, and sometimes of Spearman in 1904, who was interested in finding one factor explaining a diversity of correlated features. His work has been extended to multivariate factors by Thurstone in 1935. The coexistence of PCA and FA with the same tools has probably contributed to some controversies. An history of FA can be found in the introduction of \cite{Basilevsky1994}. A sound utilisation of notions in statistical modeling has clarified the situation. A classical presentation with this approach is \cite[sect. 4.7]{Anderson1958}, or \cite{Basilevsky1994}. It is currently developed as a statistical model with latent variables, where the correlated features are the observed variables, and the independent factors the latent variables. Classical setting is with Gaussian models. Probabilistic PCA (PPCA) has been proposed in 1999 in \cite{Tipping1999}, who have presented as well the links and differences between PCA, PPCA and FA. See as well \cite[sect. 12.2]{Bishop2006} for a clear presentation of statistical approach of \texttt{PCA}, probabilistic \texttt{PCA}, and factor analysis. We have followed \cite{Bishop2006} and \cite[chap. 12]{Murphy2012} here.

\subsection{Distribution of eigenvalues of random matrices}

A random matrix is the realization of a random variable over a space of matrices. Different spaces have been studied, and classically referred to for historical reasons as ``ensembles''. For example, the Ginibre ensemble is the set of $n \times n$ complex matrices whose entries are i.i.d. realizatons of the complex Gaussian normal law, i.e. $a_{k\ell} = x_{k\ell} + i \,  y_{k\ell}$ with $x_{k\ell}, y_{k\ell} \sim \NC(0,1)$. There are several types of results about the distribution of the eigenvalues of random matrices, like :
\begin{itemize}[label=$\rightarrow$]
\item what is the distribution of the eigenvalues of a matrix in a given ensemble for size $n$ ?
\item when $n \rightarrow \infty$, does the distribution of the eigenvalues converge to a given ``universal'' distribution ?
\end{itemize}

\paragraph{Wishart matrix:} Let $A \in \R^{n \times p}$ be a random matrix, with rows being i.i.d. realisations of a $p-$variate Gaussian distribution $X$ of mean $\mu=0$ and variance-covariance matrix $\Sigma$. Let us denote $C = A^\t A \in \R^{p \times p}$. Then, elements in $C$ follows a Wishart distribution, denoted $W(n,p,\Sigma)$, which has been analytically computed by Wishart in 1928. For $\Sigma = \I_p$, it is given by 
\begin{equation}
    p(C) = w(n,p) \, (\det C)^{(n-p-1)/2} \, \exp -\frac{1}{2}\Tr C,
\end{equation}
where $w(n,p)$ is the normalizing constant
\begin{equation}
    \frac{1}{w(n,p)} =  \pi^{p(p-1)/4} \, 2^{np/2} \, \prod_{j=1}^p \Gamma \left(\frac{n-j+1}{2}\right).
\end{equation}
The formula for a general $\Sigma$ is
\begin{equation}
    p(C) = w(n,p,\Sigma) \, (\det C)^{(n-p-1)/2} \, \exp -\frac{1}{2}\Tr \Sigma^{-1}C,
\end{equation}
with 
\begin{equation}
    \frac{1}{w(n,p,\Sigma)} =  \pi^{p(p-1)/4} \, 2^{np/2} \, (\det \Sigma)^{n/2} \, \prod_{j=1}^p \Gamma \left(\frac{n-j+1}{2}\right).
\end{equation}


\paragraph{Mar\v{c}enko-Pastur semi-circular law:} Let $A \in \R^{n \times p}$ be a random matrix as above with its rows being i.i.d. realizations of a $p-$variate Gaussian law of mean 0 and variance-covariance $\Sigma$. Let us denote by
\[
\lambda_1 > \ldots \lambda_p > 0
\]
the $p$ eigenvalues, assumed to be separated. Let the size $n \times p$ of the matrix $A$ grow with 
\[
\lim_{n \rightarrow \infty} \; \frac{p(n)}{n} = \alpha.
\]
Let
\[
C_p = \frac{1}{n}A_p^\t A_p
\] 
be an estimator of the variance covariance matrix. Let us assume that $0 < \alpha \leq 1$. Let $\Sigma = \I_p$, and define 
\[
\begin{cases}
   a &= (1-\sqrt{\alpha})^2 \\
   b &= (1+\sqrt{\alpha})^2. \\
\end{cases}
\]
Let $\mu$ be the measure such that
\begin{equation}
    \mu(\Omega) = \frac{1}{p}\# \{i \mid \lambda_i \in \Omega\}
\end{equation}
for $\Omega \subset \R$. Then, 
\begin{equation}
\lambda \in [a,b],
\end{equation}
and, if $x \in [a,b]$
\begin{equation}
    d\mu(x) = \frac{1}{2\pi\alpha}\frac{\sqrt{(b-x)(x-a)}}{x}.
\end{equation}
So, if $a \leq z < z' \leq b$
\begin{equation}
    \P(z \leq \lambda \leq z' ) = \frac{1}{2\pi\alpha} \int_z^{z'} \frac{\sqrt{(b-x)(x-a)}}{x} \, dx.
\end{equation}

\notes The Wishart matrix has been thoroughly studied, and much can be said about it. See e.g. \cite[chapter 7]{Anderson1958} devoted to it. The expressions of the joint law of the coefficients of a Wishart matrix in the case where $\Sigma=\I_p$ and for any $\Sigma$  have been borrowed from \cite[section 7.2, formula (1)]{Anderson1958} and \cite[section 2.2]{Meckes2019}. For the asymptotic distribution of the eigenvalues of a Wishart matrix through Mar\v{c}enko-Pastur theorem, we have followed \cite{Meckes2021}, which is a gem.


\subsection{Unitarily invariant norms}

The problem of \PCA as set in section \ref{sec:pca:set} can be set for any norm, and not Frobenius norm only. Most widely used norms in data analysis are $\ell^1$ and $\ell^\infty$ norms, on top of $\ell^2$ norms. However, there are very few norms for which exact solution and efficient algorithms to compute a solution are known. One exception is the family of unitarily invariant norms.

\nT{Unitarily Invariant norm (UIN)} There is a generalization of Eckart-Young theorem through unitarily invariant norms. The framework is that of vector spaces on $\C$, but it can be applied on vector spaces on $\R$ as well. $\U(\C^n)$ denotes the set of unitary matrices in $\C^n$, i.e. matrices having the property $UU^*=U^*U=\I_n$, where $U^*=\overline{U^\t}$ (and the same for $\C^p)$. The equivalent in $\R^n$ is the set $\O(\R^n)$ of orthogonal matrices such that $U^\t U = \I_n$.

\nS Let $E,F$ be two complex vector spaces, $A \in E \otimes F \simeq \L(F,E)$. A norm $\|.\|$
\[
 \begin{CD}
   E \otimes F @>\|.\|>> \R^+
 \end{CD}
\]
is said unitarily invariant if
\begin{equation}
 \forall \:
 \begin{cases}
   U &\in \U(\C^n) \\
   V &\in \U(\C^p) \\
   A & \in \C^{n \times p}
 \end{cases},
\qquad\qquad \|UAV^*\|= \|A\| 
\end{equation}
For example, spectral and Frobenius norms are unitarily invariant norms. If $E,F$ are real vector spaces, the norm is said invariant by orthogonal transformation if
\begin{equation}
 \forall \:
 \begin{cases}
   U &\in \O(\R^n) \\
   V &\in \O(\R^p) \\
   A & \in \R^{n \times p}
 \end{cases},
\qquad\qquad \|UAV^\t\|= \|A\| 
\end{equation}

\nT{Symmetric gauge function (SGF)} A norm
\[
 \begin{CD}
  \R^n @>\Phi>> \R^+
 \end{CD}
\]
is called a symmetric gauge function if it is invariant by any permutation of the coordinates in $\R^n$, i.e., if $\mathscr{P}$ is the set of permutations in $\R^n$
\begin{equation}
 \forall \: P \in \mathscr{P}, \quad \Phi(Px)=\Phi(x)
\end{equation}

\nS There is a remarkable link between UIN and SGF. Let $A \in \C^{n \times p}$ (or $\in \R^{n \times p}$) and 
\[
 \Sigma = (\sigma_1, \ldots, \sigma_p)
\]
its singular values. Let $\Phi$ be a SGF. To each SGF $\Phi$, one associates the norm $\|.\|_\Phi$ defined by
\begin{equation}
 \|A\|_\Phi = \Phi(\sigma_1, \ldots, \sigma_n)
\end{equation}
Then, $\|.\|_\Phi$ is a UIN. Let us note that if $\|.\|$ is a UIN, $A$ a matrix in $\C^{n \times p}$ and $A=U\Sigma V^*$ the SVD of $A$, then $\|A=U^*AV\|$ and, as $U^*AV=\Sigma$, $\|A\|=\|\Sigma\|$ i.e. is a function of its singular values.

\nT{Mirsky's theorem} Mirsky has shown: Let $E,F$ be two vector spaces on $\C$ or $\R$, and $A,B \in E \otimes F \simeq \L(F,E)$. Let $(\alpha_1,\ldots,\alpha_p)$ (resp. $(\beta_1,\ldots,\beta_p))$ be the singular values of $A$ (resp. $B$). Then, for any unitarily invariant norm $\|.\|$
\begin{equation}
 \| \diag (\alpha_1-\beta_1, \ldots,\alpha_p-\beta_p)\| \leq \|A-B\|
\end{equation}

\nT{Schmidt-Mirsky theorem} Let
\begin{itemize}[label=$\rightarrow$]
 \item $A \in E \otimes F \simeq \L(F,E)$
 \item $(\sigma_1,\ldots,\sigma_p)$ the singular values of $A$ in non increasing order
 \item $B \in E \otimes F$ with $\rank B =r \leq p$
 \item $\Phi$ a Symmetric Gauge Function with $\|.\|_\Phi$ as associated unitarily invariant norm
\end{itemize}
Then
\begin{equation}
 \|A-B\|_\Phi \geq \Phi(0, \ldots, 0, \sigma_{r+1},\ldots,\alpha_p)
\end{equation}
Moreover, if $A=U\Sigma V^*$ is the SVD of $A$, the equality is reached fot matrix $A_r$ defined as
\begin{equation}
 A_r=U\Sigma_rV^*
\end{equation}
where $\Sigma_r$ is the diagonal matrix obtained from $\Sigma$ be setting to 0 all singular values beyond $r$.

\nB Then, $A_r$ is the best rank $r$ approximation of $A$ for norm $\|.\|_\Phi$ as well.

\notes The link between UIN and SGF has been shown in J. von Neumann (1937), Some Matrix Inequalitues and Metrization of Matrix Spaces, \emph{Tomsk Univ. Rev.}, \textbf{1}286-300. The extension of PCA to UIN and SGF is nicely presented with many references in \cite{Claret1987}. See  \cite{Schatten1960} as well for an algebraic survey.

\clearpage 

%
\section{\PCA with Instrumental Variables}\label{sec:pcaiv}
%

\subsection*{Notations}

Here are the notations chosen for this section. They are as compatible as possible with notations selected in other sections, especially PCA. There are explained along the text. 

\begin{center}
\ovalbox{
\begin{tabular}{cll}
Symbol & in space & What it is \\
\hline 
$A$ & $\R^{n \times p}$ & matrix to be analyzed \\
$A_r$ & $\R^{n \times p}$ & best approximation of $A$ of rank $r$ with constraints \\ 
$F$ & $\subset \R^n$ & subspace of constraints on principal components \\
$H$ & $\subset \R^p$ & subspace of constraints on principal axis \\
$I$ & $\R^{n \times m}$ & matrix of instrumental variables \\
$\Lambda$ & $\R^{q \times q}$ & diagonal matrix of eigenvalues of PCA-IV \\
$m$ & $\N$ & dimension of $F$ \\
$P_\f$ & $\R^{p \times p}$ & projector on $F$ \\
$P_\h$ & $\R^{n \times n}$ & projector on $H$ \\
$\PC_{\f \otimes \h}$ & $\R^{np \times np}$ & projector on $F \otimes H$ \\
$q$ & $\N$ & dimension of $H$ \\
$r$ & $\N$ & prescribed rank for best approximation \\
$U_\f$ & $\R^{n \times m}$ & matrix of an orthonormal basis of $F$ \\
$V$ & $\R^{p \times q}$ & matrix of principal axis of PCA-IV of $(A, F,H)$ \\
$V_\h$ & $\R^{p \times q}$ & matrix of an orthonormal basis of $H$ \\
$Y$ & $\R^{n \times m}$ & principal components of PCA-IV of $(A, F,H)$ \\
$Y_r$ & $\R^{n \times r}$ & first $r$ principal components of PCA-IV of $(A, F,H)$
\end{tabular}
}
\end{center}

\vspace*{5mm}

Here, we use tensor notation for PCA. For sake of clarity for those readers not familiar with those notations, we set the problem with classical notations like  $yv^{\textsc{t}}$ for a rank one matrix as a starting point. Let $y \in \R^n$ and $v \in \R^p$. We can accept as a definition that $\otimes$ is a bilinear form
\[
\begin{CD}
 \R^n \times \R^p @>\otimes>> \R^{n \times p} \\
 (y,v) @>>> y \otimes v
\end{CD}
\]
defined by
\[
y \otimes v := yv^\t 
\]
which can be illustrated as 
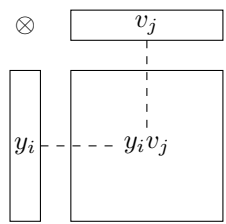
\begin{center}
 \begin{tikzpicture}
   \draw (0,0) rectangle (0.4,2) ;
   \draw (0.8,2.4) rectangle (2.8,2.8) ;
   \draw (0.8,0) rectangle (2.8,2) ;
   \node () at (0.2,2.6){$\otimes$};
   \node () at (0.2, 1){$y_i$};
   \node () at (1.8, 2.6){$v_j$};
   \node () at (1.8, 1){$y_iv_j$};
   \draw[dashed] (.4, 1) -- (1.4,1) ; 
   \draw[dashed] (1.8, 2.4) -- (1.8,1.2) ; 
 \end{tikzpicture}
\end{center}
If $y=(y_1, \ldots,y_n) \in \R^n$ and $v=(v_1,\ldots,v_p) \in \R^p$, then $y \otimes v \in \R^{n \times p}$ and
\[
 (y \otimes v)_{ij} = y_iv_j
\]

\nS Let $A \in \R^{n \times p}$. A rank$-r$  best approximation of $A$ can be written as
\[
 A_r = \sum_{a=1}^r y_av^{\textsc{t}}_a, \qquad y_a = Av_a
\]
or, equivalently
\[
 A_r = \sum_a y_a \otimes v_a, \qquad \mbox{with} \:
 \begin{cases}
  y_a & \in \R^n \\
  v_a & \in \R^p
 \end{cases}
\]
PCA with instrumental variables (\verb%PCAiv%) is setting some constraints on the principal components $(y_a)_a$ or principal axis $(v_a)_a$, i.e. that they live in given subspaces respectively $F \subset \R^n$ and $H \subset \R^p$.\\
\\
Usually, those spaces are given as spanned by a set of vectors in respectively $\R^n$ (for $F$, constraint on $y$) and $\R^p$ (for $H$, constraints on $v$). A classical situation is when there is no constraint on the axis $v$, but only on the components $y$, and $F$ is spanned by a $n \times q$ matrix denoted $I$, the columns of which are called \emph{instrumental variables}.

%
\subsection{Setting the problem}\label{sec:pcaiv:ef}
%

Let us suppose that $\dim F=m$ and $\dim H=q$. \verb%PCAiv% can be stated as
\\
\begin{center}
 \shadowbox{
  \begin{tabular}{ll}
      given & $A \in \R^{n \times p}$ \\
      & $F \subset \R^n, \quad H \subset \R^p$\\
      & $\dim F=m$, $\dim H=q$\\
      & $0 < r < \min(m,q)$ \\
      \\
      find & $A_r \in F \otimes H$ \\
      \\
      such that & $\|A-A_r\|$ minimum
  \end{tabular}
 }
\end{center}
Technically, it is possible to specify this problem with base of $F$ and $H$ having been selected. We assume here that they are orthonormal. If they are not orthonormal, it is possible to build an orthonormal base by $QR$ decomposition, for example with Gram-Schmidt orthogonalization procedure (which can be numerically unstable but is easy to implement), or using Householder reflections (more stable).

\subsection{Solving the problem}

If $F \subset \R^n$ and $H \subset \R^p$. Then $F \otimes G \subset \R^n \otimes \R^p$, where $\subset$ means ``is a vector subspace of''.

\nS Before stating the main result, we need to recall some elementary results on linear projectors. Let $F \subset \R^n$. Then, the projector on $F$ is denoted $\PC_\f$. Let $U_\f=(u_1,\ldots,u_m)$ be an orthonormal basis of $F$ columnwise. Then
\begin{equation}
 \PC_\f = U_\f \, U_\f^\t
\end{equation}
Indeed, if $x \in \R^n$, we have
\[
 \begin{aligned}
   \PC_\f x &= \sum_i \langle u_i,x\rangle \, u_i \\
   &= \sum_i (u_i \otimes u_i) \, x \\
   &= \left(\sum_iu_i \otimes u_i\right) \, x
 \end{aligned}
\]
Then
\begin{equation}
 \PC_\f = \sum_i u_i \otimes u_i = U_\f \, U_\f^\t
\end{equation}

\nS Let $\PC_\f$ be the projector on $F$, $\PC_\h$ be the projector on $H$, and $A \in \R^{n \times p}$. Then, the projector $ \PC_{F \otimes H}$ on $F \otimes H$ is defined by 
\begin{equation}\label{eq:cmplt:proj}
 \PC_{F \otimes H}A = U_\f \, U_\f^\t \, A \, V_\h \, V_\h^\t
\end{equation}

\nT{Main result} Let $A_{F \otimes H}$ be the projection of $A$ on $E \otimes F$. Then, the solution $A_r$ of \verb%PCAiv% is the PCA of $A_{F \otimes H}$.
\begin{proof}
Let us denote by $\PC$ the projection from $\R^n \otimes \R^p$ on $E \otimes F$, and by $\PC^\perp$ the projection on 
$(E \otimes F)^\perp$. Let us recall that $A_r=\sum_jy_j \otimes v_j$. We have $\I = \PC + \PC^\perp$ and 
\[
\begin{array}{lclcl}
  \|A - A_r\|^2 &=& \left\|\left(\PC + \PC^\perp\right)(A-A_r)\right\|^2 & \mathrm{as} & \I = \PC + \PC^\perp\\
  &=& \left\|\PC(A-A_r) + \PC^\perp(A-A_r)\right\|^2 & & \\
  &=& \left\|\PC(A-A_r)\right\|^2 + \left\|\PC^\perp(A-A_r)\right\|^2 & \mathrm{by} & \mathrm{Pythagore}\\
  &=& \left\|\PC(A-A_r)\right\|^2 + \left\|\PC^\perp(A)\right\|^2 & \mathrm{as} & \PC^\perp A_r=0\\
  &=& \left\|\PC(A)-A_r\right\|^2 + \left\|\PC^\perp(A)\right\|^2 & \mathrm{as} & \PC A_r= A_r  
\end{array}
\]
Let us recall that $A, F, H$ are fixed, hence so are $\PC, \PC^\perp$. Then, $A_r$ only can vary. Hence $\|A - A_r\|^2$ is minimum when $\left\|\PC(A)-A_r\right\|^2$ is minimum, and $A_r$ is the \PCA of $\PC(A)=A_{F \otimes H}$.
\end{proof}

\nS This leads to a solution for \verb!PCAiv!. Let $(u_1, \ldots,u_m)$ be an orthonormal basis for $F \subset \R^n$, and $(v_1, \ldots, v_q)$ for $H \subset \R^p$, with $m < n$ and $q < p$. Let $U_\f$ be the $n \times m$ matrix with columns $u_i$ and $V_\h$ be the $p \times q$ matrix with columns $v_j$. Then, the orthogonal projector $\PC_\h$ from $\R^p$ to $H$ is given by matrix
\[
 P_\h = V_\h \, V_h^\t
\]
and the orthogonal projector from $\R^n$ to $F$ is given by matrix
\[
 P_\f = U_\f \, U_\f^\t
\]
Let $A \in  \R^{n \times p}$. Then 
\begin{equation}
  \PC_{F \otimes H}(A) = U_\f \, U_\f^\t \, A \, V_\h \, V_\h^\t
\end{equation}
Hence the algorithm\\
\\
\begin{algorithm}[H]
\begin{algorithmic}[1]
\STATE \textbf{input:}  $A \in \R^{n \times p}$
\STATE \textbf{input:} $U_\f$ orthonormal, $F = \span U_\f \subset \R^n$, 
\STATE \textbf{input:} $V_\h$ orthonormal, $H = \span V_\h \subset \R^p$, 
\STATE \textbf{input}  $r < p$
\STATE \textbf{compute} $P_\f = U_\f \, U_\f^\t$
\STATE \textbf{compute} $P_\h = V_\h \, V_\h^\t$
\STATE \textbf{compute} $T = P_\f \, A \, P_\h$
\STATE \textbf{compute} $(Y,\Lambda,V)=  \textsc{pca\_core}(T)$
\RETURN $Y,\Lambda, V$
\end{algorithmic}
\caption{Pseudocode for \texttt{PCAiv}$(A,U_\f,V_\h)$}
\label{alg:pca:iv}
\end{algorithm}

\nS In this algorithm, the SVD is run in the PCA of $T=P_\f \, A \, P_\h \in \R^{n \times p}$, same dimensions as $A$. However, as $T \in F \otimes H$, with $\dim F=m$ and $\dim H=q$, the rank of $T$ is $q < p$ (if we assume that $q \leq m$). Here, we show how to run SVD on a matrix of dimension lower than $n \times p$. We start from 
\[ 
\begin{array}{lcl}
T  &=&  P_\f \, A \, P_\h \\
&=& (U_\f \, U_\f^\t) \, A \, (V_\h \, V_\h^\t) \\
&=& U_\f (U_\f^\t \, A \, V_\h) V_\h^\t  \\
&=& U_\f \, T' \, V_\h^\t  
\end{array}
\]
denoting
\[
T' = U_\f^\t \, A \, V_\h \qquad \in \R^{m \times q}
\]
Let $(U'_\t, \Sigma'_\t, V'_\t)$ be the SVD of $T'$
\begin{equation}
    T' = U'_\t \, \Sigma'_\t V'^{\t}_\t
\end{equation}
with (this will be useful soon), assuming $m \geq q$ here for sake of simplicity
\[
\begin{cases}
   U_\f & \in \R^{n \times m} \\
   U'_\t & \in \R^{m \times q} \\
   \Sigma'_\t & \in \R^{q \times q} \\
   V'_\t & \in \R^{q \times q} \\
   V_\h & \in \R^{p \times q}
\end{cases}
\]
Then
\begin{equation}
   \begin{array}{lcl}
        T &=& U_\f U'_\t \, \Sigma'_\t \, V'^\t_\t V_\h^\t \\
        &=& (U_\f U'_\t) \, \Sigma'_\t \, (V_\h V'_\t)^\t 
   \end{array}
\end{equation}
and $(U_\f U'_\t \, , \, \Sigma'_\t \, , \, V'^\t_\t V^\t_\h)$ is the SVD of $T$ because
\[
\left\{
  \begin{array}{lcl}
       (U_\f U'_\t)^\t(U_\f U'_\t) &=& U'^\t_\t \, U_\f^\t U_\f \, U'_\t \\
       &=& U'^\t_\t \, \I_m \, U'_\t \\
       &=& U'^\t_\t \, U'_\t \\
       &=& \I_q
  \end{array}
\right.
\]
so $U_\f U'_\t$ is orthonormal. Similarly
\[
\left\{
  \begin{array}{lcl}
       (V_\h V'_\t)^\t(V_\h V'_\t) &=& V'^\t_\t \, V_\h^\t V_\h \, V'_\t \\
       &=& V'^\t_\t \, \I_q \, V'_\t \\
       &=& V'^\t_\t \, V'_\t \\
       &=& \I_q
  \end{array}
\right.
\]
and $V_\h  V'^\t_\t$ is orthonormal too. This small calculation replaces the SVD of $T \in \R^{n \times p}$ in $\OC(n^2p)$ by the SVD of $T' \in \R^{m \times q}$ in $\OC(m^2q)$. This leads to the following algorithm:\\
\\
\begin{algorithm}[H]
\begin{algorithmic}[1]
\STATE \textbf{input:}  $A \in \R^{n \times p}$
\STATE \textbf{input:} $F = \span U_\f \subset \R^n$, $U_\f$ orthonormal 
\STATE \textbf{input:} $H = \span V_\h \subset \R^p$, $V_\h$ orthonormal
\STATE \textbf{input}  $r < p$
\STATE \textbf{compute} $T' = U^\t_\f \, A \, V_\h$
\STATE \textbf{do} the SVD of $T'$: $T'=U'_\t \Sigma'_\t V'^\t_\t$
\STATE \textbf{compute} $U = U_\f U'_\t$
\STATE \textbf{compute} $V = V_\h V'_\t$
\STATE \textbf{compute} $\Lambda = \Sigma'^2_\t$
\STATE \textbf{compute} $Y=U \, \Sigma'_\t$
\RETURN $Y,\Lambda, V$
\end{algorithmic}
\caption{Second Pseudocode for \texttt{PCAiv}$(A,U_\f,V_\h)$}
\label{alg:pca:iv2}
\end{algorithm}

%
\subsection{Interpretation of PCAiv}
%

\PCAiv is a two steps procedures:
\begin{enumerate}
 \item project the variables into the space spanned by the instrumental variables
 \item run the \PCA of the projected variables
\end{enumerate}
This can be sketched by the following diagram:\\
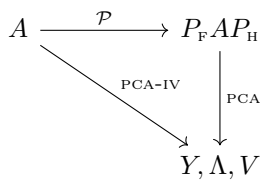
\begin{center}
 \begin{tikzcd}
   A \arrow[rr, "\PC"] \arrow[rrdd, "\textsc{pca-iv}"] & & P_\f A P_\h \arrow [dd, "\textsc{pca}"]\\
   & & \\
   & & Y,\Lambda,V
 \end{tikzcd}
\end{center}
This leads to
\begin{enumerate}
    \item another interpretation of \PCAiv
    \item an estimate of the quality of the \PCAiv
\end{enumerate}
as follows.

\nS Let us give another interpretation of \PCAiv on the example of a constraint on the principal components only. Let $a \in \R^n$. The projection of $a$ on $F$ is given by $P_\f a$. Let now $A \in \R^{n \times p}$. The matrix $\widetilde{A} = P_\f \, A$ is in $\R^{n \times p}$, like $A$, and its column $j$ is the projection of column $j$ of $A$ on $F$. It is the regression of the column $j$ by the basis vectors of $F$, given by the columns of $U_\f$. The first principal component of $\widetilde{A}$ is the best summary of $A$ by one single linear combination of the columns of $U_\f$. It is the best same linear regression applied to all columns of $A$.  In this interpretation, \PCAiv is equivalent to PLS.

\nS For the \texttt{PCA}, classical estimators of the quality of the \PCA can be used. Let us denote $Y_r$ the matrix of $r$ first principal components of $P_\f A P_\h$ on which \PCA is run, and 
\begin{equation}
 \|Y_r\| = \rho_r \|P_\f A P_\h\|
\end{equation} 
i.e. the quality of the \PCA is denoted by $\rho_r$ at rank $r$.

\nS The quality of the projection can be quantified by
\begin{equation}
 \|P_\f A P_\h\| = \theta \|A\|
\end{equation}
($\theta$ is the cosine of the angle between $A$ and $P_\f A P_\h$).

\nS We then have
\begin{equation}
 \|Y_r\|=\rho_r\theta \|A\|
\end{equation}
and the quality of the \PCAiv can be poor for two reasons:
\begin{itemize}[label=$\rightarrow$]
 \item the quality $\theta$ of the projection is poor
 \item the quality $\rho_r$ of the \PCA at rank $r$ is poor.
\end{itemize}

\nS It is essential to distinguish between the quality of the projection and the quality of the \texttt{PCA}. Let us see it on an example. We assume that the matrix $A$ is built as a low rank matrix plus some important noise. It can be expected that the noise is poorly projected (there is no specific subspace where the noise is better represented), whereas there is some specific low dimensional subspace where the low rank component of $A$ is well represented. Then, projection will filter out the noise, whereas the \PCA of the projected matrix will find the low rank property of the structure of $A$. $\theta$ will be low, but $\rho_r$ close to 1. Even if the quality $\rho_r\theta$ is poor because of the poor quality $\theta$ of the projection, \PCAiv is a success as it has filtered out the noise and exhibited the low rank of the data set.

%
\subsection{Non orthonormal basis}
%

In section \ref{sec:pcaiv:ef}, the subspaces $F$ and $H$ are given by their basis $U_\f$ and $V_\h$ respectively. Let us now suppose that $F$ and $H$ are respectively spanned by the columns of $U'_\f$ and $V'_\h$, which are no longer assumed to be orthonormal. Then
\[
 P_\f = U'_\f(U'^\t_\f U'_\f)^{-1}U'^\t_\f, \qquad P_\h = V'_\h(V'^\t_\h  V'_\h)^{-1}V'^\t_\h
\]
and equation (\ref{eq:cmplt:proj}) reads
\[
 \begin{aligned}
  T &= P_\f A P_\h\\
  &=  U'_\f(U'^\t_\f U'_\f)^{-1}U'^\t_\f \; A \; V'_\h(V'^\t_\h  V'_\h)^{-1}V'^\t_\h
 \end{aligned}
\]

\nS However, this equation, even if correct, is not efficient for numerical analysis where it should be avoided, because of the cost of three products and one inversion per projector. It is far more efficient to work with orthogonal basis of $F$ and $H$, by 
\begin{enumerate}
 \item building an orthonormal basis $U_\f$ of $F$ and an orthonormal basis $V_\h$ of $H$ by QR decomposition or Gram-Schmidt orthonormalisation,
 \item building the projectors $P_\f = U_\f U^\t_\f$ and $P_\h = V_\h V^\t_\h$
\end{enumerate}

\notes Apparently, the term ``PCA with Instrumental Variables'' appeared first in \cite{Rao1964}. It has been studied with a double set of constraints in \cite{Sabatier1984}, and widely used in ecology for example (see \cite{Lebreton1991}). \PCAiv is equivalent to PLS.

%
\section{PCA with metrics on rows and columns}\label{sec:pcamet}
%

\subsection*{Notations}

\vspace*{5mm}

\begin{center}
  \ovalbox{
    \begin{tabular}{cll}
        symbol & in space & meaning \\
        \hline
        $\|.\|_\n$ & & norm induced by $N$ in $\R^n$ \\
        $\|.\|_\p$ & & norm induced by $P$ in $\R^p$ \\
        $\|.\|_\t$ & & norm induced by $T$ in $\R^{n \times p}$ \\
        $A$ & $\R^{n \times p}$ & matrix to be analyzed \\
        $B$ & $\R^{n \times p}$ & matrix for calculations: $B=MAQ$ \\
        $\Lambda$ & $\R^{p \times p}$ & eigenvalues of $B^\t B$ \\
        $M$ & $\R^{n \times n}$ & unique SDP with $M^2=N$ \\
        $N$ & $\R^{n \times n}$ & SDP defining an inner product in $\R^n$ \\
        $\PC_v$ & $\L(\R^p)$ & projector on $\R v$ in $\R^p$ for inner product $P$\\
        $\PC_\f$ & $\L(\R^p)$ & projector on $F \subset \R^p$ in $\R^p$ for inner product $P$\\        
        $P$ & $\R^{p \times p}$ & SDP defining an inner product in $\R^p$ \\
        $Q$ & $\R^{p \times p}$ & unique SDP matrix with $P=Q^2$ \\
        $T$ & $\R^{np \times np}$ & SDP defining an inner product in $\R^{n \times p}$ \\
        $V$ & $\R^{p \times p}$ & principal axis of $A$ for inner product defined by $(N,P)$ \\
        $w$ & $\R^p$ & weights for an inner product in $R^p$ \\ 
        $W$ & $\R^{n \times p}$ & principal components of $B$ with standard inner product \\
        $X$ & $\R^{p \times p}$ & principal axis of $B$ with standard inner product \\
        $Y$ & $\R^{n \times p}$ & principal components of $A$ with inner product defined by $(N,P)$ \\
        $Z$ & $\R^{np \times np}$ & unique SDP with $Z^2=T$; $Z(A)=MAQ$
    \end{tabular}
   } 
\end{center}

\vspace*{5mm}

A matrix $A \in \R^{n \times p}$ can be considered as an element of space $\R^n \otimes \R^p$. This space is implicitly endowed with the canonical inner product
\[
\begin{aligned}
 \forall \: A, B \in \R^{n \times p}, \quad \langle A, B \rangle & = \sum_{i,j}\alpha_{ij}\beta_{ij} \\
 &= \Tr A^\t B \\
 & = \Tr B^\t A
 \end{aligned}
\]
Any $np \times np$ symmetric definite positive matrix $T$ defines an inner product on $\R^{n \times p}$. Componentwise, it is defined as
\[
\langle A \, , \, B \rangle_\t = \sum_{ij,k\ell}\,T_{ij,k\ell} \,  \alpha_{ij}\beta_{k\ell}
\]
Then, $\R^{n \times p}$ is endowed  with a Euclidean structure, which induces a norm, which induces a metric structure with
\[
d_\t(A,B)= \|A-B\|_\t
\]
Then, linear dimension reduction can be performed, e.g. on a matrix $A$.  In the algebraic framework, it consists in finding the rank $r$ matrix in $\R^{n \times p}$ which is the closest to $A$ with the distance induced by $T$. The geometric framework consists in finding an affine subspace of $\R^{n \times p}$ of dimension $r$ on which the projection of the point clod associated to $A$ is optimal. The statistical viewpoint will not be developed here.

\nB We restrict ourselves here to inner products which establish a link with possible inner products in $\R^n$ and $\R^p$.

\subsection{Metrics and weights on row and column spaces}\label{sec:pcamet:notations}

Let us define first an inner product in $\R^p$, by a SDP matrix $P$, i.e.
\begin{equation}
 \langle x, y \rangle_\p = \langle x \, , \, Py \rangle
\end{equation}
As $P$ is symmetric, we have $P=P^\t$, and $\langle x, y \rangle_\p = \langle Px \, , \, y \rangle$ as well. This means, component-wise in a given basis, that
\begin{equation}
  \langle x, y \rangle_\p = \sum_{i,j=1}^p p_{ij}\, x_i\, y_j, \qquad p_{ji} = p_{ij}
\end{equation}
if $x=(x_i)_i$, $y=(y_j)_j$ and $P=(p_{ij})_{i,j}$. If $P=\I_p$, canonical inner product is recovered, as $p_{ij} = \delta_i^j$. A case worth being studied in details is when $P$ is diagonal, i.e. $P = \diag w$ with $w = (w_1,\ldots,w_p)$, or
\[
 p_{ij} = 
 \begin{cases}
  w_i & \mbox{if}\: i=j \\
  0 & \mbox{if}\: i \neq j \\
 \end{cases}
\]
In such a case
\begin{equation}
 \langle x, y \rangle_w = \sum_{i=1}^pw_i \, x_i \, y_i
\end{equation}

\nS The norm $\|.\|_\p$ is defined as
\begin{equation}
    \|x\|_\p^2 = \langle x,x\rangle_\p = \langle x,Px\rangle
\end{equation}
If $P=\mathrm{diag}\:w$, this yields
\begin{equation}
    \|x\|_w^2 = \sum_i\, w_ix_i^2
\end{equation}

\nS As $P$ is SDP, there exists a unique SDP matrix $Q$ such that $P=Q^2$. Then
\begin{equation}
 \langle x, y \rangle_\p = \langle Qx \, , \, Qy \rangle
\end{equation}
and
\begin{equation}
 \|x\|_\p = \|Qx\|
\end{equation}

\nB One may wonder whether the metric should be given by $P$ or by $Q$. It is tempting to give it by $Q$, because $Q$ establishes an isometry $\iota_\q$ between $(\R^p,Q)$ and $(\R^p,\I)$ by
\begin{equation*}
 \begin{CD}
  (\R^p,Q) @>\iota_\q>> (\R^p,\I) \\
  x @>>> Qx
 \end{CD}
\end{equation*}
as
\[
 \|x\|_\q = \|Qx\| = \|\iota_\q x\|
\]
However, in data analysis, it is classical to use weights, i.e. define metrics with diagonal matrices. Let
\[
 w = (w_1, \ldots,w_p) \in \R^{p+}
\]
Then, a distance between $x,x' \in \R^p$ with weights $w$ is given by
\[
d_w(x,x') = \|x-x'\|_w = \sqrt{\|x-x'\|^2_w} =\sqrt{\sum_iw_i(x_i-x'_i)^2}
\]
This is consistent with an isometry defined by
\[
 Q = \diag \left(\sqrt{w_1}, \ldots,\sqrt{w_p}\right)
\]
as
\[
 d_w(x,x') = \sqrt{\sum_i\left(\sqrt{w_i}x_i-\sqrt{w_i}x'_i\right)^2}
\]
It is customary to define the weights by $w$, which are the diagonal elements of $P$, and not $\sqrt{w}$. Hence, it is consistent with this classical approach to define the metric by $P$, hence denote $\langle x,x' \rangle_P$. We will use $P$ or $Q$ indifferently, knowing that $P=Q^2$.

\nT{Projection operator with metrics} We define here the projection operator in $\R^p$ endowed with metrics defined by $P$. Let $v \in \R^p$ with $\|v\|_\p=1$. Let us denote by $\PC_v$ the projection operator in $\R^p$ on $\R v$
\[
\begin{CD}
  \R^p @>\PC_v>> \R v \\
  x @>>> \lambda v
\end{CD}
\]
$\PC_v x$ is the vector $\lambda v$ such that $\|x-\lambda v\|_\p$ is minimal. As $\|x-\lambda v\|_\p^2 = \|x\|^2_\p + \lambda^2 -2\lambda \langle x,v\rangle_\p$, and the unknown is $\lambda$, this yields 
\begin{equation}
    \lambda = \langle x,v\rangle_\p  = \langle x,Pv\rangle
\end{equation}
Then
\begin{equation}
    \PC_v x = \langle x,Pv\rangle \, v = (v \otimes Pv).x
\end{equation}
and
\begin{equation}
    \PC_v = v \otimes Pv
\end{equation}
Let us note that $\PC_v \neq Qv \otimes Qv$ (indeed, $(Qv \otimes Qv)x = \langle Qv,x\rangle Qv \in \R Qv \notin \R v$). For a general subspace 
\[
F \subset \R^p = \span (v_1, \ldots,v_m)
\]
with $(v_i)_i$ a basis of $F$ orthonormal for $P$, the same calculation leads to
\begin{equation}
    \PC_\f = \sum_i v_i \otimes Pv_i
\end{equation}

\nS A metric can be defined in the same way in $\R^n$, the column space. It is given by a SDP matrix $N \in \R^{n \times n}$. We denote $N=M^2$, and 
\[
 \langle y,y'\rangle_\n = \langle My, My'\rangle, \qquad \|y\|_\n = \|My\|
\]

\nS Let $P$ define a metric on $\R^p$ and $N$ a metric on $\R^n$. An inner product on $\R^{n \times p}$ will be defined by a SDP matrix $T$, canonically associated to $N$ and $P$ and such that the map
\[
 \begin{CD}
  (\R^{n \times p}, T) @>>> (\R^{n \times p}, \I)
 \end{CD}
\]
is an isometry. Therefore, let $a,x \in \R^n$ and $b,y \in \R^p$. Let $T=Z^2 \in \R^{np \times np}$ be an inner product on $\R^n \otimes \R^p$.  We wish to have on elementary ( $=$ rank one) matrices 
\[
 \langle a \otimes b \, , \, x \otimes y \rangle_\t = \langle a,x\rangle_\n \langle b,y\rangle_\p
\]
with
\begin{equation}\label{eq:TZ}
  \langle a \otimes b \, , \, x \otimes y \rangle_\t = \langle Z(a \otimes b) \, , \, Z(x \otimes y) \rangle
\end{equation}
As
\begin{equation}
 \begin{array}{lcl}
  \langle a \otimes b \, , \, x \otimes y \rangle_\t &=& \langle Ma,Mx\rangle \langle Qb,Qy\rangle \\
  &=& \langle Ma \otimes Qb \; , \; Mx \otimes Qy \rangle
 \end{array}
\end{equation}
(\ref{eq:TZ}) is fulfilled by selecting
\begin{equation}
    \begin{array}{lcl}
       Z(a \otimes b) &=& Ma \otimes Qb  \\
       &=& M(a \otimes b)Q
    \end{array}
\end{equation}
and, for any matrix $A \in \R^{n \times p}$
\begin{equation}
 Z(A)=MAQ
\end{equation}
by linearity.

\nS So, we have
\begin{equation}
 \|A\|_\t = \|MAQ\|
\end{equation}
This permits to solve PCA of a matrix $A \in \R^{n \times p}$ with metrics $N$ on $\R^n$ and $P$ on $\R^p$ by transporting the problem by the isometry $\iota(A)=MAQ$ to PCA of the image $\iota(A)$ and transporting back the solution into initial space by the inverse $\iota^{-1}$ of the isometry.

\nT{Remark} We might stop this section here by saying that all what has been said about \PCA makes no assumption about the inner product defining a Euclidean structure in $\R^p$, $\R^n$ or $\R^n \otimes \R^p = \R^{n \times p}$, and anything can be transported by the isometry mentionned above, period. But it may be not useless to develop this in more details. Even if in such a compact form it is exact and provides all needed information and procedures.

\subsection{Setting the problem}

Here, we use the algebraic approach of PCA\\
\\
\begin{center}
   \begin{tabular}{|ll}
     Given & a matrix $A \in \R^n \otimes \R^p$  \\
     & a rank $r < p$ \\
     \\
     find & a matrix $A_r \in \R^n \otimes \R^p$ of rank $r$ \\
     \\
     such that & $\|A-A_r\|$ is minimum
   \end{tabular}
\end{center}
to set the problem of PCA with metrics on rows and columns as\\
\begin{center}
 \shadowbox{
  \begin{tabular}{lll}
      Given & a matrix & $A \in \R^{n \times p}$ \\
      & an inner product in $\R^n$ defined by & $N \in \R^{n \times n}$ \\
      & an inner product in $\R^p$ defined by & $P \in \R^{p \times p}$ \\  
      & with & $N=M^2, \quad P=Q^2$ \\
      &a rank & $0 < r < p$ \\ 
      \\
     \\
      Define & the inner product $T$ on $\R^{n \times p}$ & $T=Z^2$\\
      & associated to $(N,P)$& $ZA=MAQ$ \\
      \\
      Find & a matrix & $A_r \in \R^{n \times p}$\\ 
      with & & $\rank A_r=r$ \\
      \\
      such that & & $ \|A-A_r\|_\t \quad \mbox{minimal}$ \\
  \end{tabular}
 }
\end{center}

\subsection{Solving the problem}

We first give a direct solution, without using the associated isometry.

\nS A matrix of rank $r$ can be written as
\[
 A_r = \sum_{i=1}^r y_i \otimes v_i
\]
with $(v_i)_i$ being an orthonormal family for the inner product in $\R^p$ induced by $P$
\[
\langle v_i \, , \, v_j \rangle_\p = \delta_i^j
\]
Then
\[
 \|A-A_r\| = \left\| A - \sum_i y_i \otimes v_i \right\|
\]

\nS Let us select the inner product $T = N \otimes P$ on $\R^{n \times p}$ as
\[
 \langle A,B\rangle_{\n ,\p} = \langle MAQ \, , \, MBQ \rangle
\]
Then
\[
\|A\|_{\n,\p} = \|MAQ\|
\]
We recall that $M(y \otimes v)Q=My \otimes Q^\t v= My \otimes Qv$ as $Q$ is symmetric. Then
\begin{equation}\label{eq:pca_met:sol1}
    \begin{array}{lcl}
        \|A-A_r\|_{\n,\p} &=& \displaystyle \left\| M\left(A_r - \sum_i y_i \otimes v_i\right)Q \right\|\\
        \\
        &=& \displaystyle \left\| MAQ - \sum_i My_i \otimes Qv_i \right\|
    \end{array}
\end{equation}

\nS Let us denote
\begin{equation}
 \begin{cases}
    My_i &= w_i \\
    Qv_i &= x_i
 \end{cases}
\end{equation}
We have $\|x_i\|=\|Qv_i\|=1$ as $\|v_i\|_\p=1$ and similarily $\langle x_i,x_j\rangle = \langle Qv_i,Qv_j\rangle = \delta_i^j$. Then, the $(x_i)_i$ are an orthonormal family for the canonical inner product. The problem can be formulated as \\
\begin{center}
 \begin{tabular}{|ll}
  find & $(x_i)_i, \quad (w_i)_i$ \\
  with & $(x_i)_i$ an orthonormal family\\
  such that & $\displaystyle \left\| MAQ - \sum_i w_i \otimes x_i \right\| \quad  \mbox{minimal}$
 \end{tabular}
\end{center}
Then, $\{(w_i)_i \, , \, (x_i)_i\}$ are the solution of the PCA of $MAQ$. Th components $(y_i)_i$ and new basis vector $(v_i)_i$ of the PCA with metrics can be recovered simply by
\begin{equation}
 \begin{cases}
   v_i &= Q^{-1}x_i \\
   y_i &= M^{-1}w_i
 \end{cases}
\end{equation}
Hence the algorithm:\\
\\
\begin{algorithm}[H]
\begin{algorithmic}[1]
\STATE \textbf{input} $A \in \R^{n \times p}$ ; $M \in \R^{p \times p}$, SDP ; $Q \in \R^{n \times n}$, SDP
\STATE \textbf{compute} $B = MAQ$
\STATE \textbf{compute} $W,\Lambda,X = \textsc{pca\_core}(B)$
\STATE \textbf{compute} $Y=M^{-1}W$
\STATE \textbf{compute} $V=Q^{-1}X$
\RETURN $Y,\Lambda,V$
\end{algorithmic}
\caption{PCA of a matrix with double metrics: \textsc{pca\_met}$(A,M,Q)$}
\label{alg:pcamet}
\end{algorithm}

\paragraph{Remark:}  The metrics on $\R^n$ and $\R^p$ are given respectively by $N$ and $P$, which are symmetric, definite and positive (SDP). The matrices involved in this algorithm are respectively $M$ and $Q$, with $M=N^{1/2}$ and $Q=P^{1/2}$. They can be computed from a SVD of respectively $N$ and $P$. As $N$ is symmetric, its SVD reads
\[
N=U\Sigma U^\t
\]
Then
\[
M = U \Sigma^{1/2}U^\t
\]
Indeed, $M^2=(U \Sigma^{1/2}U^\t)(U \Sigma^{1/2}U^\t)=U \Sigma^{1/2}U^\t U \Sigma^{1/2}U^\t=U \Sigma U^\t=N$.

\subsection{Isometry}

This result can be derived without calculation by transportation of \PCA by isometry. Let $\R^n$ (resp. $\R^p$) be embedded with a Euclidean structure defined by SDP matrix $N=M^2$ (resp. $P=Q^2$). Then, the maps
\[
 \begin{CD}
  (\R^n, N) @>>> (\R^n,\I_n) \\
  y @>>> My
 \end{CD}
\]
and
\[
 \begin{CD}
  (\R^p, P) @>>> (\R^p,\I_p) \\
  v @>>> Qv
 \end{CD}
\]
are isometries, as
\[
\left\{ 
 \begin{array}{lcl}
  \langle Qv \, , \, Qv' \rangle &=& \langle v,v' \rangle_\p \\
  \langle My \, , \, My' \rangle &=& \langle y,y' \rangle_\n
 \end{array}
 \right.
\]
This induces an isometry on $\R^n \otimes \R^p$ by
\[
 \begin{CD}
  \psi \: : \: (\R^n \otimes \R^p, N \otimes P) @>>> (\R^n \otimes \R^p, \I_n \otimes \I_p) \\
  A @>>> MAQ
 \end{CD}
\]
as
\[
   \langle MAQ \, , \, MBQ \rangle = \langle A,B\rangle_{N \otimes P} 
\]
\PCA of $A$ with metric $P$ on $\R^p$ and $N$ on $\R^n$ is finding the best rank $r$ approximation of a matrix $A$, i.e. finding $(y_j \otimes v_j)_{1 \leq j \leq k}$ such that
\[
 \Delta = \left\| A - \sum_{j=1}^r y_j \otimes v_j \right\|_{\n \otimes \p}
\]
is minimum. Then, 
\[
 \Delta_\psi = \left\| \psi(A) - \psi\left(\sum_j y_j \otimes v_j\right) \right\|
\]
is minimum ($\Delta_\psi=\Delta$ because $\psi$ is an isometry). We have
\[
 \begin{aligned}
   \psi(A) - \psi\left(\sum_j y_j \otimes v_j\right) &= MAQ - M\left(\sum_j y_j \otimes v_j\right)Q \\
   &= MAQ - \sum_j My_j \otimes Qv_j
 \end{aligned}
\]
Then, $\sum_{j \leq r} My_j \otimes Qv_j$ with $\|Qv_j\|=1 \: \forall \: j$ is the best rank $r$ approximation of $MAQ$ which can be solved by a PCA of $MAQ$. Let $\sum_j w_j \otimes x_j$ be the best rank $r$ approximation of $MAQ$. Then, by applying isometry $\psi^{-1}$, $\sum_jM^{-1}w_j \otimes Q^{-1}x_j$ is the best rank $k$ approximation of $A$ for metric defined by $N \otimes P$, and the solution is $(Y,V,\Lambda)$ with $Y=M^{-1}W$ and $V=Q^{-1}X$.

%
\subsection{Interpretation and plotting}
%

\nS A common situation is when metrics are given as weights on the columns only. Then, $M=\I_n$ and 
\[
 MAQ = AQ
\]
Hence
\begin{equation}
 \begin{array}{lcl}
   B^\t B&=& (MAQ)^\t (MAQ) \\
   &=& (AQ)^\t (AQ) \\
   &=& Q^\t A^\t A Q
 \end{array}
\end{equation}
Similarly, if the metrics are weights on the rows only, $Q=\I_p$ and 
\[
 MAQ=MA
\]
Hence
\begin{equation}
 \begin{array}{lcl}
   B^\t B&=& (MAQ)^\t (MAQ) \\
   &=& (MA)^\t (MA) \\
   &=& A^\t MMA \\
   &=& A^\t NA
 \end{array}
\end{equation}
If metrics are given on both rows and columns, we have
\begin{equation}
 B^\t B = QA^\t NAQ
\end{equation}

\nT{A remark about the calculation} Principal components $(y_i)_i$ and principal axis $(v_i)$ are solution of
\begin{equation}
 \left\{
    \begin{array}{ccl}
     B &=& MAQ \\
     B^\t Bw_i &=& \lambda_iw_i  \\
     x_i &=& Bw_i \\
     v_i &=& Q^{-1}w_i \\
     y_i &=& M^{-1}x_i
    \end{array}
 \right.
\end{equation}
with
\[
 x_i,y_i \in \R^n, \qquad v_i,w_i \in \R^p
\]
We have
\begin{equation}
 \left. 
    \begin{array}{lcl}
     B^\t B &=& QA^\t NAQ \\
     B^\t Bw_i &=& \lambda_iw_i \\
     w_i &=& Qv_i
    \end{array}
    \right\}    
    \quad \Longrightarrow \quad 
    QA^\t NAPv_i =\lambda_iQv_i
\end{equation}
and, as $Q$ is invertible
\begin{equation}\label{pcamet:eq:res}
 A^\t NAPv_i = \lambda_iv_i
\end{equation}
This might lead to a way of computing the principal axis $(v_i)_i$ directly without computing the $(w_i)_i$ before. However, the matrix $B^\t B$ is symmetric, whereas matrix $A^\t NAP$ is not. It is known that numerical computation of eigenvectors and eigenvalues of symmetric matrices is more accurate and robust than of non-symmetric matrices. Hence, it is recommended to compute first $(w_i)_i$ as solutions of $B^\t Bw_i=\lambda_iw_i$ and then the principal axis as $v_i=Q^{-1}w_i$.

\nT{Centering the cloud} As for \PCA, it is advised to center the cloud before analysing it when there are some weights on rows. Let us recall that if each point $a_i \in \R^p$ is given a weight $w_i$, the barycenter $g \in \R^p$ is given by
\begin{equation}
 \left(\sum_iw_i\right)g = \sum_iw_ia_i
\end{equation}
or
\begin{equation}
 g = \frac{1}{w}\sum_iw_ia_i, \qquad w = \sum_iw_i
\end{equation}
The centered cloud is the cloud with points
\begin{equation}
 \overline{a}_i = a_i-g
\end{equation}
One checks that
\[
 \begin{aligned}
  \sum_iw_i\overline{a}_i &= \sum_iw_ia_i - \left(\sum_iw_i\right)g \\
  &= wg - wg \\
  &= \zeros
 \end{aligned}
\]

\nT{Geometric approach: attached point cloud} Let $\A$ be the point cloud of $n$ points in $\R^p$ attached to matrix $A$. Distances between points do not reflect the distances induced by the inner products $(M,Q$). Let us denote by $\B$ the point cloud in $\R^p$ attached to matrix $B=MAQ$. Points $b_i,b_k$ have as coordinates respectively the rows $i$ and $k$ of $B$. If $M,Q$ are diagonal matrices with weights $(\sqrt{\nu_i},\sqrt{\pi_j})$ respectively, then
\[
 MAQ = [\sqrt{\nu_i\pi_j}\, \alpha_{ij}]_{i,j}
\]
and, in $\R^p$
\begin{equation}
 \begin{array}{lcl}
   d^2(b_i,b_k) &=& \displaystyle \sum_j (\sqrt{\nu_i\pi_j}\, \alpha_{ij} - \sqrt{\nu_k\pi_j}\, \alpha_{kj})^2 \\
   &=& \displaystyle \sum_j \pi_j\, (\sqrt{\nu_i}\, \alpha_{ij} - \sqrt{\nu_k}\, \alpha_{kj})^2 
 \end{array}
\end{equation}
This is the distance between points of the point cloud in $\R^p$ attached to matrix $MA$ with the inner product in $\R^p$ defined by weight matrix $P$. So, \PCA of matrix $A \in \R^{n \times p}$ with inner product defined by $N$ in $\R^n$ and $P$ in $\R^p$ is \PCA of point cloud $\A_\m$ in $\R^p$ attached to matrix $MA$ with inner product defined by $P$ in $\R^p$ for computing distances. This will be useful for Correspondance Analysis (see section \ref{sec:coa}).

\nT{Scaled \PCA} A straightforward and standard application of \PCA with metrics is scaled PCA. Let $A \in \R^{n \times p}$ be a columnwise centered matrix, i.e.
\begin{equation}
 \sum_ia_i=0
\end{equation}
The variances (or norms) of columns of $A$ can vary significantly. In such a case, the variance/covariance matrix $\Sigma=A^\t A$ can be dominated by rows and columns corresponding to the variable with largest variance. Scaled PCA is clipping this uninteresting result off, by giving equal weights to each variable. The technical trick is to equalize variances between columns, by dividing each column $j$ by its standard deviation (or norm). If $a_{\bullet j} \in \R^n$ is column $j$ of $A$, this reads
\begin{equation}
 \begin{CD}
  a_{\bullet j} @>>> \displaystyle a'_{\bullet j} = \frac{a_{\bullet j}}{\|a_{\bullet j}\|}
 \end{CD}
\end{equation}
Hence
\begin{equation}
 \begin{array}{lcl}
   \|a'_{\bullet j}\| = 1
 \end{array}
\end{equation}
This can be read as a \PCA with inner product $N=\I_n$ in $\R^n$ and $P = \diag \left(\frac{1}{\|a_{\bullet j}\|}\right)$ in $\R^p$. So $MAQ=A'$, and \PCA of $A'$ is run.

\notes The problem (and solution) of PCA with weights on rows and columns can be found in \cite{Rao1964} or \cite{Greenacre1984}. It is presented in \cite[sect. 14.2]{Jolliffe2002}. The algebraic approach with generalization to metrics in Euclidean spaces has been proposed as a general method with the notion of duality diagram in \cite{Cailliez1979} which has been at the root of many works (see \cite{Pages1979}). The formalism presented here can be found in \cite{Franc1992}.

%
%

%
\subsection{Analysis of a matrix with metrics and weights: a geometric approach}\label{sec:pcamet:wP}
%

The geometric school of linear dimension reduction has developed a framework to analyse a matrix $A \in \R^{n \times p}$ with a metric defined by $P$ in $\R^p$ and weights $(w_i)_i$ on the individual. The point cloud is made by rows $(a_{i*})_i \in \R^p$ of $A$. It is presented here and will be very useful for understanding the geometric approach of \CoA. 

\nT{Setting the problem} Rows of $A$ (points of the point cloud) are in $\R^p$. We set the problem for the best projection of the point cloud on a one-dimensional space in $\R^p$ spanned by a vector $v$ with $\|v\|_\p=1$. The aim is to compute $v$. The projection of $a_{i*}$ on $\R v$ is given by
\[
\PC_v a_{i*} = \langle Pa_{i*},v\rangle v
\]
and its norm is
\[ 
\|\PC_v a_{i*}\|_\p^2 = \langle Pa_{i*},v\rangle^2
\]
Let
\[
\IC = \sum_i \; \left\|\PC_v a_{i*}\right\|_\p^2 = \sum_i \;\langle Pa_{i*},v\rangle^2
\]
be the inertia of the point cloud attached to $A$ for inner product $P$. So, first step without weights but with metrics defined by $P$ is to find $v$ with $\|v\|_\p=1$ such that $\IC$ is maximal. Final step is to introduce the weights, and define
\[
\IC = \sum_i \;w_i \; \langle Pa_{i*},v\rangle^2
\]
as the inertia of the point cloud with inner product defined by $P$ and weights on rows defined by $w$. Then, the geometrical approach can be set as\\
\begin{center}
    \ovalbox{
        \begin{tabular}{ll}
          given & a matrix $A \in \R^{n \times p}$ \\
          with & row $i$ denoted $a_{i*} \in \R^p$ \\
          an inner product & $P$ in $\R^p$ \\
          a set of row weights & $w \in \R^n$ \\
          &\\
          find & a vector $v \in \R^p$ \\
          with & $\|v\|_\p=1$ \\
          & \\
          such that & $ \displaystyle \IC = \sum_i \;w_i \; \left\| \PC_v a_{i*}\right\|^2$ is maximal \\
          where & $\PC_v$ is the projection on $\R v$ in $\R^p$
        \end{tabular}
    }
\end{center}

\nT{Solving the problem} Let us note first that
\begin{equation}
    \left\| \PC_v a_{i*}\right\|^2 = \langle  Pa_{i*},v\rangle^2
\end{equation}
Then
\begin{equation}
\IC = \sum_i \;w_i \; \langle Pa_{i*},v\rangle^2 = \sum_i \langle \sqrt{w_i}\, P\,a_{i*} \; , \; v\rangle^2
\end{equation}
Let us observe that
\begin{itemize}[label=$\rightarrow$]
    \item the vectors $Pa_{i*} \in \R^p$ are the row vectors of matrix $AP$ (indeed, $P$ is symmetric by definition)
    \item the terms $\sqrt{w_i}Pa_{i*}$ are the rows of matrix $D_w^{1/2}AP$ if $D_w$ is the $n \times n$ diagonal matrix with $w$ in the diagonal
\end{itemize}
Then $\IC$ can be rewritten as
\begin{equation}
    \IC = \|D_w^{1/2}APv\|^2
\end{equation}
So, the problem can be restated as
\begin{center}
    \begin{tabular}{|ll}
        given & $A,P,w$ as above,\\
        find & $v \in \R^p$ \\
        with & $\|v\|_\p^2 = 1$ \\
        such that & $\IC = \|D_w^{1/2}APv\|^2$ is maximal
    \end{tabular}
\end{center}
To solve this with Lagrange multipliers (an optimum with constraints), let us denote 
\begin{equation}
    H=D_w^{1/2}AP
\end{equation}
Then, $\IC = \|Hv\|^2$, and 
\begin{equation}
    \frac{\partial \IC}{\partial v} = 2H^\t H, \qquad \frac{\partial \|v\|_\p^2}{\partial v} = 2Pv
\end{equation}
so
\begin{equation}
    H^\t H v = Pv
\end{equation}
This can be written
\begin{equation}
    PA^\t D_w^{1/2}D_w^{1/2}APv = PA^\t D_wAPv = \lambda Pv
\end{equation}
and $Pv$ is an eigenvector of $PA^\t D_wA$. As $P$ is invertible (it is a SDP matrix), this yields
\begin{equation}
    A^\t D_wAPv = \lambda v
\end{equation}
and $v$ is an eigenvector of $A^\t D_wAP$.

\notes This section is adapted from \cite[1.1.6]{LMP2000} where it is presented as a diversification of the general analysis (PCA). Similar approaches can be found in \cite{Lebart1977,LMF82}.

%
\subsection{\PCA with metrics and instrumental variables}
%

Those methods, \PCAmet and \PCAiv can be associated like  pieces of puzzle to build a chain of treatments.

\nS Let us recall (see section \ref{sec:pcaiv}) that \PCAiv is running the \PCA of $A$ with constraints on principal axis and components which must belong to subspaces of respectively $\R^p$ and $\R^n$:
\begin{equation}
 \begin{cases}
  y_j \in E \subset \R^n \\
  v_j \in F \subset \R^p
 \end{cases}
\end{equation}
If $U$ (resp. $V$) is an orthonormal matrix with a basis of $E$ (resp. $F$) as column vectors, this is done by building the projectors 
\begin{equation}
 \begin{CD}
  \R^n @>R=UU^\t>> E, \qquad \R^p @>S=VV^\t>>F
 \end{CD}
\end{equation}
and running the \PCA of $A'=RAS$, the projection of $A$ on $E \otimes F$.

\nS If $\R^n$ and $\R^p$ are endowed with metrics given by $N$ and $P$ respectively, running the \PCA of $A'$ with those metrics is running the \PCA of $MAQ$ with $M=N^{1/2}$ and $Q=P^{1/2}$. However, the projectors $R$ and $S$ depend on the metrics $N$ and $P$. 

\nS Let us write the projector on $F \subset \R^p$ with the inner product defined by $P$ first. Let $x \in \R^p$ and $v \in F$ with $\|v\|_\p=1$. The projection $x'$ of $x$ on $F=\R v$ is given by
\begin{equation}
  \begin{array}{lcl}
    x' &=& \langle x,v\rangle_\p v \\
    &=& \langle x,Pv\rangle v \\
    &=& \langle Pv,x\rangle v \\
    &=& (v \otimes Pv).x
  \end{array}
\end{equation}
Hence, the projector on $\R v$ with the inner product defined by $P$ is $v \otimes Pv$. If $F = \span (v_1,\ldots,v_r)$, we have
\begin{equation}
 \begin{array}{lcl}
   x' &=& \displaystyle \sum_a \langle x,Pv_a\rangle v_a \\
   &=& \displaystyle \sum_a (v_a \otimes Pv_a)x 
 \end{array}
\end{equation}
and the projector is
\begin{equation}
  \begin{array}{lcl}
    S &=& \displaystyle \sum_a v_a \otimes Pv_a \\
    &=& V(PV)^\t \\
    &=& VV^\t P \\
    &=& PVV^\t
  \end{array}
\end{equation}
if $V \in \R^{p \times r}$ is the matrix with $v_a$ in column $a$. The last equality comes from the observation that $P$ and $VV^\t$ are symmetric, hence $(VV^\t)P$ is symmetric and $VV^\t P=(VV^\t P)^\t=PVV^\t$.

\nS Similarly, we have
\begin{equation}
 R = UU^\t N
\end{equation}
and 
\begin{equation}
 \begin{array}{lcl}
   A' &=& RAS \\
   &=& UU^\t NAPVV^\t
 \end{array}
\end{equation}

\nS Let us now write the \PCA of $A'$ with inner products defined by $N$ on $\R^n$ and $P$ on $\R^p$. It is the \PCA of
\[
 A'' = MA'Q, \qquad \mbox{with} \quad N=M^2, \quad P=Q^2
\]

%
\section{Correspondence Analysis}\label{sec:coa}
%

\subsection*{Notations}

Some notations for Correspondence Analysis are standard and specific to this method. They are given here and explained along the chapter.\\
\\
\begin{center}
\ovalbox{
\begin{tabular}{ccl}
symbol & in space & meaning \\
\hline
$A$ & $\R^{I \times J}$ & matrix of frequencies $(T/n_{++})$ \\
$\alpha_{ij}$ & $[0,1] \subset \R$ & general term of $A$ \\ 
$c$ & $\R^J$ & vector of marginals of $A$ by columns: $(c_1, \ldots, c_J)$ \\
$c_j$ & $\R$ & marginal of column  $j$ of $A$: $(\sum_i\alpha_{ij})$\\
$D_c$ & $\R^{J \times J}$ & diagonal matrix with elements $c_j$ \\
$D_r$ & $\R^{I \times I}$ & diagonal matrix with elements $r_i$ \\
$i$ & $\N$ & indices of rows \\
$I$ & $\N$ & number of rows of $T$ \\
$j$ & $\N$ & indices of columns \\
$J$ & $\N$ & number of columns of $T$ \\
$n_{ij}$ & $\N$ & number of item in class $i$ (row) and $j$ (column) \\
$n_{i+}$ & $\N$ & sum of terms in row $i$ of $T$\\
$n_{+j}$ & $\N$ & sum of terms in column $j$ of $T$\\
$n_{++}$ & $\N$ & sum of terms in $T$ \\
$r$ & $\R^I$ & vector of marginals of $A$ by rows: $(r_1, \ldots,r_I)$ \\
$r_i$ & $\R$ & marginal of row $i$ of $A$: $(\sum_j\alpha_{ij})$\\
$T$ & $\R^{I \times J}$ & contingency table \\
\end{tabular}
}

\end{center}

\vspace*{5mm}

\noindent A remarkable application of the \PCA with weights on rows and columns is the development of Correspondence Analysis as the 
analysis of a contingency table with metrics associated to its margins.

\nB Let us adopt here some standard notations for contingency tables. A contingency table $T$ is a table of counts of $n$ items allocated to categories of two variables. Indices of the values of the first variable are denoted $i$, and of the second variable $j$. The value $n_{ij}$ in row $i$ and column $j$ of $T$ is the number of items in category $i$ for the first variable and $j$ for the second. It is standard to denote that $i \in \llbracket 1,I\rrbracket$ and $j \in \llbracket 1,J \rrbracket$. Then
\[
 T \in \R^{I \times J} \simeq \R^I \otimes \R^J
\]
It is standard to denote
\[
n_{i+}= \sum_jn_{ij}, \quad n_{+j} = \sum_{i}n_{ij}, \qquad n_{++}= \sum_{i,j}n_{ij} =\sum_in_{i+}=\sum_jn_{+j}
\]
%
\subsection{Link with $\chi^2$ distance}
%

We first establish a link between the norm of a contingency table with metrics associated to margins on rows and columns on one hand and the $\chi^2$ of the table on the other.

\nS Let $T$ be a contingency table of two discrete variables observed on $n$ individuals, with $n_{ij}$ being the number of individuals with observation $i$ for first variable and $j$ for the second. Then
$T \in \R^{I \times J}$ if first variable has $I$ values and second $J$.\\
\\
Let us denote
\begin{equation}
 A = \frac{T}{n_{++}}
\end{equation}
The general term in $A$ is denoted $\alpha_{ij}$ and we have 
\begin{equation}
 A \in \R^{I \times J} \simeq \R^I \otimes \R^J
\end{equation}
Let us denote respectively by $r \in \R^I$ and $c \in \R^J$ the marginal sums of $A$ on rows and columns
\begin{equation}
 \left\{
   \begin{array}{lclcl}
    r_i &=& \displaystyle \sum_j \alpha_{ij} &=& \alpha_{i+}\\
    c_j &=& \displaystyle \sum_i \alpha_{ij} &=& \alpha_{+j}
   \end{array}
 \right.
\end{equation}
Let us denote by $D_r$ and $D_c$ the square diagonal matrices with diagonal respectively $r$ and $c$:
\begin{equation}
 D_r = \diag r, \qquad D_c = \diag c
\end{equation}

\nS In case of independence between both variables, the expectation for $A$ is
\begin{equation}
 \widetilde{A} = r \otimes c
\end{equation}
Indeed, we have, ignoring the value of the other variable, if $X_r$ is the ransom variable for rows and $X_c$ for columns
\[
 P(X_r=i) = r_i, \qquad P(X_c=j)=c_j
\]
Then, in case of independence
\[
 P(X_r=i \, ; \, X_c=j)=r_ic_j
\]

\nS Then, a first observation is that
\begin{equation}
 \chi^2(A) = \|A - \widetilde{A}\|_{D_r^{-1} \otimes D_c^{-1}}
\end{equation}
\begin{proof}
Indeed, we have
\begin{equation}
  \begin{array}{lcl}
    \chi^2(A) &=& \displaystyle \sum_{i,j} \frac{(\alpha_{ij}-r_ic_j)^2}{r_ic_j} \\
    &=& \displaystyle \sum_{i,j} \left(\frac{\alpha_{ij}-r_ic_j}{\sqrt{r_ic_j}}\right)^2 \\
    &=& \displaystyle \sum_{i,j} \left(\frac{1}{\sqrt{r_i}}\, (\alpha_{ij}-r_ic_j)\, \frac{1}{\sqrt{c_j}}\right)^2 \\
    &=& \displaystyle \left\|D_r^{-1/2}(A - r \otimes c)D_c^{-1/2}\right\|^2 \\
    &=& \|A - r \otimes c\|^2_{D_r^{-1} \otimes D_c^{-1}}
  \end{array}
\end{equation}
\end{proof}

%
\subsection{Description of the method}
%

Then, Correspondence Analysis is a partition of the variance $\|A - r \otimes c\|_{D_r^{-1} \otimes D_c^{-1}}$ concentrated on the first axis. It is henceforth a \PCA of $A - r \otimes c$ with metrics defined by $D_r^{-1}$ on rows and $D_c^{-1}$ on columns, i.e. with weights $1/r_i$ on row $i$  and $1/c_j$ on column $j$.\\
\\
\begin{center}
 \ovalbox{
  \begin{tabular}{ll}
      given & $T = (n_{ij})_{i,j}$ (a contingency table)\\
      & \\
      compute & $n_{++} = \sum_{i,j}n_{ij}$ \\
      & $A = \frac{T}{n_{++}}$ \\
      & $r_i = \sum_j\alpha_{ij}$ \\
      & $c_j = \sum_i\alpha_{ij}$ \\
      \\
      run & \texttt{\PCAmet} \\
      on & $A - r \otimes c$ \\
      & \\
      with diagonal metrics & $1/r_i$ on row $i$ \\
      & $1/c_j$ on column $j$
  \end{tabular}
 }
\end{center}
We have
\[
 \left\{ 
    \begin{aligned}
    M & = \diag 1/\sqrt{r_i} \\
    Q & = \diag 1/\sqrt{c_j} \\  
    \end{aligned}
 \right.
\]
Then
\begin{equation}
 \begin{aligned}
   A_{\textsc{m},\textsc{q}} &= M(A-r\otimes c) Q \\
   &= \left[ \frac{\alpha_{ij} - r_ic_j}{\sqrt{r_ic_j}} \right]_{i,j}
 \end{aligned}
\end{equation}
This yields the following algorithm:\\
\\
\begin{algorithm}[H]
\begin{algorithmic}[1]
\STATE \textbf{input} $T \in \R^{I \times J}$, a contingency table
\STATE \textbf{compute} $A = T/T_{++}, \quad \mbox{with} \quad T_{++} = \sum_{i,j}T_{ij}$ \\
\STATE \textbf{compute} $r_i = \sum_j\alpha_{ij}, \quad D_r = \diag r$
\STATE \textbf{compute} $c_j = \sum_i\alpha_{ij}, \quad D_c = \diag c$
\STATE \textbf{compute} $M = \diag (1/\sqrt{r_i})$
\STATE \textbf{compute} $Q = \diag (1/\sqrt{c_j})$
\STATE \textbf{compute} $A_{\textsc{m},\textsc{q}} = M(A-r\otimes c)Q = \left[\frac{\alpha_{ij}-r_ic_j}{\sqrt{r_ic_j}}\right]_{i,j}$
\STATE \textbf{compute} $Z,X,\Lambda = \textsc{pca\_core}(A_{\textsc{m},\textsc{q}} )$
\STATE \textbf{compute} $Y_r=M^{-1}Z$
\STATE \textbf{compute} $Y_c=Q^{-1}X$
\RETURN $Y_r,Y_c, \Lambda$
\end{algorithmic}
\caption{Correspondence Analysis of a contingency table: \textsc{coa}$(T)$}
\label{alg:coa}
\end{algorithm}

%
\subsection{\CoA and geometry of point clouds}
%

Let us consider the point cloud $\X$ of $I$ points in $\R^\j$ of coordinates
\begin{equation}
 X_{ij} = \frac{\alpha_{ij}}{\sqrt{r_ic_j}}
\end{equation}
in a Euclidean space with standard inner product (i.e. $\I_\j$). Then, if $x_i,x_{i'}$ are two points in $\X$, we have
\begin{equation}
\begin{aligned}
 d(x_i,x_{i'})^2 &= \sum_j\left(\frac{\alpha_{ij}}{\sqrt{r_ic_j}} - \frac{\alpha_{i'j}}{\sqrt{r_{i'}c_j}}\right)^2 \\
 &= \sum_j\frac{1}{c_j}\left(\frac{\alpha_{ij}}{\sqrt{r_i}} - \frac{\alpha_{i'j}}{\sqrt{r_{i'}}}\right)^2 
\end{aligned}
\end{equation}
It is standard to say that point cloud $\X=(x_i)_i$ is embedded with weights $1/c_j$ for column ( $=$ variable) $j$ and metrics defined by diagonal matrix of term $1/r_i$ in $\R^I$.

%
\subsection{Classical presentation: geometric approach}
%
There is a classical presentation of Correspondance Analysis as an analysis of two point clouds associated to a contingency table, one for the rows and one for the columns (which play a symmetric role), each with weights and metrics. This is the geometric approach. It emphasizes that two point clouds, and not simply one, can be built: one for rows, and one for columns, and \CoA can be seen as a simultaneous analysis of both.

\nS Let us recall that if $T$ is a contingency table, $F$ (the matrix of frequencies) is built from $T$ by dividing it by the sum of all its elements:
\[
 \begin{CD}
   T = (n_{ij})_{i,j} @>>> \displaystyle n_{++} = \sum_{i,j}n_{ij} @>>> \displaystyle F \: : \: f_{ij} = \frac{n_{ij}}{n_{++}}
 \end{CD}
\]
To comply with standard notations in classical textbooks (see below), we denote by $F$, and not $A$, the matrix of frequencies. We will denote
\begin{center}
   \begin{tabular}{lcccl}
     for row & & & & for columns \\
     \hline
     \\
     $\displaystyle f_{i+} = \sum_jf_{ij}$ & $\qquad$ & and  & $\qquad$  & $\displaystyle f_{+j} = \sum_if_{ij}$ \\ 
     \\
     $ \displaystyle f_{i*} = (f_{i1}, \ldots,f_{iJ}) \in \R^J$ & $\qquad$  & and & $\qquad$ & $\displaystyle f_{*j} = (f_{1j}, \ldots,f_{Ij}) \in \R^I$ \\
     \\
     $r = (f_{1*},\ldots,f_{I*}) \in R^I$ & $\qquad$  & and & $\qquad$ & $c = (f_{*1}, \ldots,f_{*J}) \in \R^J$  
   \end{tabular}
\end{center}

\nS Two point clouds are classically attached to $A$
\begin{itemize}[label=$\rightarrow$]
 \item a cloud of row profiles, as $I$ points $(x_i)_i$ in $\R^\j$, with point $x_i$ of coordinates
 \begin{equation}
  x_i = \left[\frac{f_{ij}}{f_{i+}}\right]_j \in \R^\j
  \end{equation}
 \item a cloud of column profiles, as $J$ points $y_j$ in $\R^\i$, with point $y_j$ of coordinates
 \begin{equation}
  y_j = \left[\frac{f_{ij}}{f_{+j}}\right]_i \in \R^\i
  \end{equation}
\end{itemize}

\nS Then, metrics with diagonal matrices respectively $D_c^{-1}$ for $R^\j$ and $D_r^{-1}$ for $\R^\i$ are selected, with $D_c$ being the $J \times J$ diagonal matrix of term $1/f_{+j}$ and $D_r$ the $I \times I$ diagonal matrix of term $f_{i+}$. Hence, distances between points $x_i$ and $x_k$ in $\R^\j$ are computed as
\begin{equation}
 d^2(x_i,x_k) = \sum_j \frac{1}{f_{+j}}\left(\frac{f_{ij}}{f_{i+}}- \frac{f_{kj}}{f_{k+}}\right)^2
\end{equation}
and between points $y_j$ and $y_{\ell}$ in $\R^\i$ as
\begin{equation}
 d^2(y_j,y_{\ell}) = \sum_i \frac{1}{f_{i+}}\left(\frac{f_{ij}}{f_{+j}}- \frac{f_{i\ell}}{f_{+\ell}}\right)^2
\end{equation}
These weights tend to give an equal importance to all modalities of a variable, whatever their size. It is analogous to weighting by the inverse of the variance in scaled \PCA. A further property often invoked is that, if two categories of a same variable have the same profile (say $i$ and $i'$ for which $\forall \j, \:\:f_{ij}/f_{i+}= f_{i'j}/f_{i'+}$), then it is logical to lump them together into a single category, and this must not modify the distances between row profiles. 

\nS  Let us recall that 
\[
r = \left(f_{1+},\ldots,f_{I+}\right) \in \R^I, \qquad \mbox{and} \qquad c = \left(f_{+1},\ldots,f_{+J}\right) \in \R^J, 
\]
so
\[
r_i = f_{i+} \qquad \mbox{and} \qquad c_j =f_{+j}
\]
It is then standard to define both point clouds as (see \cite[sect. 4.1]{Greenacre1984}, \cite[sect. IV.5.]{LMF82}, \cite[sect. 10.1]{Saporta1990},\cite[sect. 1.3.3]{LMP2000})\\
\begin{center}
 \ovalbox{
 \begin{tabular}{l|c|c}
   & R: row profiles in $\R^J$& C: column profiles in $\R^I$\\
   \hline 
   point cloud & $R = D_r^{-1}F$ & $C = D_c^{-1}F^\t$ \\
   metric & $D_c^{-1}$ & $D_r^{-1}$ \\
   weights & $r$ & $c$ \\   
 \end{tabular}
 }
\end{center}

\nB $R$ is the point cloud of row profiles of $F$, and $C$ of its columns profiles. It is easy to get lost in the indices, the rows, the columns, spaces $\R^I$ and $\R^J$. A row profile is in $\R^J$ and its coefficients are indexed by $j$; a column profile is in $\R^I$ and its coefficients are indexed by $i$. With coefficients, this yields
\begin{equation}
    R_{ij} = \frac{f_{ij}}{f_{i+}}, \qquad C_{ij} = \frac{f_{ij}}{f_{+j}}
\end{equation}

\nS The centroids $g^{(\r)}$ of $R$ and $g^{(\c)}$ of $C$ are respectively $c$ and $r$ (and not $r$ and $c$). Indeed, 
\begin{equation}
    g^{(\r)}_j = \sum_i f_{i+}\,\frac{f_{ij}}{f_{i+}} = \sum_if_{ij} = f_{+j} = c_j
\end{equation}
and
\begin{equation}
    g^{(\c)}_i = \sum_j f_{+j}\,\frac{f_{ij}}{f_{+j}} = \sum_jf_{ij} = f_{i+} = r_i
\end{equation}
Then, the inertia $I_\textsc{r}$ of centered point cloud $R$, i.e. of $R - \ones_I \otimes c$ with weights $r$ on rows and metric $D_c^{-1}$ in $\R^J$ is,  (explained step by step ...)
\begin{equation}
 \begin{array}{lclcl}
   I_\textsc{r}^2 &=& \displaystyle \sum_i r_i \left\|R_{i*} - c\right\|_{D_c^{-1}}^2 & \qquad & R_{i*} = (R_{i1},\ldots,R_{iJ})\\
   &=& \displaystyle \sum_i f_{i+} \left\| \frac{f_{i*}}{f_{i+}}-c\right\|_{D_c^{-1}}^2, &\qquad &  R_{i*} =  \frac{f_{i*}}{f_{i+}}, \: r_i=f_{i+}\\
   &=& \displaystyle \sum_i f_{i+} \left(\frac{1}{f_{+j}}\sum_j \left(\frac{f_{ij}}{f_{i+}}-f_{+j}\right)^2\right), &\qquad& c_j = f_{+j}, \:D_c^{-1} = \diag \left(\frac{1}{f_{+j}}\right)\\
   &=& \displaystyle \sum_i f_{i+} \left(\sum_j \left(\frac{f_{ij}}{f_{i+}\sqrt{f_{+j}}}-\sqrt{f_{+j}}\right)^2\right) &\qquad& \sqrt{f_{+j}} = \frac{f_{+j}}{\sqrt{f_{+j}}}\\
   &=& \displaystyle \sum_{i,j} \left(\sqrt{f_{i+}}\frac{f_{ij}}{f_{i+}\sqrt{f_{+j}}}-\sqrt{f_{i+}}\sqrt{f_{+j}}\right)^2\\
   &=& \displaystyle  \sum_{i,j} \left(\frac{f_{ij} - f_{i+}f_{+j}}{\sqrt{f_{i+}f_{+j}}}\right)^2 
 \end{array}
\end{equation}
where we recognize the $\chi^2$ norm of $F=(f_{ij})_{i,j}$, or the norm of $F$ with metrics defined by $D_c^{-1}$ in $\R^J$ (space of rows) and by $D_r^{-1}$ in $\R^I$ (space of columns). Its partition with concentration of the inertia on the first components is the \CoA of $F$.

\nB If $I_\textsc{c}^2$ is the inertia of centered point cloud $C$ in $\R^I$, i.e. of $C - \ones_p \otimes r$ with weights $c$ on rows and metric $D_r^{-1}$ in $\R^I$, a similar calculation yields
\begin{equation}
 I_\textsc{c}^2 = \sum_{i,j} \left(\frac{f_{ij} - f_{i+}f_{+j}}{\sqrt{f_{i+}f_{+j}}}\right)^2 = I_\textsc{r}^2 
\end{equation}
Both inertia are equal, and the analysis of both point clouds is one geometric guise of the \CoA of $F$.

\nS To see this, one can use the developments presented in section \ref{sec:pcamet:wP}. Let us first present a translation of the notations between section \ref{sec:pcamet:wP} and here:\\
\begin{center}
    \begin{tabular}{ccccc|c}
     \PCAmet with weight & & & & \CoA on $R$ & \CoA on $C$\\
     \hline
     $A$ & & $\longleftrightarrow$ & & $R=D_r^{-1}F$ & $C = D_c^{-1}F^\t$ \\
     $P$ & & $\longleftrightarrow$ & & $D_c^{-1}$ & $D_r^{-1}$ \\
     $w$ & & $\longleftrightarrow$ & & $r$ & $c$ \\
    \end{tabular}
\end{center}
Let us recall that the first principal axis $u$ is solution of 
\begin{equation}
    H^\t H u = Pu
\end{equation}
with
\begin{equation}
    H = D_w^{1/2}AP
\end{equation}
so
\begin{equation}
    H^\t Hu = PA^\t D_w APu=\lambda Pu
\end{equation}
and, as $P$ is ivertible as a SDP matrix, $u$ is solution of
\begin{equation}
    A^\t D_w APu=\lambda u
\end{equation}
We then have, for \CoA on row profiles
\begin{itemize}[label=$\rightarrow$]
    \item $H=D_r^{1/2}.D_r^{-1}F.D_c^{-1}= D_r^{-1/2}FD_c^{-1}$
    \item $H^\t H = (D_c^{-1}F^\t D_r^{-1/2}).(D_r^{-1/2}FD_c^{-1}) = D_c^{-1}F^\t D_r^{-1}FD_c^{-1}$
\end{itemize}
and, after simplification by $D_c^{-1}$, $u$ is solution of
\begin{equation}\label{eq:coa:1}
    F^\t D_r^{-1}FD_c^{-1}u = \lambda u
\end{equation}
and for \CoA on column profiles
\begin{itemize}[label=$\rightarrow$]
    \item $H=D_c^{1/2}.D_c^{-1}F^\t .D_r^{-1}= D_c^{-1/2}F^\t D_r^{-1}$
    \item $H^\t H = (D_r^{-1}F D_c^{-1/2}).(D_c^{-1/2}F^\t D_r^{-1}) = D_r^{-1}F D_c^{-1}F^\t D_r^{-1}$
\end{itemize}
and, after simplification by $D_r^{-1}$, first axis $v$ is solution of
\begin{equation}\label{eq:coa:2}
    F D_c^{-1}F^\t D_r^{-1}v = \lambda v
\end{equation}
(see \cite[p. 83]{LMP2000}, with, here again, a translation of notations). 

\nS If we multiply equation (\ref{eq:coa:1}) on the left by $D_c^{-1}$ and set $u'=D_c^{-1}u$, we have 
\begin{equation}
D_c^{-1}F^\t D_r^{-1}FD_c^{-1}u = \lambda D_c^{-1}u
\end{equation}
or
\begin{equation}
  CRu'=\lambda u'  
\end{equation}
Similarily, multiplying equation (\ref{eq:coa:2}) left by $D_r^{-1}$ and setting $v'=D_r^{-1}v$ yields
\begin{equation}
    RCv'=\lambda v'
\end{equation}
This yields
\begin{equation}
    \left\{ 
      \begin{array}{lcl}
        v' &=& Ru' \\
        u' &=& Cv'
      \end{array}
    \right.
\end{equation}
or
\begin{equation}
    \left\{ 
      \begin{array}{lcl}
        v &=& D_r^{-1}RD_c^{-1}u \\
        u &=& D_c^{-1}CD_r^{-1}v
      \end{array}
    \right.
\end{equation}

\notes Correspondence Analysis has a long history, and has been object of long debates, renaming and rediscoveries between different schools. Correspondence Analysis has been proposed first by Hirshfeld in 1935 (Hirschfeld, M. O. - 1935 - A connection between correlation and contingency - \emph{Proc. Camb. Phil. Soc.}, \textbf{31:}520-524). It has been rediscovered by Guttman in 1959 (Guttman, L. - 1959 - Metricizing rank ordered and unordered data for a linear factor analysis. \emph{Sankhy\={a}}, \textbf{21:}257-268). The link between CA and reciprocal averaging has been presented in \cite{Hill1974}. Correspondence Analysis has been rediscovered and studied independently by several researchers, as J.-P. Benzecri, in 1962 in the context of mathematical linguistics inspired by the works of Chomsky; J. de Leeuw in Netherlands and C. Hayashi in Japan. His type III Quantification methods, published in the 50's (Hayashi, C. (1954). Multidimensional quantification
with applications to analysis of social phenomena. \emph{Annals of the Institute of Statistical Mathematics}, \textbf{5(2):}121--143.), is equivalent to Correspondence Analysis. A historical background 
and synthesis is given in \cite{Tenenhaus1985}. One of the early work in the French school is Cordier, B. - Sur l'analyse Factorielle des Correspondances. \emph{PhD, Rennes}, 1965. This approach has been developed by Greenacre in Pretoria who has studied with J.-P. Benzecri \cite{Greenacre1984}. The algorithm presented here is the one presented in \cite{Nenadic2007}. We have used as well \cite[chapter 10]{Saporta1990} and \cite[sect. 1.3]{LMP2000} for the geometric interpretation.

%
\section{Canonical Correlation Analysis}\label{sec:cca}
%

\subsection*{Notations}

The following notations have been chosen to be as compatible as possible with those in other sections, especially section \ref{sec:pcamet} (\PCA with metrics).\\
\begin{center}
  \ovalbox{
   \begin{tabular}{cll}
      symbol & space & meaning \\
      \hline 
      $A$ & $\R^{n \times p}$ & one of the two matrices to be analysed \\
      $B$ & $\R^{n \times q}$ & one of the two matrices to be analysed \\
      $M$ & $\R^{p \times p}$ & $M=N^{1/2}$ \\
      $N$ & $\R^{p \times p}$ & SDP matrix defining an inner product in $\R^p$; $N = (A^\t A)^{-1}$ \\
      $P$ & $\R^{q \times q}$ & SDP matrix defining an inner product in $\R^q$; $P = (B^\t B)^{-1}$ \\      
      $Q$ & $\R^{q \times q}$ & $Q = P^{1/2}$ \\
      $R$ & $\R^{p \times q}$ & for calculation; $R=MTQ$ \\
      $T$ & $\R^{p \times q}$ & for calculation; $T=A^\t B$ \\
      $v_\a$ & $\R^p$ & a vector in $\R^p$ \\
      $v_\b$ & $\R^q$ & a vector in $\R^q$ \\
      $V_\a$ & $\R^{p \times q}$ & matrix with axis $v_\a$ columnwise \\
      $V_\b$ & $\R^{q \times q}$ & matrix with axis $v_\b$ columnwise \\
      $y_\a$ & $\R^n$ & a vector in $\span A$ \\
      $y_\b$ & $\R^n$ & a vector in $\span B$ \\
      $Y_\a$ & $\R^{n \times q}$ & matrix with components $y_\a$ columnwise \\
      $Y_\b$ & $\R^{n \times q}$ & matrix with components $y_\b$ columnwise \\
   \end{tabular}
  } 
\end{center}

\vspace*{5mm}

\noindent Let us have two data sets as two sets of variables on the same set of items, as $A,B$, with $A \in \R^{n \times p}$ and $B \in \R^{n \times q}$.  We assume that $p,q < n$. The set of columns of each matrix spans a subspace in $\R^n$. If a column of $A$ belongs to the space spanned by the columns of $B$, than there exists a linear regression on the columns of $B$ which explains this column of $A$, and both sets of columns are correlated in $\R^n$. Canonical Correlation Analysis (\CCA) is about finding sets of vectors ( $=$ components) in the spaces spanned by the columns of each matrix with greatest correlation. In this section, the problem is stated (section \ref{sec:cca:set}) and solved (section \ref{sec:cca:solve}) first in an algebraic way, and a second effort must be done to implement involved linear algebra calculations with a reasonable costs (section \ref{sec:cca:numerical}).

\subsection{Stating the problem}\label{sec:cca:set}

We will denote by $v_\a$ a vector in $\R^p$ and by $v_\b$ in $\R^q$. A vector $y_\textsc{a} \in \span A$ (resp. $y_\textsc{b} \in \span B$) can be written as $Av_\a$ (resp. $Bv_\b$). The correlation between $y_\textsc{a}$ and $y_\textsc{b}$ is
\[
 \corr (y_\textsc{a},y_\textsc{b}) = \frac{\langle y_\textsc{a},y_\textsc{b}\rangle}{\|y_\textsc{a}\|\|y_\textsc{b}\|}
\]
Then, Canonical Correlation Analysis of $(A,B)$ for first canonical components can be stated as\\
\\
\begin{center}
 \shadowbox{
  \begin{tabular}{ll}
      Given & $A \in \R^{n \times p}$, $B \in \R^{n \times q}$\\
      \\
      Find & $v_\a \in \R^p$, $v_\b \in \R^q$ \\
      \\
      such that & $\displaystyle \frac{\langle Av_\a,Bv_\b\rangle}{\|Av_\a\|\|Bv_\b\|}$ is maximal \\
  \end{tabular}
 }
\end{center}

\nB As such the problem is difficult to solve. One reason is that $\|v_\a\|$ and $\|v_\b\|$ can take any non zero value (the correlation remains unchanged by a rescaling of $\|v_\a\|$ or $\|v_\b\|$). One could add a constraint like $\|v_\a\|=\|v_\b\|=1$, but the problem still is difficult to solve. There is an equivalent formulation leading to easier calculations for the solution, by setting a  constraint on $y_\a=Av_\a$ (resp. $y_\b = Bv_\b$):
\\
\begin{center}
 \shadowbox{
  \begin{tabular}{ll}
      Given & $A \in \R^{n \times p}$, $B \in \R^{n \times q}$\\
      \\
      Find & $v_\a \in \R^p$, $v_\b \in \R^q$ \\
      with & $\|Av_\a\|^2=1$, $\|Bv_\b\|^2=1$ \\
      \\
      such that & $\langle Av_\a, Bv_\b \rangle$ is maximal \\
  \end{tabular}
 }
\end{center}

\subsection{Solving the problem}\label{sec:cca:solve}

This is an optimization problem with constraints, which can be solved by Lagrange multipliers. Let us recall that if
\[
 \begin{CD}
   \R^n @>f,g>> \R
 \end{CD}
\]
an optimum of $f(x)$ under the constraint $g(x)=0$ is obtained at some points $x$ satisfying
\[
 \nabla\, f - \lambda \nabla \, g=\zeros
\]
where
\[
 \nabla f = \left(\frac{\partial f}{\partial x_1}, \ldots, \frac{\partial f}{\partial x_n}\right)
\]
It remains to check that such a solution is a maximum (it can be a minimum or a saddle point).

\nS Here, the unknowns are $(v_\a, v_\b)$, the function $f$ is $f(v_\a,v_\b)=\langle Av_\a,Bv_\b\rangle$ and $g$ is $\|Av_\a\|^2=\|Bv_\b\|^2=1$. One computes separately the partial derivatives with respect to $v_\a$ and $v_\b$, denoting them $\nabla_{v_\a}$ and $\nabla_{v_\b}$. One has
\begin{equation}
\left\{
 \begin{array}{lclclcl}
 \nabla_{v_\a} \langle Av_\a, Bv_\b \rangle &=& A^\t Bv_\b & \quad& \nabla_{v_\b}\langle Av_\a, Bv_\b \rangle &=& B^\t Av_\a \\ 
 \nabla_{v_\a} \|Av_\a\|^2 &=& 2 A^\t Av_\a &\quad& \nabla_{v_\b} \|Bv_\b\|^2 &=& 2B^\t Bv_\b
 \end{array}
 \right.
\end{equation}
Then, the solution satisfies to
\begin{equation}\label{eq:nabla_can}
 \left\{ 
  \begin{array}{lclcl}
    \nabla_{v_\a}: &\: & A^\t Bv_\b &=& \lambda \, A^\t Av_\a \\
    \nabla_{v_\b}: &\: & B^\t Av_\a &=& \mu \, B^\t Bv_\b
  \end{array}
 \right.
\end{equation}

\nS We first show that $\lambda=\mu$. 
\begin{proof}
Therefore, we observe that
\[
 \left\{ 
  \begin{array}{lcl}
    \langle A^\t Bv_\b,v_\a\rangle  &=& \lambda \langle A^\t Av_\a, v_\a\rangle \\
    \langle B^\t Av_\a, v_\b\rangle &=& \mu \langle B^\t B , v_\b \rangle
  \end{array}
 \right.
\]
and that
\[
 \begin{array}{lclcl}
   \langle A^\t Av_\a, v_\a\rangle &=& \langle Av_\a,Av_\a\rangle &=& 1 \\
   \langle B^\t Bv_\b,v_\b\rangle &=& \langle Bv_\b,Bv_\b\rangle &=& 1
 \end{array}
\]
Then
\[
 \lambda = \langle A^\t Bv_\b,v_\a\rangle = \langle B^\t Av_\a, v_\b\rangle = \mu
\]
\end{proof}

\nS Thus, equation (\ref{eq:nabla_can}) reads
\begin{equation}\label{eq:nabla_can2}
 \left\{ 
  \begin{array}{lcl}
    A^\t Bv_\b &=& \lambda A^\t Av_\a \\
    B^\t Av_\a &=& \lambda B^\t Bv_\b
  \end{array}
 \right.
\end{equation}
Multiplying leftwise the first equation by $(A^\t A)^{-1}$ and the second by $(B^\t B)^{-1}$ yields 
\begin{equation}
 \left\{ 
  \begin{array}{lcl}
    (A^\t A)^{-1}A^\t Bv_\b &=& \lambda v_\a \\
    (B^\t B)^{-1}B^\t Av_\a &=& \lambda v_\b
  \end{array}
 \right.
\end{equation}
and, having in mind that $Av_\a=y_\textsc{a}$ and $Bv_\b=y_\textsc{b}$
\begin{equation}\label{eq:nabla_can3}
 \left\{ 
  \begin{array}{lcl}
    A(A^\t A)^{-1}A^\t y_\textsc{b} &=& \lambda y_\textsc{a} \\
    B(B^\t B)^{-1}B^\t y_\textsc{a} &=& \lambda y_\textsc{b}
  \end{array}
 \right.
\end{equation}

\nS One recognizes in the l.h.s. of (\ref{eq:nabla_can3})
\[
 \begin{cases}
  A(A^\t A)^{-1}A^\t  &= \PC_\textsc{a} \\
  B(B^\t B)^{-1}B^\t  &= \PC_\textsc{b}
 \end{cases}
\]
i.e. the projectors on the spaces spanned by the columns of $A$ and of $B$ respectively. This leads to
\[
 \left\{ 
  \begin{array}{lcl}
    \PC_\textsc{a}y_\textsc{b} &=& \lambda y_\textsc{a} \\
    \PC_\textsc{b}y_\textsc{a} &=& \lambda y_\textsc{b}
  \end{array}
 \right.
\]
or
\begin{equation}\label{eq:nabla_can:res}
 \left\{ 
  \begin{array}{lcl}
    \PC_\textsc{a}\PC_\textsc{b} y_\textsc{a} &=& \lambda^2 y_\textsc{a} \\
    \PC_\textsc{b}\PC_\textsc{a} y_\textsc{b} &=& \lambda^2 y_\textsc{b}
  \end{array}
 \right.
\end{equation}

\nT{Interpretation} The interpretation is quite natural. $\span A$ and $\span B$ are two vector subspaces in $\R^n$ of dimension $p$ and $q$ respectively. Let us have $y_\a \in \span A$. It is projected as $y'_\b \in \span B$ by $y'_b = \PC_\b y_\a$. $y'_b$ itself is projected as $y''_a \in \span A$ by $y''_a= \PC_\a y'_b$. One has 
\[
 \begin{CD}
   \span A @>\PC_\b>> \span B @>\PC_\a>> \span A \\
   y_\a @>>> y'_\b @>>> y''_\a
 \end{CD}
\]
The same can be written for $y_\b$:
\[
 \begin{CD}
   \span B @>\PC_\a>> \span A @>\PC_\b>> \span B \\
   y_\b @>>> y'_\a @>>> y''_\b
 \end{CD}
\]
Equation (\ref{eq:nabla_can:res}) says that when the correlation between $y_\a$ and $y_\b$ is maximal, then $y_\a$ (resp. $y_b$) and $y''_a$ (resp. $y''_b$) are colinear, and eigenvectors of $\PC_\textsc{a}\PC_\textsc{b}$ (resp. $\PC_\textsc{b}\PC_\textsc{a})$.

\notes The computation of the solution in this section is classical, and has been borrowed from \cite[section IV.6.4]{LMF82}.

%
\subsection{Computing the solution}\label{sec:cca:numerical}
%

This solution is geometrically speaking very intuitive, but it does not lead to the most efficient way to compute a solution. We start from
\[
 \left\{ 
  \begin{array}{lcl}
    \PC_\textsc{a}\PC_\textsc{b} y_\textsc{a} &=& \lambda^2 y_\textsc{a} \\
    \PC_\textsc{b}\PC_\textsc{a} y_\textsc{b} &=& \lambda^2 y_\textsc{b}
  \end{array}
 \right.
\]
Let us recall that
\[
 \left\{ 
  \begin{array}{lcl}
    \PC_\textsc{a} &=& A(A^\t A)^{-1}A^\t \\
    \PC_\textsc{b} &=& B(B^\t B)^{-1}B^\t 
  \end{array}
 \right.
\]
Then
\begin{equation}
 \left\{ 
  \begin{array}{lcl}
    A(A^\t A)^{-1}A^\t B(B^\t B)^{-1}B^\t  y_\textsc{a} &=& \lambda^2 y_\textsc{a} \\
    B(B^\t B)^{-1}B^\t A(A^\t A)^{-1}A^\t  y_\textsc{b} &=& \lambda^2 y_\textsc{b}
  \end{array}
 \right.
\end{equation}
We show here how it is possible to avoid the complexity of computing $\PC_\a =  A(A^\t A)^{-1}A^\t$ and $\PC_\b = B(B^\t B)^{-1}B^\t$. Indeed, computing\\
\begin{center}
\begin{tabular}{ccc|ccc}
$A^\t A$ & is in & $\OC(np^2)$ & $B^\t B$ & is in & $\OC(nq^2)$ \\
$(A^\t A)^{-1}$ & & $\OC(p^3)$ & $(B^\t B)^{-1}$ & & $\OC(q^3)$ \\
$A(A^\t A)^{-1}A^\t$ & & $\OC(np^2)$ & $B(B^\t B)^{-1}B^\t$ & & $\OC(nq^2)$
\end{tabular}
\end{center}
Hence, computing $\PC_\a$ (resp. $\PC_\b$) is in $\OC(np^2)$ (resp. $\OC(nq^2)$). We have $\PC_\a, \PC_\b \in \R^{n \times n}$, so $\PC_\a \PC_\b, \PC_\b \PC_\a \sim \OC(n^3)$. However, $\rank \PC_\a = m$ and $\rank \PC_\b=q$. So, it is possible to reduce the complexity of the calculation. 

\nS One possibility is to run a SVD of $A$ and $B$. If the matrices of left singular vectors are denoted respectively $U_\a$ and $U_\b$, the same calculation holds for the projection on the space spanned by the columns of $U_\a$ (resp. $U_\b$) as they are the same as those spanned by $A$ (resp. $B$). We have $U^\t_\a \, U_\a = \I_p$, and $U^\t_\b \, U_\b = \I_q$. Then $\PC(U_\a)= U_\a U^\t_\a$ (resp. $\PC(U_\b)= U_\b U^\t_\b$) and the complexity of the calculation of the projectors is reduced. However, we still have $\PC(U_\a), \PC(U_\b) \in \R^{n \times n}$, and  the calculations of $\PC_\a \PC_\b$ and $\PC_\b \PC_\a$ still are in $\OC(n^3)$. Let us try another way.

\nS Without loss of generality, we assume that $p \geq q$. Let us denote
\begin{equation}
 \begin{cases}
  T &= A^\t B \\
  N &= (A^\t A)^{-1} \\
  P &= (B^\t B)^{-1}
 \end{cases}
\end{equation}
Then
\[
 \left\{ 
  \begin{array}{lcl}
    ANTPB^\t  y_\textsc{a} &=& \lambda^2 y_\textsc{a} \\
    BPT^\t NA^\t  y_\textsc{b} &=& \lambda^2 y_\textsc{b}
  \end{array}
 \right.
\]
Let us recall that
\[
 y_\textsc{a} = Av_\a, \qquad y_\textsc{b} = Bv_\b
\]
Then
\[
 \left\{ 
  \begin{array}{lcl}
    ANTPT^\t  v_\a &=& \lambda^2 Av_\a \\
    BPT^\t NTv_\b &=& \lambda^2 Bv_\b
  \end{array}
 \right.
\]
We can ``simplify'' by $A$ and $B$ by left multiplication by $(A^\t A)^{-1}A^\t $ and $(B^\t B)^{-1}B^\t $. We have
\begin{equation}\label{eq:cca:sol1}
  \left\{ 
  \begin{array}{lcl}
    NTPT^\t  v_\a &=& \lambda^2 v_\a \\
    PT^\t N T v_\b &=& \lambda^2 v_\b
  \end{array}
 \right.
\end{equation}
This reminds of the type of equation of a \PCA with metrics (see equation (\ref{pcamet:eq:res})).

\nS Here, we show how the solution of \CCA can be read as the solution of a \PCA with metrics. Therefore, let us denote, as in section    \ref{sec:pcamet}, 
\[
M=N^{1/2}, \qquad Q=P^{1/2}, \qquad R=MTQ \in \R^{p \times q}
\]
Then, 
\[
\begin{array}{lcl}
   NTPT^\t &=& M^2TQ^2T^\t \\
   &=& M(MTQ)(QT^\t M)M^{-1} \\
   &=& MRR^\t M^{-1} \\
\end{array}
\]
and (the same for $PT^\t N T$) equation (\ref{eq:cca:sol1}) reads
\begin{equation}
 \left\{ 
  \begin{array}{lcl}
    MRR^\t M^{-1}v_\a &=& \lambda^2 v_\a \\
    QR^\t R Q^{-1}v_\b &=& \lambda^2 v_\b
  \end{array}
 \right.
\end{equation}
Let us denote
\[
  w_\textsc{a} = M^{-1}v_\a, \qquad w_\textsc{b} = Q^{-1}v_\b
\]
Then, by left multiplication by $M^{-1}$  of the first equation and by $Q^{-1}$ of the second, we have
\begin{equation}
 \left\{ 
  \begin{array}{lcl}
    RR^\t w_\textsc{a} &=& \lambda^2 w_\textsc{a} \\
    R^\t Rw_\textsc{b} &=& \lambda^2 w_\textsc{b}
  \end{array}
 \right.
\end{equation}
where we recognize the \PCA of $R=MTQ$. The complexity of this calculation is:\\
\begin{center}
    \begin{tabular}{clcl}
    $T$ & $=A^\t B$ & is in & $\OC(npq)$ \\
    $N$ & $=(A^\t A)^{-1}$ & & $O(np^2)$ \\
    $P$ & $=(B^\t B)^{-1}$ & & $O(nq^2)$ \\
    $M$ & $=N^{1/2}$ & & $\OC(p^3)$ \\
    $Q$ & $=P^{1/2}$ & & $\OC(q^3)$ \\
    $R$ & $=MTQ$ & & $\OC(p^2q)$ \\
    \end{tabular}
\end{center}

\nS The \PCA of $R$ is the \PCA of $T$ with metrics defined by $N$ on $\R^p$ and $P$ on $\R^q$. Let us note as well that $w_\a$ is a principal axis of the \PCA of $R^\t$, and $w_\b$ is a principal axis of the \PCA of $R$. As $R \in \R^{p \times q}$ and as we have assumed that $p \geq q$, it is natural to run the \PCA of $R$, hence compute the $w_\b$ as a principal axis of $R$. Then, $w_\a$ as a principal axis of $R^\t$ is a principal component of $R$ and related to $w_\b$ by
\begin{equation}
 w_\a=Rw_\b, \qquad \mbox{with} \quad 
 \begin{cases}
  R & \in \R^{p \times q} \\
  w_\a & \in \R^p \\
  w_\b & \in \R^q
 \end{cases}
\end{equation}
Hence the SVD or search for eigenvalues and eigenvectors will be done once only.

\nT{Interpretation}  Hence the result: the solution of the CCA of two arrays $A$ and $B$ is the solution of the \PCA of $T=A^\t B$ with an inner product defined by $N$ on rows and by $P$ on columns where $N = (A^\t A)^{-1}$ and $P= (B^\t B)^{-1}$. We then have 
\begin{equation*}
 v_\a = Mw_\a, \qquad v_\b = Qw_\b
\end{equation*}
and 
\begin{equation*}
 y_\a = Av_\a, \qquad y_\b = Bv_\b
\end{equation*}

\nT{Algorithm} There are several ways to write an algorithm for this calculation. Here is a direct one, without calling $\textsc{pca}$ or $\textsc{pca-met}$.\\
\begin{algorithm}[H]
\begin{algorithmic}[1]
\STATE \textbf{input} $A \in \R^{n \times p}$, $B \in \R^{n \times q}$ with $p \geq q$
\STATE \textbf{compute} $T = A^\t B$
\STATE \textbf{compute} $N = (A^\t A)^{-1}$
\STATE \textbf{compute} $P = (B^\t B)^{-1}$
\STATE \textbf{compute} $M=N^{1/2}$
\STATE \textbf{compute} $Q=P^{1/2}$
\STATE \textbf{compute} $R = MTQ$
\STATE \textbf{compute} $W_\a, \Lambda, W_\b = \textsc{svd}(R)$
\STATE \textbf{compute} $V_\a = MW_\a$
\STATE \textbf{compute} $V_\b = QW_\b$
\STATE \textbf{compute} $Y_\a = AV_\a$
\STATE \textbf{compute} $Y_\b = BV_\b$
\RETURN $Y_\a, Y_\b , V_\a, V_\b, \Lambda$
\end{algorithmic}
\caption{Canonical Correlation Analysis: \textsc{cca}$(A,B)$}
\label{alg:can}
\end{algorithm}

\paragraph{Comments:} Here are some comments on the algorithm:
\begin{itemize}[label=$\rightarrow$]
    \item The components $y_\a$ are the columns of $Y_\a$
    \item The components $y_\b$ are the columns of $Y_\b$
    \item the axis $v_\a$ are the columns of $V_\a$
    \item the axis $v_\b$ are the columns of $V_\b$
    \item the singular values of $R$ are the correlation coefficients $\lambda$, because $R^\t R w_\a = \lambda^2 w_\a$, and the eigenvalues of $R^\t R$ are the square of the singular values of $R$
    \item it is possible to compute $M=N^{1/2}$ through a SVD of $N$: if $N=U_\m \Sigma_\m V_\m$, then $M = U_\m \Sigma_\m^{1/2} V_\m$ and $\Sigma_\m$ is diagonal. The same for computing $Q$ from $P$.
\end{itemize}

\notes Canonical Analysis seems to have been proposed by Hotelling in 1936 in Hotelling, H. - Relation between two sets of variables, \emph{Biometrika}, \textbf{28:}361-377. It is presented with a statistical approach in  \cite[chap. 12]{Anderson1958} or \cite[Section 8f]{Rao1973}.  A review of Canonical Analysis with (debated) applications in ecology can be found in \cite{Gittins1985}. It is presented with a more algebraic approach in classical textbooks of the French school of data analysis, e.g.  \cite{Lebart1977,LMF82,Escofier1990,Saporta1990} from which it has been borrowed and adapted here. Canonical Analysis is often referred to as Canonical Correlation Analysis (\CCA). Both denominations will be used here.

%
\section{Multiple Correspondence Analysis}\label{sec:mcoa}
%

Here, we develop a way to extend \CCA with more than 2 arrays $A$ and $B$. There are different ways to do it, which are not equivalent. Indeed, let us have 3 arrays, $A$, $B$ and $C$. Let us select one component per array, respectively $y_\a=Au_\a$, $y_\b=Bu_\b$ and $y_\c=Cu_c$ with $\|y_\a\|=\|y_\b\|=\|y_\c\|=1$. One can define three-ways \CCA as finding a triplet $(y_\a, y_\b, y_\c)$ such that their correlation is maximal. There is however no canonical way to define the correlation between 3 vectors. It can be defined from $\|y_\a \wedge y_\b \wedge y_\c\|$ (vectors are independent when their wedge product is maximal), or $\|y_\a + y_\b + y_\c\|$, or other ways. Here, we extend \CCA to more than two arrays along a way which passes through an equivalence between \CCA and \CoA. We then extend \CoA to more than two variables, and come back to \CCA through the equivalence between \CoA and \CCA.

\subsection{A tight link between Canonical Analysis and Correspondence Analysis}

Let us consider a set of two qualitative variables $A$ and $B$ on a same set of items $\llbracket 1, n \rrbracket$. We build the so called indicator array of each variable, called as well completely disjunctive table. If there are $n$ items, $p$ modalities for $A$ and $q$ for $B$, it is a $n \times p$ array for $A$ and $n \times q$ for $B$, with, in each row $i$, zero in each column but 1 in column $j$ if the modality $j$ of the variable has been observed for item $i$. 
\begin{equation*}
 \alpha_{ij} = 
 \begin{cases}
    1 & \mbox{if modality}\: j \: \mbox{has been observed for item}\: i \\
    0 & \mbox{otherwise}
 \end{cases}
\end{equation*}
\begin{equation*}
 \beta_{ik} = 
 \begin{cases}
    1 & \mbox{if modality}\: k \: \mbox{has been observed for item}\: i \\
    0 & \mbox{otherwise}
 \end{cases}
\end{equation*}
Let us do the Canonical Analysis of both arrays $A$ and $B$. The solution is given by
\[
w = M^\t Mw, \qquad  v= D_\textsc{b}^{1/2}w
\]
with
\[
 M = D_\textsc{a}^{1/2} T D_\textsc{b}^{1/2}, \qquad T=A^\t B, \qquad D_\textsc{a} = (A^\t A)^{-1}, \qquad D_\textsc{b} = (B^\t B)^{-1}
\]

\paragraph{Key observations:} Let us note three things:
\begin{itemize}[label=$\rightarrow$]
    \item $(A^\t A)^{-1} \in \R^{p\times p}$ (resp. $(B^\t B)^{-1} \in \R^{q \times q}$) is the diagonal matrix with in position $j$ (resp. $k$) the inverse of the number of rows in $A$ (resp. $B$) where the modality $j$ (resp. $k$) of the variable has been observed.
    \item One recognizes in $T$ the contingency table between $A$ and $B$, %
    \item and in the PCA of $M$ the \PCA of $T$ with inner product defined by $D_{\textsc{a}}^{1/2}$ in $\R^p$ and by $D_{\textsc{b}}^{1/2}$ in $\R^q$, i.e. correspondence analysis of $T$.
\end{itemize}
This leads to an intimate link between Correspondence Analysis and Canonical Analysis, which permits further an extension of Correspondence Analysis to more than two variables.

%
\subsection{Link between Canonical Analysis and PCA with metric on rows}
%

Let $A \in \R^{n \times p}$ and $B \in \R^{n \times q}$ with $q \leq p$ and $p+q \leq n$ be two data sets, each a set of quantitative variables (the columns of $A$ and $B$) on a same set of items (the rows of $A$ and $B$). Let us have in mind the Canonical Analysis of $(A,B)$, which will be developed here along another calculation. 

\nB Therefore, let us consider the data set
\[
 X = [A|B] \in \R^{n \times (p+q)}
\]
built by columnwise concatenation of $A$ and $B$. Let us define
\[
\left|
 \begin{array}{lclcl}
  D_\textsc{a} &=& (A^\t A)^{-1} & \quad & \in \R^{p \times p} \\
  D_\textsc{b} &=& (B^\t B)^{-1} & \quad & \in \R^{q \times q} \\  
 \end{array}
 \right.
\]
and
\[
 D = 
 \begin{pmatrix}
   D_\textsc{a} & 0 \\
   0 & D_\textsc{b}
 \end{pmatrix}
 \in \R^{(p+q) \times (p+q)}
\]
Let us perform the \PCA of $X$ with metrics on rows given by $D$ (a row of $X$  belongs to $\R^{p+q}$). Let us denote
\[
 T = A^\t B, \qquad M = D_{\textsc{a}}^{1/2}TD_{\textsc{b}}^{1/2}
\]
Then, performing the \PCA of $X$ with row distances given by $D$ is performing the \PCA of $M$. If $A$ and $B$ are indicator arrays, this is performing the \CoA of $T$, which is the contingency table of $A$ and $B$. 

\begin{proof}
The solution is performing a \PCA of
\[
 R = XD^{1/2}, \qquad R \in \R^{n \times (p+q)}
\]
i.e. performing an EVD of $R^\t R$. We will denote the matrix dimensions under each block. We have
\[
 \underset{n \times (p+q)}{R} = \left[ \underset{n \times p}{AD_{\textsc{a}}^{1/2}} \; , \; \underset{n \times q}{BD_{\textsc{b}}^{1/2}} \right],
 \qquad
 \underset{(p+q) \times n}{R^\t } = \left[ 
 \begin{array}{c}
 \underset{p \times n}{D_{\textsc{a}}^{1/2}A^\t } \\
 \\
 \underset{q \times n}{D_{\textsc{b}}^{1/2}B^\t } 
 \end{array}
 \right]
\]
So
\[
 R^\t R = \left[
 \begin{array}{ccc}
   \underset{p \times p}{D_{\textsc{a}}^{1/2}A^\t AD_{\textsc{a}}^{1/2}} & & \underset{p \times q}{D_{\textsc{a}}^{1/2}A^\t BD_{\textsc{b}}^{1/2}} \\
   &&\\
   \underset{q \times p}{D_{\textsc{b}}^{1/2}B^\t AD_{\textsc{a}}^{1/2}} & & \underset{q \times q}{D_{\textsc{b}}^{1/2}B^\t BD_{\textsc{b}}^{1/2}} \\
 \end{array}
 \right]
\]
Let us observe that
\[
 A^\t A = D_{\textsc{a}}^{-1}, \qquad B^\t B = D_{\textsc{b}}^{-1}
\]
Hence
\[
 R^\t R = \left[
 \begin{array}{ccc}
  \I_p & & \underset{p \times q}{D_{\textsc{a}}^{1/2}A^\t BD_{\textsc{b}}^{1/2}} \\
   &&\\
   \underset{q \times p}{D_{\textsc{b}}^{1/2}B^\t AD_{\textsc{a}}^{1/2}} & & \I_q \\   
 \end{array}
 \right]
\]
Let
\[
 x = \left[\begin{array}{c}u\\v\end{array}\right]
\]
be such that
\[
 R^\t Rx = \lambda x
\]
This yields
\[
 \left\{ 
  \begin{array}{lcl}
   D_{\textsc{a}}^{1/2}A^\t BD_{\textsc{b}}^{1/2} v &=& (\lambda-1)u \\
   D_{\textsc{b}}^{1/2}B^\t AD_{\textsc{a}}^{1/2} u &=& (\lambda-1)v
  \end{array}
 \right.
\]
Let us denote
\[
 T = A^\t B, \qquad M = D_{\textsc{a}}^{1/2}TD_{\textsc{b}}^{1/2}
\]
Then
\[
 \left\{ 
  \begin{array}{lcl}
   Mv &=& (\lambda-1)u \\
   M^\t u &=& (\lambda-1)v
  \end{array}
 \right.
\]
or
\[
 \left\{ 
  \begin{array}{lcl}
   M^\t Mv &=& (\lambda-1)^2v \\
   MM^\t u &=& (\lambda-1)^2u
  \end{array}
 \right.
\]
where we recognize the \PCA of $M$, hence the \PCA of $T=A^\t B$ with metrics defined by $(A^\t A)^{-1}$ on columns and $(B^\t B)^{-1}$ on rows, hence the \verb!CoA! of $T$ as a contingency table.
\end{proof}

\subsection{Multiple Canonical Analysis}

Knowing that, the extension if Canonical Analysis to more than two quantitative variables is straightforward.

\nS Let us have $m$ qualitative variables on $n$ items, each with indicator matrix $A_\ell$ $(\ell \in \llbracket 1,m \rrbracket)$ with
\[
 A_\ell \in \R^{n \times p_\ell}, \qquad \sum_\ell p_\ell \leq n
\]
Let us define
\[
 D_\ell = (A^\t_\ell A_\ell)^{-1}, \: \in \R^{p_\ell \times p_\ell}
\]
Let us build
\[
 X = [A_1 | \ldots | A_m]
\]
and 
\[
 D = \left[
 \begin{array}{ccccc}
  D_1 & 0   & \ldots & \ldots &  0\\
  0   & D_2 & 0      & \ldots & \vdots \\
  \vdots & \ddots & \ddots & \ddots & \vdots \\
   \vdots & \ddots & \ddots & \ddots & 0 \\
   0 & \ldots & \ldots & 0 & D_m\\
 \end{array}
 \right]
\]
Then, the \MCoA of $(A_1, \ldots,A_m)$ is the \PCA of $X$ with distances on rows given by $D$. It is the \PCA of
\begin{equation}
 R = [A_1D_1 | \ldots | A_mD_m]
\end{equation}
We have
\[
 R^\t = \left[ 
 \begin{array}{c}
   D_1A^\t_1 \\
   \vdots \\
   D_mA^\t_m
 \end{array}
 \right]
\]
and $M=R^\t R$ can be written blockwise as
\[
 \begin{cases}
   M_{\ell\ell} = \I_{p_\ell} & \mbox{if}\:\:\ell = \ell'\\
   M_{\ell\ell'} = D_\ell A^\t_\ell A_{\ell'}D_{\ell'} & \mbox{if}\:\:\ell \neq \ell'
 \end{cases}
\]

%
\subsection{Summary of relationships between some methods}
%

Let us recall that, in this document, we denote:\\
\begin{center}
    \begin{tabular}{lcl}
         \verb!CoA! &:& Correspondence Analysis \\
         \verb!Can! &:& Canonical Analysis \\
         \verb!PCAmet! &:& Principal Component Analysis with metrics 
    \end{tabular}
\end{center}
Let us recall that:
\begin{itemize}[label=$\rightarrow$]
    \item \verb!CoA! is the PCA of a contingency table $T$ with metrics on rows and columns as the inverse of the row and column marginals,
    \item \verb!Can! is the analysis of two quantitative arrays $(A,B)$ in order to find the most common components,
    \item \verb!PCAmet!  the PCA of a quantitative array $A$ with metrics defined by $N$ on columns and $P$ on rows.
\end{itemize}

\nD Let $(A,B)$ be two indicator arrays of one categorical variable each on the same items. Then, The Canonical Analysis of $(A,B)$ is equivalent to the Correspondence Analysis of $X=A^\t B$
\begin{equation}
    \texttt{Can}(A,B) = \texttt{CoA}(A^\t B)
\end{equation}

\nD Let $(A,B)$ be two quantitative arrays. Let us consider their concatenation $X = [A|B]$, and the \verb!PCAmet! of $X$ with metrics defined by $D$ on rows, with
\[
 D = 
 \begin{pmatrix}
   D_\textsc{a} & 0 \\
   0 & D_\textsc{b}
 \end{pmatrix}
 \qquad \mbox{with} \quad 
 \left\{
 \begin{array}{lclcl}
  D_\textsc{a} &=& (A^\t A)^{-1} & \quad & \in \R^{p \times p} \\
  D_\textsc{b} &=& (B^\t B)^{-1} & \quad & \in \R^{q \times q} \\  
 \end{array}
 \right.
\]
Then, the \verb!PCA! of $X$ with metrics on rows defined by $D$ is equivalent to performing the Canonical Analysis of $(A,B)$
\begin{equation}
    \texttt{Can}(A,B) = \texttt{PCAmet}(X,D)
\end{equation}

\nD Let us now assume that $A$ and $B$ are each the indicator array of a categorical variable. Then, $\texttt{Can}(A,B)$ is equivalent to the Correspondence Analysis of $T=A^\t B$. Then, $\texttt{PCAmet}(X,D)$, equivalent to $\verb!Can!(A,B)$, is equivalent to the Correspondence Analysis of $T=A^\t B$ as well by transitivity. 

\nD The Canonical Analysis of two tables $(A,B)$ as the \verb!PCAmet! of $X=[A|B]$ with metric defined by $D$ on the rows of $X$ has been extended to the analysis of $m$ arrays $A_\ell$ as the \verb!PCAmet! of table
\[
X = [A_1| \ldots|A_m]
\]
with metric defined by matrix $D$ blockwise diagonal defined as
\[
D = \mathrm{Diag}(D_\ell) \qquad \mbox{with} \quad D_\ell = (A_\ell^\t A_\ell)^{-1}
\]
This leads to Multiple Correspondence Analysis (\verb!MCA!) as follows.

%
\subsection{Multiple Correspondence Analysis}
%

Let us have $m$ categorical variables observed each on $n$ items. Let us denote by $A_\ell$ with $\ell \in \llbracket 1,m\rrbracket$ the indicator array of variable $\ell$, i.e. is a $n \times p_\ell$ binary array with
\begin{equation*}
 \alpha_{ij_\ell} = 
 \begin{cases}
    1 & \mbox{if modality}\: j_\ell \: \mbox{has been observed for item}\: i \\
    0 & \mbox{otherwise}
 \end{cases}
\end{equation*}
Let us define
\[
X = [A_1| \ldots|A_m]
\]
and the metric on rows of $X$ defined by the matrix $D$ blockwise diagonal defined as
\[
D = \mathrm{Diag}(D_\ell) \qquad \mbox{with} \quad D_\ell = (A_\ell^\t A_\ell)^{-1}
\]
Then, the Multiple Correspondence Analysis of $(A_1, \ldots,A_m)$ is equivalent to the \verb!PCAmet! of $X$ with metric defined by $D$ on rows.

\notes These links between different treatments of a same dataset which establish some dependencies between methods have been subject to thorough studies and presentations by the French school of multivariate analysis, more algebraic and geometrical than statistical, in the 70's. with the names of B. Escofier, J.-P. Fénelon, L. Lebart, A. Morineau, J. Pagès, among others. It is presented in all textbooks of this school in multivariate data analysis, under chapters called ``other methods and complements'', like \cite{Lebart1977,LMF82,LMP2000}. Since then, the diversity of related methods have flourished, and a more recent panorama is given in \cite{Greenacre2006}. According to the introduction of \cite{Greenacre2006}, L. Guttman should be credited for all the basic ideas at the root of Multiple Correspondence Analysis, in Guttman, L. 1941. \emph{The quantification of a class of attributes: A theory and method of scale construction}, Chapter in \emph{The Prediction of Personal Adjustment}, P. Horst, eds. New York : Social Science Research Council. The article \cite{Tenenhaus1985} is a thorough survey of different methods associated with Multiple Correspondence Analysis, with both a historical background (and they point out several independent beginnings in the 30's and early 40's) and, long before the present notes, an organization of methods to show that they all lead to the same equation  to analyze the data. The unifying formalism selected in \cite{Tenenhaus1985} is the duality diagram.

%
\section{Multidimensional Scaling}\label{sec:mds}
%

Multidimensional Scaling is a technique to map a discrete metric space into a Euclidean space. Let $(M,d)$ be a metric space with $|M|=n$. It is given by a pairwise distance matrix
\[
 D = [d_{ij}]_{i,j} \in \R^{n \times n} \qquad \mbox{with} \quad 1 \leq i,j \leq n,
\]
such that 
\[
 d_{ij} = d(i,j).
\]
MDS at dimension $r$ addresses the question of finding a point cloud 
\[
 \X = (x_i)_{1 \leq i \leq n}, \qquad x_i \in \R^r
\]
such that the distance between points $x_i$ and $x_j$ is as close as possible from $d_{ij}$ or, loosely speaking
\[
 \|x_i-x_j\| \approx d_{ij}
\]
A matrix $X \in \R^{n\times r}$ is attached to the point cloud $\X$, with $x_i$ being row $i$ of $X$.

\nS Then, two situations may occur
\begin{enumerate}
 \item\label{mds:class} either the distances $d_{ij}$ come from a Euclidean distance between (unknown) points, and the problem is to recover them, i.e. produce an isometry, and a best approximation of it at dimension $r$
 \item\label{mds:ls} or they do not come from a Euclidean distance, and a best approximation is sought for, knowing there exists no exact isometry
\end{enumerate}

\noindent Problem \ref{mds:class} is known as \emph{classical MDS} and problem \ref{mds:ls} as \emph{Least Square Scaling}. In this note, we address the problem of classical MDS only. Least square scaling is a delicate and non trivial optimization problem. It appears however that, in practical cases, one never knows whether the distances come from a Euclidean point cloud or not, and the procedure is ``as if''.

\nS One thing to understand is that classical \MDS assumes that the distances in metric space $(M,d)$ are $\ell^2-$norm distance in an (unknown) Euclidean space. If it is not the case, a numerical problem occurs (but a \emph{ad hoc} standard solution is provided). Three points with pairwise distances can be arranged as a triangle in $\R^2$, and four points with pairwise distances can be arranged as a tetrahedron in $\R^3$. More generally, $n$ points with pairwise Euclidean distances can be isometrically embedded in $\R^{n-1}$ at most. Then, \MDS at rank $r$ is made in two steps:
\begin{enumerate}
 \item find a point cloud $X$ in $\R^{n-1}$ with distances matching exactly the values of $d(i,j)$ for each pair
 \item reduce the dimension $r < n$ of the space where to build the point cloud $X_r$ as close as possible to $X$. This can be solved by \PCA of $X$. It will appear that principal axis and components of the \PCA of $X$ can be given by \MDS without further calculations.
\end{enumerate}

The procedure to build matrix $X$ attached to  point cloud $\X$ is in four steps, presented hereafter\\
\begin{algorithm}[H]
\begin{algorithmic}[1]
\STATE \textbf{given} a distance matrix $D \in \R^{n \times n}$ and a dimension $r< n$
\STATE \textbf{compute} the Gram matrix $G$ associated to $D$ 
\STATE \textbf{compute} the EVD of $G$ with eigenvalues $\Lambda$
\STATE \textbf{compute} the coordinates $X \in \R^{n \times n}$ of point cloud $\X$ from the EVD of $G$
\STATE \textbf{compute} the best low rank approximation $X_r \in \R^{n \times r}$ of $X$
\RETURN $X_r, \Sigma_r$
\end{algorithmic}
\caption{General procedure for classical MDS}
\label{alg:mds}
\end{algorithm}
\noindent and developed further in this section. 

\notes There exists many excellent textbooks presenting MDS, classical MDS or LSS. We can recommend \cite{Cox2001} or \cite[chap. 13]{Izenman2008}. A comprehensive reference is \cite{Borg2005}. A classical and rigorous reference with many results, their demonstration and history is \cite{Mardia1979} which we highly recommend for those enjoying a mathematical based approach. Classical MDS has been proposed by Torgerson in 1952 \citep{Torgerson1952}. Here, we have followed \cite[chap. 2]{Cox2001}.

%
\subsection{The Gram matrix}
%

Let $X=(x_i)_i$ be a cloud on $n$ points in $\R^r$, such that the pairwise distances only are known.
\[
 \|x_i-x_j\| = d_{ij}
\]
and not the coordinates of the $x_i$. The Gram matrix $G \in \R^{n \times n}$ of $X$ is the matrix with elements
\[
 G_{ij} = \langle x_i,x_j\rangle
\]

\nS There is a well known correspondence between the Gram Matrices $G$ of inner products $\langle x_i,x_j\rangle$ and the Euclidean Distance 
Matrix of quantities $\|x_i-x_j\|$, both $n \times n$. If
\begin{equation*}
 \left|
    \begin{array}{lcl}
      g_{ij} &=& \langle x_i,x_j\rangle \\
      d_{ij}^2 &=& \|x_i-x_j\|^2
    \end{array}
 \right. 
\end{equation*}
Then 
\begin{equation}\label{eq:gramdis}
  \left\{
  \begin{array}{lcl}
   d_{ij}^2 &=& g_{ii} + g_{jj} -2g_{ij} \\
g_{ij} &=&  \displaystyle -\frac{1}{2}\left(d_{ij}^2-d_{i\bullet}^2-d_{j\bullet}^2 + d_{\bullet\bullet}^2 \right)
  \end{array}
  \right.
\end{equation}
with
\begin{equation}
  \left\{
    \begin{array}{lcl}
     d_{i\bullet}^2 &=& \displaystyle \frac{1}{n}\sum_jd_{ij}^2 \\
     \\
     d_{\bullet\bullet}^2 &=& \displaystyle  \frac{1}{n^2}\sum_{i,j}d_{ij}^2 = \frac{1}{n}\sum_id_{i\bullet}^2
    \end{array}
  \right.  
\end{equation}
Such a correspondence has been studied for decades (see e.g. \cite{Schoenberg1938,Laurent1998}). We then have the scheme (with an arrow meaning ``built from''):\\
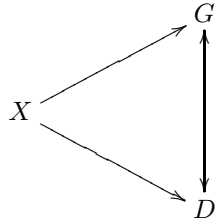
\[
\xymatrix{
& & G \ar @{<->} [dd]\\
X \ar [urr] \ar [drr] & & \\
& & D 
}
\]

\nS A key question is to know whether one-way arrows $X \longrightarrow G$ and $X \longrightarrow D$ can be reversed, i.e. whether $X$ can be computed knowing $G$ or knowing $D$. This amounts to answering to the question: a pairwise distance matrix $D$ being given, is there a dimension $m$ and a point cloud $X$ in $\R^m$ such that the distance between $x_i$ and $x_j$ is precisely $d_{ij}$?

\nS A matrix $D$ being given, it is always possible to compute a matrix $G$ by equation (\ref{eq:gramdis}). But $G$ is not necessarily positive, i.e. it is not necessarily the Gram matrix of a point cloud $X$. The conditions on $D$ for $G$ to be positive, i.e. a Gram matrix have been thoroughly studied. In most of the case, when the Gram matrix is not positive, the negative eigenvalues are just ignored. This leads to some subtleties when connecting the EVD and the SVD of the Gram matrix to compute the coordinates.

\nS The coordinates of the point cloud $X$ can be computed from the eigenvectors and eigenvalues, or the Singular Value Decomposition of the Gram matrix. The recipe is given here.

%
\subsection{Eigendecomposition of the Gram Matrix}\label{sec:mds:svd}
%

Let $G$ be the Gram matrix. If (this is an hypothesis) there exists a set of $n$ points in $\R^m$ such that
\[
 \forall \: i,j, \quad \|x_i-x_j\| = d_{ij}
\]
then
\begin{equation}
 G_{ij} = \langle x_i,x_j\rangle
\end{equation}
and $G$ is positive. We assume here that $G$ as computed from equation (\ref{eq:gramdis}) is positive.

\nS The objective is to associate to $(M,d)$ a point cloud $\X$ in a Euclidean space such that, if possible, $d(i,j) = \|x_i-x_j\|$. Let $X \in \R^{n \times m}$ be the matrix with row $i$ being $x_i$. Then
\begin{equation}
 G = XX^\t
\end{equation}
Next step is about computing $X$ knowing $XX^\t$.

\nS Let $(u_\alpha,\lambda_\alpha)_\alpha$ be the set of eigenpairs of $G$
\begin{equation}
 Gu_\alpha = \lambda_\alpha u_\alpha
\end{equation}
with
\[
 \lambda_1 \geq \lambda_2 \geq \ldots \geq \lambda_n \geq 0
\]
As $G$ is symmetric, the eigenvectors if normed form an orthonormal family. If $U=[u_1|\ldots|u_n]$ is the matrix with $u_\alpha$ in column $\alpha$ and $\Lambda$ the diagonal matrix with $\lambda_\alpha$ in its diagonal, we have
\begin{equation}
 GU=U\Lambda
\end{equation}
As $U$ is orthonormal, $UU^\t=\I_n$, and we have by right multiplication by $U^\t$
\begin{equation}
G = U\Lambda U^\t
\end{equation}
We recognize here the SVD of $G$ if $G$ is positive. The case where $G$ is not positive is handled by designing a quadratic embedding, and is addressed in section \ref{sec:quadratic_embedding}. Let us note that the standard practice when $G$ is not positive is not to design a quadratic embedding, but to clip to zero the non positive eigenvalues and eigenvectors, i.e. to keep track of positive eigenvalues only and associated eigenvectors. If $G$ is positive, the eigenvalues of $G$ are its singular values, and if $G$ is not positive, the negative eigenvalues are singular values up to their sign (if $\lambda < 0$ is an eigenvalue of $G$, $-\lambda >0$ is a singular value of $G$). 

\nS As $G$ is definite positive, let
\[
 \Sigma = \Lambda^{1/2}
\]
Then
\begin{equation}
 G = (U\Sigma)(U\Sigma)^\t
\end{equation}
and we can select
\begin{equation}
 X = U\Sigma, \qquad \Sigma = \Lambda^{1/2}
\end{equation}
It is not the only solution, because any matrix $X'=X\Omega$ where $\Omega$ is a rotation in $\R^m$ is a solution too.

\nS Hence the algorithm:\\
\begin{algorithm}[H]
\begin{algorithmic}[1]
\STATE \textbf{input:} Gram matrix $G$
\STATE \textbf{compute} eigendecomposition of $G$: $GU=U\Lambda$, with eigenvectors $u_\alpha$ (columns of $U$) and eigenvalues $\lambda_\alpha$: $Gu_\alpha = \lambda_\alpha u_\alpha$
\STATE \textbf{compute} $\Sigma = \Lambda^{1/2}$
\STATE \textbf{compute} $X=U\Sigma$
\RETURN $X,\Lambda$
\end{algorithmic}
\caption{MDS with eigendecomposition of Gram matrix}
\label{alg:mds:evd}
\end{algorithm}
\noindent One notices that for $X$ to be computed this way, one must have $\Lambda \geq 0$ for $\Lambda^{1/2}$ to be real. It is one condition for an isometry between $(M,d)$ and a Euclidean space to exist.

\subsection{Dimension reduction}

Once matrix $X$ has been computed, finding a point cloud $X_r$ of $n$ points in $\R^r$ as close as possible to $\X$ is done by the \PCA of $X$. 

\nS Let us recall that
\begin{equation}
 X = U\Sigma^{1/2} \qquad \mbox{with} \quad U^\t U = \I_n
\end{equation}
We recognize here the SVD of $X$ as $X=U\Sigma^{1/2} V^\t$, with singular values being the diagonal of $\Sigma^{1/2}$ and $V=\I_n$. Hence, $X=U\Sigma^{1/2}$ is the matrix of principal components of $X$, and dimensionality reduction of $X$ by PCA is given as a cherry on top.

\subsection{\MDS algorithm}

Wrapping all this together yields the following algorithm for  \MDS\\
\begin{algorithm}[H]
\begin{algorithmic}[1]
\STATE \textbf{input:} a distance matrix $D \in \R^{n \times n}$; a dimension $r < n$
\STATE \textbf{compute} the Gram matrix of $D$: $G=\textsc{gram}(D)$
\STATE \textbf{compute} the eigenpairs $(u_\alpha,\lambda_\alpha)$ such that $Gu_\alpha = \lambda_\alpha u_\alpha$
\STATE keep in $U$ the columns associated to non-negative eigenvalues of $G$; clip them off in $\Lambda$
\STATE \textbf{compute} $\Sigma = \Lambda^{1/2}$
\STATE \textbf{compute} $X=U\Sigma$
\STATE keep in $X_r$ the $r$  first columns of $X$ only, and in $\Lambda_r$ the first $r$ non negative eigenvalues of $G$
\RETURN $X_r,\Lambda_r$
\end{algorithmic}
\caption{Classical MDS: $X, \Sigma = \textsc{D},r$}
\label{alg:mds:global}
\end{algorithm}

\nT{Note} If the user selects the computation of the eigenvalues of $G$, it is simple to detect those which are negative, and clip the corresponding columns in $U$ (and values in $\Lambda$). If the SVD is selected, we have $G=U\Sigma V^\t$ (classical notation $(U,\Sigma,V$, here $\Sigma \neq \Lambda^{1/2}$), with columns of $V$ (resp. singular values) being the same as the columns of $U$ (resp. eigenvalues of $G$) up to the sign, depending on the sign of the corresponding eigenvalue of $G$: 
\begin{center}
\begin{tabular}{|llll}
$\lambda_\alpha >0$ & $\Longrightarrow$ & $v_\alpha=u_\alpha$ & $\sigma_\alpha=\lambda_\alpha$ \\ 
$\lambda_\alpha <0$ & $\Longrightarrow$ & $v_\alpha=-u_\alpha$ & $\sigma_\alpha=-\lambda_\alpha$ \\ 
\end{tabular}
\end{center}

%
\subsection{Quadratic embedding}\label{sec:quadratic_embedding}
%

There is still one point to look at. In all these developments, we have assumed that the Gram matrix $G$ built from distance matrix $D$ with recipe in equation (\ref{eq:gramdis}) is positive, i.e. that all eigenvalues of $G$ (the $\lambda_\alpha$) are non negative. This is a property of a Gram matrix, but is not always the case for the matrix $G$ computed with real data. If one eigenvalue at least is negative, the matrix $G$ is more strictly called a \emph{kernel matrix}, and there is no isometry between $(M,d)$ and a Euclidean space, whatever its dimension. However, it can be shown that there exists a quadratic embedding between $(M,d)$ and a pseudo-euclidean space, i.e. a vector space with a quadratic form with a signature $(p,m)$ (see appendix \ref{sec:quadratic} for a short introduction to quadratic forms and spaces). Such an embedding is rarely done, and most of the time the axis (and components) associated to negative eigenvalues of $G$ simply are ignored, or clipped to zero. We show in this section how to take into account the axis associated to negative eigenvalues of the kernel matrix.

\paragraph{A simple example of a discrete metric space without isometry in a Euclidean space:} A discrete metric space $(M,d)$ being given, there exists not systematically an isometry 
(metric embedding) into an Euclidean space preserving the distances. A simple example is
\[
M = \{i,j,k,c\}
\]
which is the set of vertices of a graph with $E = \{(i,c),(j,c),(k,c)\}$ (there is no conflict of notations between $E$ which is here classically the set of edges and $E$ which is in this section a vector space) with the distance between two vertices being the length of the shortest path between them. We then have
\[
D = \quad 
\begin{tabular}{c|cccc}
& $i$ & $j$ & $k$ & $c$ \\
\hline
$i$ & 0 & 2 & 2 & 1 \\
$j$ & 2 & 0 & 2 & 1 \\
$k$ & 2 & 2 & 0 & 1 \\
$c$ & 1 & 1 & 1 & 0
\end{tabular}
\]
There exists however a map into a pseudo-euclidean space $(E,q)$ as host space, 
preserving a quadratic form.

\paragraph{Quadratic embedding:}  Let $(M,d)$ be a discrete metric space, with $M = \{1, \ldots,i, \ldots, n\}$. Let $(E,q)$ be a quadratic space, with $\dim E=n$, and $\varphi$ a map 
\begin{equation*}
    \begin{CD}
       M @>\varphi>> E
    \end{CD}
\end{equation*}
Let us denote 
\[
a_i = \varphi(i)
\]
So, $a_i \in E$, and we have a cloud of $n$ points in $E$, each one corresponding to an element in $M$. Then, $\varphi$ is a quadratic embedding, or preserves the distances, if 
\begin{equation}
    q(a_i-a_j) = d^2(i,j) 
\end{equation}
We show next
\begin{itemize}[label=$\diamond$]
\item that such an embedding exists (and is not unique),
\item how to build it by specifying $q$ and $\varphi$.
\end{itemize}
We first build $q$, and next $\varphi$.

\paragraph{Building a quadratic form:} The key is that the kernel matrix built with recipe in equation (\ref{eq:gramdis}) is the polar form of $q$, which permits to build $q$ knowing $d$. Let us denote $G= (g_{ij})_{i,j}$ with $1 \leq i,j \leq n$, and
\[
G(a_i,a_j) = g_{ij}.
\]
Then, $G$ can be computed from the distances $d(i,j)$ by
\begin{equation}
    g_{ij} = - \frac{1}{2}\, \left(d_{ij}^2 - \Delta_i - \Delta_j + \Delta\right)
\end{equation}
with
\[
\Delta_i = \frac{1}{n}\sum_kd_{ik}^2, \qquad \Delta_j = \frac{1}{n}\sum_kd_{kj}^2, \qquad \Delta = \frac{1}{n^2}\sum_{k,\ell}d_{k\ell}^2
\]
$G$ is a symmetric bilinear form. Let us assume for sake of simplicity that it is definite (the extension to the case where it is non definite, i.e. is not full rank, is straightforward). The quadratic form is given by
\[
q(x) = B(x,x), \qquad \forall \: x \in E
\]
or, if $x = \sum_i \, x_i\, a_i$,
\[
q(x) = \sum_{i,j=1}^n \, g_{ij}\, x_ix_j.
\]
Now that we have the quadratic form, we need to construct the embedding $\varphi$.

\paragraph{Building the quadratic embedding:} Let us now denote by $(\omega_k, z_k)$ the eigenpairs of $G$, with $\omega_k\in \R$ and $z_k \in E$, i.e. pairs $(\omega_k,z_k)$ with
\begin{equation}\label{eq:mds:svdG}
    Gz_k = \omega_k \, z_k.
\end{equation}
or
\[
G\, Z = Z \, \Omega
\]
($Z$ is the $n \times n$ matrix with eigenvectors as columns, and $\Omega$ is the $n \times n$ diagonal matrix with $\omega_k$ on the diagonal). We have
\begin{equation}
    G = Z \, \Omega \, Z^\t
\end{equation}
(indeed, if $G'=Z\Omega Z^\t$, we have $G'Z=Z\Omega=GZ$ as $Z^\t Z=\I_n$, and as $Z$ is invertible, $G'=G$). Let us separate the positive from the negative eigenvalues of $G$, i.e. denote
\begin{equation}
    \omega_1 \geq \ldots \geq \omega_p > 0 > \omega'_1 \geq \ldots \geq \omega_m'
\end{equation}
(we assume for sake of simplicity that 0 is not an eigenvalue of $G$) or
\[
\left\{
\begin{array}{lcl}
  \omega_k > 0 & \mbox{if} & k \leq p \\
  \omega_k < 0 & \mbox{if} & k > p
\end{array}
\right.
\qquad \mbox{with} \qquad \omega'_j = \omega_{p+j}
\]
Let us denote
\begin{equation}
    \omega_k = \sigma_k^2, \qquad \omega'_j = - \theta^2_j
\end{equation}
and
\[
\begin{cases}
G \, u_k &= \sigma_k^2\, u_k \\
G \, v_j &= -\theta_j^2\, v_j
\end{cases}
\]
Let us use blockwise notations 
\begin{equation}
    Z = [U,V], \qquad \Omega = 
    \begin{pmatrix}
        \Sigma^2 & 0 \\
        0 & -\Theta^2
    \end{pmatrix}
\end{equation}
Then
\[
\begin{cases}
G \, U &= U \, \Sigma^2 \\
G \, V &= - V \, \Theta^2
\end{cases}
\]
Rewriting equation (\ref{eq:mds:svdG}) $G=Z\Omega Z^\t$ with this blockwise decomposition leads to 
\begin{equation}
    G =  U \, \Sigma^2 \, U^\t - V \, \Theta^2 \, V^\t
\end{equation}
illustrated by
\begin{center}
\begin{tikzpicture}
  \node () at (-.5,1) {$G=$} ; 
  \draw (0,0) rectangle (2,2) ; 
  \draw (1,0) -- (1,2) ;
  \node () at (0.5,1) {$U$} ;
  \node () at (1.5,1) {$V$} ;  
  \draw (2.5,2.5) rectangle (4.5,4.5) ; 
  \draw (3.5,2.5) -- (3.5,4.5) ;
  \draw (2.5,3.5) -- (4.5,3.5) ; 
  \node () at (3,4) {$\Sigma^2$} ;
  \node () at (4,3) {$-\Theta^2$} ; 
  \node () at (3,3) {$0$} ;
  \node () at (4,4) {$0$} ;  
  \draw (2.5,0) rectangle (4.5,2) ; 
  \draw (3.5,0) -- (3.5,2) ;
  \node () at (3,1) {$U\Sigma^2$} ;
  \node () at (4,1) {$-V\Theta^2$} ; 
  \draw (5,2.5) rectangle (7,4.5) ; 
  \draw (5,3.5) -- (7,3.5) ;
  \node () at (6,4) {$U^\t$} ;
  \node () at (6,3) {$V^\t$} ;
  \draw (5,0) rectangle (7,2) ;
  \node () at (6,1.3) {$U \, \Sigma^2 \, U^\t$} ; 
  \node () at (6,0.7) {$-V \, \Theta^2 \, V^\t$} ;   
\end{tikzpicture} 
\end{center}
Let us denote
\begin{equation}
    \begin{cases}
      X &= U \,\Sigma \\
      Y &= V \, \Theta
    \end{cases}
\end{equation}
Then
\begin{equation}
    G = XX^\t - YY^\t
\end{equation}
One recovers in $X$ the components obtained by clipping to 0 the eigenvalues $\omega'_j$, and complement this classical result with a second point cloud associated with the negative eigenvalues of the Gram matrix.

\paragraph{Algorithm:} This leads to the following algorithm (for sake of clarity, we denote here by $\lambda$ the positive eigenvalues of $G$ and by $\psi$ the negative ones, what was denoted $\omega$ and $\omega'$ in the text. \\
\\
\begin{algorithm}[H]
\begin{algorithmic}[1]
\STATE \textbf{input} $D \in \R^{n \times n} \: : \: D[i,j] = d_{ij}$
\STATE \textbf{compute} $\Delta_i = \frac{1}{n}\sum_kd_{ik}^2, \quad \Delta_j = \frac{1}{n}\sum_kd_{kj}^2, \qquad \Delta = \frac{1}{n^2}\sum_{i,j}d_{ij}^2$
\STATE \textbf{compute} $G \in \R^{n \times n} \: : \: G[i,j] = -\frac{1}{2}\left(d_{ij}^2 - \Delta_i - \Delta_j + \Delta\right)$
\STATE \textbf{compute} $(\omega_k,z_k) \: : \: Gz_k=\omega_kz_k$
\STATE \textbf{denote} $\Lambda, \Psi$, diagonal matrices of eigenvalues $\lambda >0$ and $ \psi<0$
\STATE \textbf{compute} $\Sigma=\Lambda^{1/2}, \: \Theta = (-\Psi)^{1/2}$
\STATE \textbf{denote} $U, V$, eigenvectors of $G$ with $\lambda >0$, and $ \psi < 0$
\STATE \textbf{compute} $X=U\Sigma, \: Y=V\Theta$
\RETURN $X,Y, \Sigma, \Theta$
\end{algorithmic}
\caption{Quadratic embedding of a discrete metric space}
\label{alg:quadratic_embedding}
\end{algorithm}

\paragraph{Summary and quality of the low rank approximation:} We have a metric space $(M,d)$ with $|M|=n$. There exists a quadratic space $(E,q)$ with $\dim E=n$ and a map
\[
\begin{CD}
M @>\varphi>> E
\end{CD}
\]
with $a_i=\varphi(i)$ such that
\begin{equation*}
    d^2(i,j) = q(a_i-a_j).
\end{equation*}
Let $G$ be the polar form of $q$, i.e.
\[
g_{ij} = G(a_i,a_j) = \frac{1}{2} (q(a_i+a_j) - q(a_i) - q(a_j)), 
\]
It is called the kernel matrix of $q$, and can be computed knowing the distances by 
\begin{equation*}
    g_{ij} = - \frac{1}{2}\, \left(d_{ij}^2 - \Delta_i - \Delta_j + \Delta\right).
\end{equation*}
with
\[
\Delta_i = \frac{1}{n}\sum_jd_{ij}^2, \qquad \Delta_j = \frac{1}{n}\sum_id_{ij}^2, \qquad \Delta = \frac{1}{n^2}\sum_{i,j}d_{ij}^2
\]
The quadratic form $q$ is defined by its polar form: $q(x)=G(x,x)$. The signature of $q$ is $(p,m)$ with $p$ being the number of positive eigenvalues of $G$ and $m$ the number of negative eigenvalues. We denote by $E_+$ (resp. $E_-$) the subspace of $E$ spanned by the eigenvectors of $G$ associated to a positive (resp. negative) eigenvalue of $G$. Then, one can write
\[ 
E = E_+ \oplus E_-
\]
and, for any point $a_i \in E$,
\[
a_i = x_i \oplus y_i, \qquad \mbox{with} \quad 
\begin{cases}
x_i & \in E_+ \\
y_i & \in E_-.
\end{cases}
\]
Two point clouds $X$ and $Y$ are built with algorithm \ref{alg:quadratic_embedding}, with $X$ being made of $p$ points in $\R^p$ and $Y$ of $m$ points in $\R^n$ such that
\[
G = XX^\t - YY^\t.
\]
As $G$ is symmetric, its eigenvectors are orthogonal. The columns of $X$ (resp. $Y$) are the eigenvectors associated with positive (resp. negative) eigenvalues of $G$. Then $X^\t Y = \zeros$. Hence
\begin{equation}
    \begin{array}{lcl}
         \|G\|^2 &=& \|XX^\t - YY^\t\|^2 \\
         &=& \|XX^\t\|^2 + \|YY^\t\|^2 - 2 \langle XX^\t \, , \, YY^\t \rangle \\
         &=& \|XX^\t\|^2 + \|YY^\t\|^2 
    \end{array}
\end{equation}
because
\begin{equation}
    \begin{array}{lcl}
         \langle XX^\t \, , \, YY^\t \rangle &=& \Tr \{XX^\t (YY^\t)^\t\} \\
         &=& \Tr \{XX^\t \, YY^\t\} \\
         &=& \Tr \{X \, (X^\t Y) \, Y^\t\} \\
         &=& 0
    \end{array}
\end{equation}
If negative eigenvalues are clipped to 0, one has $G \approx XX^\t$, and the quality of representation of $G$ by $XX^\t$ is $\|XX^t\|/\|G\|$. Similarly, knowing that 
\[
    g_{ij} = \langle x_i \, , \, x_j\rangle - \langle y_i \, , \, y_j\rangle,
\]
(this is another formulation of $G=XX^\t - YY^\t$), one has, 
\begin{equation}
    \begin{array}{lcl} 
    d^2(i,j) &=& q(a_i-a_j) \\
    &=& G(a_i-a_j \, , \, a_i-a_j) \\
    &=& G(a_i,a_i) + G(a_j,a_j) - 2 G(a_i,a_j) \\
    &=& g_{ii} + g_{jj} - 2 \, g_{ij}\\
    &=& \|x_i\|^2 - \|y_i\|^2 + \|x_j\|^2 - \|y_j\|^2 - 2 \langle x_i \, , \, x_j\rangle + 2  \langle y_i \, , \, y_j\rangle \\
    &=& \|x_i-x_j\|^2 - \|y_i-y_j\|^2.
    \end{array}
\end{equation}
So, $\|y_i-y_j\|^2$ quantifies the discrepancy between $d^2(i,j)$ and $\|x_i-x_j\|^2$ in classical MDS (where negative eigenvalues and eigenectors are clipped to 0) with full rank. This shows as well that distances in the point cloud built by classical MDS with clipping negative eigenvalues to 0 are overestimated as $\|x_i-x_j\|^2 = d^2_{ij} + \|y_i-y_j\|^2 \geq d^2_{ij}$.

\notes See appendix \ref{sec:quadratic} for a short presentation of quadratic forms and spaces. The classification of quadratic forms when $\K = \R$ or $\C$ depends on the matrix of the polar form and is well understood and given for example in \cite{SeguinPazzis2010}. The classification on a arbitrary field is said to be immensely difficult \cite{Berger1987a}. Embedding from a metric space into a quadratic space have been studied e.g. in \cite{Goldfarb1984, Gower1985} and is presented in \cite[section 3.5]{Pekalska2005} where it is called \emph{pseudo-euclidean embedding}.

%
\clearpage
\section{Summary}\label{sec:summary}
%

Here is a summary of the methods with 
\begin{itemize}
 \item the name of the function
 \item the call of the function
 \item the calculations to produce the result
\end{itemize}

\nB

\begin{center}
\ovalbox{
\begin{tabular}{cll}
\textbf{Method} & \textbf{Call} & \textbf{Computation (simplified)}\\
& & \\
\PCA & $(Y,V,\Lambda)=\textsc{pca\_core}(A)$ & $ C = A^\t A$\\
& & $(\Lambda, V) = \textsc{eig}(C)$ \\
& & $Y=AV$ \\
& & or \\
& & $(U,\Sigma,V) = \textsc{svd}(A)$\\
& & $\Lambda=\Sigma^2$ \\
& & $Y=U\Sigma$ \\
& & \\
\PCAiv & $(Y,V,\Lambda)=\textsc{pca-iv}(A,U,V)$ & $\PC_E=U(U^\t U)^{-1}U^\t$ \\
& & $\PC_F=V(V^\t V)^{-1}V^\t$ \\
& & $A_{\textsc{u},\textsc{v}} = \PC_\textsc{e}A\PC_\textsc{f}$ \\
& & $(Y,\Lambda,V)=\textsc{pca\_core}(A_{\textsc{u},\textsc{v}} )$ \\
& & \\
\PCAmet & $(Y,V,\Lambda)=\textsc{pca-met}(A,M,Q)$ & $A_{\textsc{m,q}}=MAQ$\\
& & $(Z,\Lambda,X)= \textsc{pca\_core}(A_{\textsc{m,q}})$ \\
& & $Y=M^{-1}Z$ \\
& & $V=Q^{-1}X$ \\
& & \\
\CoA & $(Y_r,Y_c,\Lambda ) = \CoA(T)$ & $A = T/T_{++}, \quad T_{++}=\sum_{i,j}T_{ij}$\\
& & $r_i = \sum_j\alpha_{ij}, \quad c_j = \sum_i\alpha_{ij}$ \\
& & $M = \diag (1/\sqrt{r_i}), \quad Q = \diag (1/\sqrt{c_j})$ \\
& & $A_{\m,\q}= M(A-r \otimes c)Q = \frac{\alpha_{ij}-r_ic_j}{\sqrt{r_ic_j}}$ \\
& & $Z,\Lambda,X = \textsc{pca\_core}(A_{\m,\q})$\\
& & $Y_r = M^{-1}Z, \quad Y_c = Q^{-1}X$ \\
& & \\
\CCA & $(Y_\a, Y_\b, U_\a, U_\b, \Lambda) = \textsc{cca}(A,B)$ & $T=A^\t B$ \\
& & $M = (A^\t A)^{-1/2}$, \quad $Q = (B^\t B)^{-1/2}$ \\
& & $R = MTQ$ \\
& & $W_\a, \Psi, \texttt{}W_\b = \textsc{pca\_core}(R)$ \\
& & $\Lambda = \sqrt{\Psi}$ \\
& & $U_\a = MW_\a, \quad U_\b=QW_\b$ \\
& & $Y_\a= AU_\a, \quad Y_\b=BU_\b$ \\
& & \\
\MDS & $Y,\Lambda = \textsc{mds}(D,r)$ & $G = \textsc{gram}(D)$ \\
& & $(U,\Sigma) = \textsc{svd}(G)$: $G=U\Sigma U^\t$ \\
& & $Y=U\Sigma^{1/2}$ \\
\end{tabular}
}
\end{center}

%
\newpage
\section{References in textbooks}\label{sec:reftextbooks}
%

In this section, we indicate where, and under which name, some techniques are presented in most classical textbooks.

\nB

\begin{center}
\ovalbox{
 \begin{tabular}{lll}
  Canonical Correlation Analysis & \cite{Anderson1958} & chapter 12 \\
  & \cite{Rao1973} & section 8f1 \\
  & \cite{Mardia1979} & chapter 10 \\
  & \cite{Jolliffe2002} & section 9.3 \\
  & \cite{Izenman2008} & section 7.3 \\
  & \cite{Hastie2009} & section 14.5.1 \\
  & \cite{Murphy2012} & section 12.5.3 \\
  Correspondence Analysis & \cite{Greenacre1984} & the book \\
  & \cite{Izenman2008} & chapter 17 \\
  Factor Analysis & \cite{Anderson1958} & section 14.7 \\
  & \cite{Rao1973} & section 8f4 \\
  & \cite{Mardia1979} & chapter 9 \\
  & \cite{Izenman2008} & secion 15.4 \\
  & \cite{Hastie2009} & section 14.7.1 \\
  & \cite{Murphy2012} & section 12.1 \\
  Independent Component Analysis & \cite{Izenman2008} & section 15.3 \\
  & \cite{Hastie2009} & section 560 \\
  Karhunen Loeve transform & \cite{Murphy2012} & section 12.2, p. 387 \\
  Latent variables & \cite{Bishop2006} & chapter 12 \\
  & \cite{Izenman2008} & chapter 15 \\
  & \cite{Hastie2009} & section 14.7.1 \\
  & \cite{Murphy2012} & chapter 12 \\
  Multiple Correspondence Analysis & \cite{Izenman2008} & section 17.4 \\
  Non metric scaling & \cite{Cox2001} & chapter 3 \\
  & \cite{Izenman2008} & section 13.9 \\
  Principal Component Analysis & \cite{Anderson1958} & chapter 11 \\
  & \cite{Mardia1979} & chapter 8 \\
  & \cite{Jolliffe2002} & the book \\
  & \cite{Bishop2006} & section 12.1 \\
  & \cite{Izenman2008} & section 7.2 \\
  & \cite{Murphy2012} & section 12.2 \\
  Probabilistic PCA & \cite{Bishop2006} \\
  & \cite{Murphy2012} & section 12.2.4 \\
  Sensible PCA & \cite{Murphy2012} & section 12.2, p. 387 \\
  Sparse PCA & \cite{Hastie2009} & section 14.5.5 \\
  Supervised PCA & \cite{Murphy2012} & section 12.5.1 \\
 \end{tabular}
 }
\end{center}

\begin{center}
\ovalbox{
 \begin{tabular}{lll}
  Classical Scaling & \cite{Cox2001} & section 2.2 \\
  & \cite{Izenman2008} & section 13.6 \\ 
  Multidimensional Scaling & \cite{Cox2001} & the book \\
  & \cite{Mardia1979} & chapter 14 \\
  & \cite{Izenman2008} & chapter 13 \\
  & \cite{Hastie2009} & section 14.8 \\ 
  Non metric scaling & \cite{Cox2001} & chapter 3 \\
  & \cite{Izenman2008} & section 13.9 \\
 \end{tabular}
 }
\end{center}

\section*{Abbreviations}

\begin{center}
 \ovalbox{
  \begin{tabular}{cl}
   CCA & Canonical Correlation Analysis \\
   CoA & Correspondance Analysis \\
   FA & Factor Analysis \\
   GRP & Gaussian Random Projection \\
   IV & Instrumental Variables \\
   MCoA & Multiple Correspondence Analysis \\
   MDA & Multivariate Data Analysis \\
   MDS & Mulidimensional Scaling \\
   PCA & Princial Component Analysis \\
   PPCA & Probabilistic PCA \\
   RP & Random Projection \\
   rSVD & Randomized SVD \\
   SDP & Symmetric Definite Positive \\
   SGF & Symmetric Gauge Function \\
   SVD & Singular Value Decomposition \\
   SVM & Support Vector Machine \\
   UIN & Unitarily Invariant Norm \\
  \end{tabular}
 }
\end{center}

%

\appendix
\section{Preliminaries in linear algebra}\label{sec:la}
%

Let us recall here some basic facts in linear algebra. Linear algebra can have an abstract setting, with the notions of vector spaces and linear maps, or numerical setting, working with arrays (vectors and matrices). A matrix $A$ is the expression of a linear map $E \longrightarrow F$ by an array once basis have been selected in $E$ and $F$.  

%
\subsection{Vector space and linear map}
%

\nD Let $\K$ be a field, usually $\R$ or $\C$, and here $\R$ unless otherwise stated. A vector space $E$ over $\K$ is a set endowed with two operations:
\begin{itemize}[label=$\rightarrow$]
\item an addition, $+$, such that $(E,+)$ is an Abelian group
\item a multiplication by a scalar
\[
\begin{CD}
\K \times E @>>> E \\
(\alpha,u) @>>> \alpha u. 
\end{CD}
\]
\end{itemize}
Multiplication by a scalar verifies the following properties:
\[
\begin{array}{|ccc}
   \alpha (u + v) &=& \alpha u  + \alpha v \\
   (\alpha + \beta) u &=& \alpha u + \beta u \\
   \alpha(\beta u) &=& (\alpha\beta)u \\
   1 u &=&  u.
\end{array}
\]

\nD A vector subspace of $E$ is a subset of $E$ which is closed under the addition or the multiplication by a scalar. It is itself a vector space on $\K$. 

\nD A basis of a vector space $E$ is a collection $(u_k)_k$ of vectors $u_k \in E$ such that any vector $u \in E$ can be written as
\[
u = \sum_k \alpha_k u_k,
\]
in a unique manner where $\alpha_k \in \K$. If a basis exists, all basis have the same cardinality. This cardinality is the dimension of the vector space. Usually, the dimension is either an integer $n \in \N$ or $\aleph_0$, the cardinality of $\N$. The vector spaces we will work with in these notes are finite dimensional. 

\nD Let $E,F$ be two vector spaces on the same field $\K$. A linear map 
\[
\begin{CD}
   E @>L>> F,
\end{CD}
\]
from $E$ on $F$ is a map which preserves the linear structure, i.e. 
\[
\begin{array}{ccc}
     L(u+v) &=& Lu + Lv \\
     L(\alpha u) &=& \alpha \, Lu
\end{array}
\]
(it is customary for linear maps to write $Lu$ instead of $L(u)$ when there is no ambiguity). The image of a linear subspace $E' \subset E$ of $E$ is a vector subspace $F' \subset F$.

\nD Let $n = \dim E \in \N$ and $m = \dim F \in \N$. Let $(u_j)_j$ with $1 \leq j \leq n$ be a basis of $E$ and $(v_i)_i$ with $1 \leq i \leq m$ a basis in $F$. Let
\[
Lu_j = \sum_i \alpha_{ij} v_i.
\]
Then, the matrix $(\alpha_{ij})_{i,j} \in \R^{m \times n}$ is the matrix of the linear map $L$ in basis $(u_j)_j$ for $E$ and $(v_i)_i$ for $F$. We have, for any vector $u = \sum_j\beta_j u_j \in E$
\begin{equation}
   \begin{array}{lcl}
      Lu &=& \displaystyle L\left(\sum_j\beta_j u_j\right) \\
      &=& \displaystyle \sum_j \beta_j Lu_j \\
      &=& \displaystyle \sum_j \beta_j \left(\sum_i \alpha_{ij} v_i\right) \\
      &=& \displaystyle \sum_i\left(\sum_j \alpha_{ij}\beta_j\right)v_i.
   \end{array}
\end{equation}
Hence, the $j-$th column of $A$ is the vector $Lu_j \in \R^m$ with coordinates $(\sum_{\ell=1}^n\alpha_{i\ell}\beta_\ell)_{1 \leq i \leq m}$.

\nD The kernel of $L$ is the set of vectors $u \in E$ the image of which is $0 \in F$
\begin{equation}
    \ker L = \{u \in E \mid Lu=0\}.
\end{equation}
It is a vector subspace of $E$.

\nD The image of $L$ is the set of vectors $v \in F$ which have a preimage in $E$
\begin{equation}
    \im L = \{v \in F \mid \exists \: u \in E \quad \st \quad v=Lu\}.
\end{equation}
It is a vector subspace of $F$. 

\nD The rank nullity theorem states that
\begin{equation}
    \dim \ker L + \dim \im L = \dim E
\end{equation}

\nD A linear map from $E$ to $E$ is called an endomorphism. If it is an isomorphism, it is called an automorphism. The set of all endomorphisms of $E$, denoted $\mathrm{end\:}E$, is an associative algebra for the operation $+$ and multiplication by a scalar for the vector space structure, and the composition $\circ$ for it to be an algebra.

%
\subsection{Eigenspace, eigenvector, eigenvalue}
%

\nD Let  $E$ be a finite dimensional vector space on a field $\K$, and $L$ an endomorphism in $E$
\[
\begin{CD}
    E @>L>> E.
\end{CD}
\]
An eigenvector of $L$ is a vector $u \neq 0 \in E$ such that
\begin{equation}
    Lu = \lambda u,
\end{equation}
for some $\lambda \in \K$. $\lambda$ is called an eigenvalue of $L$. If $A$ is the matrix of $L$ in a given basis, one writes as well $Ax=\lambda x$ if $x$ is the expression of $u$ in the same basis, and $(x,\lambda)$ is an eigenpair of $A$. The set of eigenvalues of $A$ is called the spectrum of $A$:
\[
\Sp A = \{\lambda \in \K \mid \exists \: x \neq 0 \quad \st \quad Ax=\lambda x\}.
\]
The eigenspace of $A$ associated to eigenvalue $\lambda$ is the set of all eigenvectors associated $\lambda$. This space completed with $\{0\}$ is the kernel of $A-\lambda \I$. In order to avoid technicalities, one includes $\{0\}$ in the eigenspace associated to an existing eigenvalue. The eigenspace is then $\ker (A-\lambda \I)$, and is a linear subspace of $E$.  

\nD Even if matrix $A$ is real, its eigenvalues can be complex. For example, if
\[
A = \begin{pmatrix}
    0 & 1 \\
    -1 & 0
\end{pmatrix},
\]
we have
\[
\left\{
  \begin{array}{lcl}
    z &=& \lambda y \\
    -y &=& \lambda z,
  \end{array}
\right.
\]
hence $z = -\lambda^2 \, z$, so $\lambda^2=-1$ and $\lambda = \pm i \in \C$. The eigenvalues of $A$, if they exist, are the root of the characteristic polynomial of $A$
\begin{equation}
    \lambda \in \Sp A \quad \Longrightarrow \quad \det (A - \lambda \I)=0,
\end{equation}
where $\I$ is the identity matrix. The roots of a real polynomial can be complex.

\nD A square matrix $A$ acting on $E$ is diagonalisable if $E$ has a basis of eigenvectors of $A$. If $A \in \K^{n \times n}$, the sum of the dimensions of its eigenspaces is $n$. Not all matrices are diagonalisable. For example, if
\[
A = \begin{pmatrix}
    0 & 1 \\
    0 & 0
\end{pmatrix},
\qquad 
x = \begin{pmatrix}
    y \\
    z 
\end{pmatrix},
\]
$Ax=\lambda x$ leads to
\[
\left\{
  \begin{array}{lcl}
    0.y + z &=& \lambda y \\
    0.y + 0.z &=& \lambda z,
  \end{array}
\right.
\]
then $\lambda=0$, $z=0$ and $y \in \K$. The vector space spanned by the eigenvectors of $A$ has dimension 1, and $A$ is not diagonalisable. A matrix which is not diagonalisable is called defective. A nilpotent matrix is a matrix $A$ such that there exists an integer $m>0$ with $A^m=0$. A nilpotent matrix  is defective (except the zero matrix $\mathbf{0}$, because any vector of $E$ is an eigenvector of $\mathbf{0}$ associated to eigenvalue $\lambda=0$).

\nD A matrix $A$ is diagonalizable, or non defective, if there exists an invertible matrix $P$ such that $A=P\Lambda P^{-1}$ where $\Lambda$ is a diagonal matrix (all matrices in this formula are in $\K^{n \times n}$). $A=P\Lambda P^{-1}$ is called the eigendecomposition of $A$. The columns of $P$ are eigenvectors of $A$. It can be seen through $AP=P\Lambda$. If $A$ is real symmetric, its eigendecomposition exists, with $P$ orthogonal $(P^{-1}=P^\t)$, and its spectrum is real.

\subsection{Perturbation of eigenvalues}

A matrix $A$ being given, there are two sources of errors when commputing its eigenvalues:
\begin{itemize}[label=$\rightarrow$]
\item a numerical error while using rounding during the calculation: this is addressed by numerical analysis of eigendecomposition;
\item when the matrix is the outcome of an experiment, i.e. a dataset, the data can be corrupted, which leads to a corruption of the eigenvalues as well.
\end{itemize} In this section, we address the second source of errors: possible corruption of a dataset (which is dealing with uncertainty rather than error). It will be referred to as a perturbation. The theory of perturbations of the eigenvalues of a given matrix is the study of the variations of the eigenvalues of a given matrix under a perturbation of its coefficients.\\
\\
We know by Rouché theorem that the eigenvalues of a matrix $A$ are continuous functions of the coefficients of $A$. Pointedly speaking, let $A, H \in \K^{n \times n}$ where $\K$ is $\R$ or $\C$, and $\epsilon \in \R$. What can be said
on the localization of eigenvalues of $A + \epsilon H$ knowing the spectrum of $A$ ? A first observation is that the eigenvalues of $A$ are the roots of a 
polynomial (the characteristic polynomial, of degree $n$ if $\dim A=n$). As, in general, the roots of a polynomial as functions of its coefficients are unstable, it is likely that the eigenvalues of a matrix are unstable.
This will be shown, but a good news is that the eigenvalues of a symmetric matrix are stable under a perturbation by a symmetric matrix, i.e. vary with the same order of magnitude that the coefficients of the matrix. This is Weyl's theorem which dates back to 1912. \\
\\
Let us have (this example is borrowed from \cite[chapter IV]{Stewart1990})
\[
 A = \begin{pmatrix}
      0 & 1 & 0 & 0 \\
      0 & 0 & 1 & 0 \\
      0 & 0 & 0 & 1 \\
      0 & 0 & 0 & 0
     \end{pmatrix}
\]
It is easy to show that $\Sp A = \{0\}$. Let us now have
\[
 A + \epsilon H = \begin{pmatrix}
      0 & 1 & 0 & 0 \\
      0 & 0 & 1 & 0 \\
      0 & 0 & 0 & 1 \\
      \epsilon & 0 & 0 & 0
     \end{pmatrix}
\]
Then, 
\begin{equation}
 \Sp (A + \epsilon H) = \{ \pm \epsilon^{1/4}, \pm i \epsilon^{1/4} \}
\end{equation}
This can be generalized to any dimension $n$, which shows that the perturbation of eigenvalues can be in $\epsilon^{1/n}$
if $A \in \M(n,n)$. If $n=100$ and $\epsilon = 10^{-2}$, then $\epsilon^{1/n} \approx 0.95$. If $\epsilon'$ is the magnitude of the perturbation of the eigenvalues, we have $\epsilon'/\epsilon = 0.95/10^{-2} \approx 95.49$. The perturbation is amplified about 100 times.\\
\\
The first tool needed is a way to compare the spectrum of the initial and the perturbed matrix. The tools therefore are the spectral variation, and distances, like Hausdorff distance and matching distance between spectra. They are presented hereafter. Then, some bounds are given on distances between spectra of $A$, $B$ and $A+B$, and used where $B$ is considered as a perturbation of $A$. This is a key question in numerical analysis which has been thoroughly studied (see references and notes section).\\
\\
\textbf{Spectral variation:} Let $A, B \in \M(n,n)$ with $\Sp A = \{\lambda_k\}_k$ and $\Sp B = \{\mu_k\}_k$ with $1 \leq k \leq n$. Then, the
spectral variation of $B$ relatively to $A$ is the number
\[
 \mbox{sv}_A(B) = \underset{i}{\max} \left(\underset{j}{\min} \: |\lambda_i-\mu_j|\right)
\]
It measures the maximum possible gap between an eigenvalue in $A$ and the closest eigenvalue in $B$.\\
\\
\textbf{Matching distance:} This gap is not a distance. However, a distance between spectra can be built as
\[
 d_H(\Sp A, \Sp B) = \max \{\mbox{sv}_A(B) \, , \, \mbox{sv}_B(A)\}
\]
It is the Hausdorff distance between the spectra of $A$ and of $B$. However, this distance is not term by term comparison of eigenvalues. Hence the matching distance is defined as
\begin{equation}
 \mbox{md}\:(\Sp A, \Sp B) = \underset{\pi}{\min} \left( \underset{i}{\max} |\lambda_i - \mu_{\pi(i)}|\right)
\end{equation}
where $\pi$ runs over all permutations of $\Sp B$.\\
\\
\textbf{Ostrowski and Elsener theorem:} It can be shown that
\begin{equation}\label{eq_elsner90}
 \mbox{md}\: (\Sp A, \Sp B) \leq 4 \times 2^{-1/n} \left(\|A\| + \|B\|\right)^{1-1/n} \|A-B\|^{1/n}
\end{equation}
or
\begin{equation}
 \mbox{md}\: (\Sp A, \Sp B) \leq 4 \times 2^{-1/n} \left(\|A\| + \|B\|\right)\left( \frac{\|A-B\|}{\|A\| + \|B\|}\right)^{1/n}
\end{equation}
where
\begin{equation}
 \|X\| = \sqrt{\lambda_{\max}(X^*X)} = \underset{\|u\|=1}{\max} \:Xu
\end{equation}
Hence, if $B = A + \epsilon H$, this yields
\begin{equation}
  \begin{aligned}
    \mbox{md}\: (\Sp A, \Sp B) & \leq 4 \times 2^{-1/n} \left(\|A\| + \|A+\epsilon H\|\right)
 \left( \frac{\|\epsilon H\|}{\|A\|+ \|A+\epsilon H\|}\right)^{1/n} \\
 & \leq C \epsilon^{1/n}
  \end{aligned}
\end{equation}
The bound $h =1/n$ in $\epsilon^h$ is reached in the example.\\
\\
\textbf{Henrici theorem:} There exists a general and powerful result for the localization of eigenvalues of a matrix
which is the sum of two symmetric matrices. Let $A,B,C \in \M(n,n)$, symmetric, such that
\begin{equation}
 B = A+C
\end{equation}
Let us denote
\begin{equation}
 \left\{
 \begin{array}{lcl}
   \Sp A &=& \{\alpha_i\} \\
   \Sp B &=& \{\beta_j\} \\
   \Sp C &=& \{\gamma_k\}
 \end{array}
 \right.
\end{equation}
Let us assume that
\begin{equation}
 \left\{
 \begin{array}{l}
   \alpha_1 \geq \ldots \geq \alpha_n \\
   \beta_1 \geq \ldots \geq \beta_n \\
   \gamma_1 \geq \ldots \geq \gamma_n \\  
 \end{array}
 \right.
\end{equation}
Then (recall that $C = B-A$)
\begin{equation}
 \forall \: i, \quad \gamma_n \leq \beta_i - \alpha_i \leq \gamma_1
\end{equation}

\noindent This theorem has a nice consequence for symmetric perturbation of eigenvalues of a symmetric matrix. Let us assume
that $C=\epsilon H$, or $B = A + \epsilon H$. Let us denote
\begin{equation}
 \overline{h} = \underset{i,j}{\mbox{sup}}\: |h_{ij}|
\end{equation}
Then
\begin{equation}
 \forall \: i, \quad |\gamma_i| \leq \epsilon \overline{h}
\end{equation}
and 
\begin{equation}
  \forall \: i, \quad |\beta_i - \alpha_i| \leq \epsilon \overline{h}
\end{equation}
Hence the response  of any eigenvalue of a symmetric matrix by a symmetric perturbation is bounded by a term linear
with the perturbation, and not in $\epsilon^{1/n}$ as for  any matrix or perturbation. Eigenvalues of symmetric matrices 
are much more stable under symmetric perurbations.  \\
\\
\textbf{Application to PCA:} Let us now have a matrix $A \in \R^{n \times p}$. PCA involves the computations of eigenvalues and eigenvectors of $C = A^\t A$.
Let us now have a small perturbation of $A$: $A \longrightarrow B=A+\epsilon H$. Then $C \longrightarrow C'$ with
\begin{equation}
 \begin{aligned}
   C' &= B^\t B \\ 
   &= (A^\t + \epsilon H^\t)(A + \epsilon H) \\
   &= C + \epsilon (H^\t A + A^\t H) + \epsilon^2 \; H^\t H \\
   &= C + \epsilon (H^\t A + A^\t H + \epsilon \;  H^\t H)
 \end{aligned}
\end{equation}
with $H^\t A + A^\t H + \epsilon \; H^\t H$ being symmetric. Then, $C'$ is obtained from $C$ by a symmetric perturbation,
and Henrici theorems apply. The variation of the eigenvalues is bounded by a term linear with $\epsilon$.\\
\\
This can be derived directly from Weyl's 1912 theorem (see references and notes), knowing that $\forall \: i, \quad
|\lambda_i| \leq \underset{i,j}{\sup}\: |\alpha_{ij}|$.

\notes Here, we have followed \citep[chapter IV]{Stewart1990}, with historical notes p. 176 \& sq. $\,$.  Perturbation theory of eigenvalues of matrices or linear operators has a long history, from a result by H. Weyl in 1912. It states that if $A,B$ are self-adjoint matrices (a self-adjoint matrix is a matrix $A \in \M_\C(n,n)$ with $A=A^*$ where $A^* = \overline{A^\t}$) with spectra $\Sp A = \{\alpha_i\}$ and $\Sp B = \{\beta_j\}$, with $\alpha_1 \geq \ldots \geq \alpha_n$ and $\beta_1 \geq \ldots \geq \beta_n$,
then $\underset{k}{\max}\: |\alpha_k-\beta_k| \leq \|A-B\|_\sp$, where the norm is the spectral norm: $\|A\|_\sp = \underset{\|x\|=1}{\max} \|Ax\|$. Several decades or efforts have aimed at finding similar bounds in more general situations, i.e. non self-adjoint matrices or other norms. Much work has concerned so called normal matrices (a normal matrix is a matrix $A$ such that $AA^*=A^*A$). A matrix is normal if, and only if, it is diagonal in some orthogonal basis. As eigenvalues of normal matrices can be complex, they cannot be ordered as in $\R$. The notion of matching distance permits to compare spectra with complex values. A result to the question whether $\mbox{md}\: (\Sp A, \Sp B) \leq \|A-B\|$ for normal matrices have been object of intensive researches for decades. Hoffman and Wielandt have proved in 1953 (40 years after Weyl's result) a similar result for normal matrices, but where the norm is Frobenius norm $\|A\|_F= \sqrt{\Tr A^*A} = \left(\sum_{i,j}|\alpha_{ij}|^2\right)^{1/2}$: $\mbox{md}_F\: (\Sp A, \Sp B) \leq \|A-B\|_F$, where $\mbox{md}_F$ is a matching distance adapted to Frobenius norm ($\ell^2$ norm). See Bathia \citep{Bhatia2007} and \cite{Stewart1990} for historical notes, from which those few milestones have been borrowed.\\
\\
Due to its role in numerical analysis, spectral variation has been thoroughly studied over several decades (see e.g. Henrici \citep{Henrici1962}, Bhatia \citep{Bhatia1982}, Elsner \citep{Elsner1982}). Classical books on bounds for eigenvalue perturbation theory are \citep{Bhatia1987,Stewart1990}. The inequality (\ref{eq_elsner90}) appears in Bhatia \citep{Bhatia1990}. Several recent results for upper bounds of matching distance between two spectra appear in Galant\'ai \citep{Galantai2008}. Ostrowski theorem dates from 1940, and has been published in Ostrowski, A. (1940) Recherches sur la m\'ethode de Gr\"affe et les z\'eros des polyn\^omes et des series de Laurent. \emph{Acta Math.}, \textbf{72:}99-257. (see \cite{Holbrook1992}).

%
\section{Quadratic forms}\label{sec:quadratic}
%

\subsection{Quadratic and polar forms}

Let $B$ be a definite bilinear for on a finite dimensional vector space $E$. Then, the map
\[
\begin{CD}
E @>q>> \K \\
\end{CD}
\]
with
\[ 
q(x)=B(x,x)
\]
is called a quadratic form. We have
\[
B(x,y)=\frac{1}{2}(q'x+y)-q(x)-q(y))= \frac{1}{4}(q(x+y) - q(x-y))
\]
and $B$ is called the polar form of the quadratic form $q$. Let $\K=\R$, and $q$ be a quadratic form on $E$ with $\dim E=n$. 

\subsection{Signature of a quadratic form}

There exists a basis in $E$ in which the matrix of the polar form $B$ of $q$ is
\begin{equation}
    B = 
    \begin{pmatrix}
        \I_p & 0 & 0 \\
        0 & -\I_m & 0 \\
        0 & 0 & 0
    \end{pmatrix}
\end{equation}
The pair $(p,m)$ is called the signature of the quadratic form, and does not depend on the basis (Sylvester law of intertia). It is classical to denote $\R^n_{p,m}$ or $\R_{p,m}$ the space $\R^n$ equipped with the quadratic form
\[
q(x) = \sum_{i=1}^p x_i^2 - \sum_{j=p+1}^{n'}x_j^2
\]
with $n'=p+m \leq n$.

\subsection{Geometry in quadratic spaces}

Euclidean spaces are those for which $p=n$. Non Euclidean quadratic spaces $(p < n)$ are called pseudo-euclidean spaces.  In a Euclidean space, the map
\[
\begin{CD}
(x,y) @>>> d(x,y) = \sqrt{q(x-y)}
\end{CD}
\]
defines a distance. It is no longer the case in a pseudo-euclidean space.  For example, if $n=2, p=m=1$, the pseudo-sphere of radius 1 is defined by $x^2-y^2=1$ and is an hyperbola (hence the name hyperbolic geometry for geometry in pseudo-euclidean spaces).

%
\section{Random projection}\label{sec:random_projection}
%

Vector spaces of very large dimension exhibit properties we are not familiar with from our geometric intuition developed in $\R^2$ or $\R^3$. For our purpose of linear dimension reduction, this leads to the approaches relying on so called ``random projection'', which is gently introduced here. This is only the tip of an iceberg, addressing issues in measure concentration and geometry of Banach spaces. We restrict ourselves here on what is focused towards linear dimenson reduction, i.e. SVD with random projection.

%
\subsection{Isometries, orthogonal matrices and rotations}
%

Random projection consists in selecting randomly a subspace of small dimension, say $k$, in a space of large space, say $p$, and project a point cloud on it. One question is to characterize and select randomly such a subspace. It can be done by selecting an orthonormal basis of it, by selecting vectors with a uniform measure on the sphere (the Haar measure) with constraints of orthogonality. This boils down to selecting randomly a rotation, i.e. a $p \times p$ matrix in special orthogonal group, and keeping the first $k$ columns only.\\
\\
Let $n \in \N$. The set of invertible matrices in $\R^{n \times n}$ is the linear group of $\R^n$, denoted $\GL(n)$.

\nD A matrix $A \in \R^{n \times n}$ is orthogonal if
\begin{equation}
    AA^\t = A^\t A = \I_n.
\end{equation}
It is then invertible, and
\begin{equation}
    A^{-1}= A^\t.
\end{equation}
The set of all orthogonal matrices on $\R^n$ is called the orthogonal group, and is denoted $\O(n)$:
\begin{equation}
    \O(n) = \{A \in \R^{n \times n} \mid AA^\t = A^\t A = \I_n\}.
\end{equation}
It is one of the classical compact Lie group in $\GL(n)$
We have
\begin{equation}
    \det A = \pm 1.
\end{equation}
The Special Orthogonal Group $\S\O(n)$ is the set of the orthogonal matrices with determinant equal to one
\begin{equation}
    \S\O(n) = \{A \in \O(n) \mid \det A = 1\}
\end{equation}
Its elements are called rotations. The set of orthogonal matrices such that $\det A=-1$ is often denoted $\S\O^-(n)$.

\nD A matrix $A$ is orthogonal if, and only if, 
\begin{itemize}[label=$\rightarrow$]
\item its columns form an orthonormal basis of $\R^n$,
\item it acts as an isometry on $\R^n$, i.e.
\begin{equation}
    \forall \: x,y \in \R^n, \quad \langle Ax, Ay\rangle = \langle x, y\rangle.
\end{equation}
\end{itemize}

%
\subsection{Concentration of the measure on the sphere}
%

Concentration of the measure is an unexpected result in spaces of very large dimension about the global variations of a function whose local variations are kept small. The control of local variations is given by Lipschitz property. It is an immense domain, which is merely touched here for measure concentration on the sphere which will be useful for a proof of Johnson-Lindenstrauss lemma.

\paragraph{Lipschitz function:} Let $E, F$ be two metric spaces (here, we will be concerned by distances being associated to a norm in a finite dimensional vector space). A function
\[
\begin{CD}
  E @>f>> F
\end{CD}
\]
is said Lipschitz if there exists a constant $C$ such that
\begin{equation}
    \forall \: x,y \in E, \quad d(f(x)-f(y)) \leq C d(x,y).
\end{equation}
A linear map is a Lipschitz function. Indeed, Let $L \in \L(E,F)$ be a linear map. The spectral norm $\|.\|_\sp$ is defined as
\begin{equation}
    \| L\|_\sp = \max_{x \neq 0} \, \frac{\|Lx\|}{\|x\|}
\end{equation}
where $\|.\|$ is the Frobenius norm. Hence, for all $x$,  $\|Lx\| \leq \|L\|_\sp\, \|x\|$. Letting $x \rightarrow  x-y$ leads to 
\begin{equation}
    \|Lx-Ly\| = \|L(x-y)\| \leq \| L\|_\sp \, \|x-y\|.
\end{equation}


\paragraph{Measure concentration on the sphere:} Let $\S^{n-1}$ denote the sphere in $\R^n$:
\begin{equation}
    \S^{n-1} = \{x \in \R^n \mid \|x\|=1\}.
\end{equation}
Let
\[
\begin{CD}
\S^{n-1} @>f>> \R
\end{CD}
\]
be a Lipschitz function with constant $C$. Let $m$ be the median of this function on the sphere (i.e. a value such that for a random vector $x$ on the sphere, the probability that $f(x) \geq m$ is larger or equal to $1/2$, as well as the probability that $f(x) \leq 1/2$). Let $x$ be a uniformly selected random vector on the sphere. Then, Levy's lemma, or measure concentration on the sphere, asserts that, for any $t>0$,
\begin{equation}
    \P(|f(x)-m| \geq tC) \leq 2 \, \exp -(n-2)t^2.
\end{equation}
If $n \rightarrow \infty$, the quantity on the r.h.s. $\rightarrow 0$. Hence, for very large dimensions, the function $f$ is essentially equal to its median.


\nD Here is an example. Let us select randomly a vector $a \in \S^{n-1}$. It will be called the ``north pole'' of the sphere. Let us define the map
\begin{equation}
    \begin{CD}
      \S^{n-1} @>f>> \R \\
      x @>>> \langle a,x\rangle
    \end{CD}
\end{equation}
It is a linear form, hence is Lipschitz. We have
\begin{equation}
    |f(x)-f(y)| = |\langle a,x\rangle - \langle a,y\rangle| = \langle a,x-y\rangle| \leq \|a\| \, \|x-y\| = \|x-y\|
\end{equation}
Hence the constant $C$ is $C=1$. The median is $m=0$. Indeed, let us consider the hyperplane $H$ orthogonal to $a$. All points on the sphere in the upper half-space defined by it (the part which contains $a$) are such that $f(x) \geq 0$, and those in the lower half-space which contains $-a$ are such that $f(x) \leq 0$ ($f(x)=0$ for the points on the equator, defined as $\S^{n-1} \cap H$). We then have
\begin{equation}
     \P(|f(x)| \geq t) \leq 2 \, \exp -(n-2)t^2.
\end{equation}
Then, $f(x)= \langle a,x\rangle$ is zero almost everywhere (recall that $m=0$), and nearly all points are close to the equator. However, the fraction of points out of the $t-$neighborhood of the equator becomes negligeable  for very large dimensions only. For this fraction to be lower than a given $\epsilon$ for a given $t$, one must have
\begin{equation}
    n > 2 + \frac{1}{t^2}\, \Log \frac{2}{\epsilon} \approx \frac{1}{t^2}\, \Log \frac{2}{\epsilon}
\end{equation}
Selecting $t=\epsilon= 10^{-2}$ yields $n \approx 5.3 \times 10^4$. 

\notes A very clear and accessible introduction to the concentration of measure is \cite{Ledoux2001}. A basic example already known to Borel and mentioned in the introduction of \cite{Ledoux2001} is the geometric interpretation of the law of large numbers, given as follows. Let us consider the hypercube $\K_n=[0,1]^n \subset \R^n$ and $H$ be its intersection with the hyperplane orthogonal to the diagonal from $(0, \ldots,0)$ to $(1, \ldots,1)$. Let $H_t$ be the set of points in $\K_n$ at distance $\leq t$ from  $H$. Then, if $\mu$ is the uniform measure in $\K_n$ (noticing that $\mu(\K_n)=1$), i.e. $d\mu = dx_1\ldots dx_n$, then $\mu(H_{t\sqrt{n}}) \rightarrow 1$ when $n \rightarrow \infty$. All points in $\K_n$ are concentrated in the $t-$neighborhood of $H$ for $n$ sufficiently large. Their projection on the diagonal is concentrated on the segment $[1/2 - t\sqrt{n}, 1/2 + t\sqrt{n}]$. This is the law of large numbers. 

%
\subsection{The Johnson-Lindenstrauss lemma}
%

Johnson-Lindenstrauss lemma is about a good surprise for linear dimension reduction: loosely speaking, for a random point cloud $\X$ of size $m$ in a large dimension space $\R^n$, and an accuracy $\epsilon$, there exists a dimension $k$ and a subspace $E \subset \R^n$ of dimension $k$ such that the distances between the projections of points in $\X$ on $E$ approximates the distances between the points in $\X$ with accuracy $\epsilon$. For $k$ to be significantly smaller that $n$ when $\epsilon$ is small, $n$ must be very, very large. 

\nD Let us first give some notations
\begin{itemize}[label=$\rightarrow$]
\item $n \in \N$ and $\epsilon \in [0, 1/2]$,
\item $\X = (x_i)_i$ is a point cloud with $1 \leq i \leq n$, $x_i \in \R^n$
\item Let
\[
k \geq \frac{4 \, \Log n}{\epsilon^2/2 - \epsilon^3/3}
\]
\end{itemize}
Then, there exists a linear map
\[
\begin{CD}
\R^n @>f>> \R^k
\end{CD}
\]
such that
\begin{equation}
    \forall \: i,j, \quad (1-\epsilon)\|x_i-x_j\| \leq \|f(x_i)-f(x_j)\| \leq (1+\epsilon)\|x_i-x_j\|.
\end{equation}
This is an amazing lemma, which tells that in very large dimensions $(n \gg 1)$, a cloud $\X$ of $n$ points in $\R^n$ is sharply concentrated on a vector subspace of dimension $k$.

\nD Here is a sketch of a classical demonstration. It is a consequence of the measure concentration on the sphere. Let us recall that a mapping
\[
\begin{CD}
\R^n @>f>>\R^k
\end{CD}
\]
is $L-$Lipschitz if
\begin{equation}
    \forall \: x, y \in \R^n, \quad \|f(x) - f(y)\| \leq L \, \|x-y\|.
\end{equation}
Let us consider the sphere
\[
\S^{n-1} = \{x \in \R^n \mid \|x\|=1\}
\]
and the mapping
\begin{equation}
    \begin{CD}
      \S^{n-1} @>>> \R \\
      (x_1, \ldots, x_n) @>>> \sqrt{x_1^2 + \ldots + x_k^2}.
    \end{CD}
\end{equation}
It is easy to see that $f$ is $1-$Lipschitz. Then, measure concentration of the sphere leads to
\begin{equation}
    \P(|f(x)-m| \geq t) \leq 2\, \mathrm{e}^{-nt^2/2},
\end{equation}
where $m$ is a ``suitable number'' $m=m(n,k)$ which satisfies to $m > \frac{1}{2}\sqrt{k/n}$ under some conditions on $n$ and $k$.

\nD There are two ways to read this lemma: $(i)$ either the subspace is fixed, as above, and $x$ is a random vector on the sphere, or $(ii)$ $x$ is fixed, and the $k-$dimensional subspace is randomly chosen. Selecting a subspace $E$ of dimension $k$ at random can be done by selecting at random a rotation $R \in \SO(n)$ with uniform distribution in $\SO(n)$. The prof of the flattening lemma follows. First select such a subspace $E$ at random. Then, it can be shown that for any pair $x,y \in \R^n$
\begin{equation}
    \left(1-\frac{\epsilon}{3}\right)m\|x-y\| \leq \|p(x)-p(y)\| \leq \left(1+\frac{\epsilon}{3}\right)m\|x-y\|
\end{equation}
is violated with probability at most $1/n^2$. Second, the flattening lemma follows (with some technicalities omitted here) from the observation that there are less than $n^2$ pairs of points. 

\notes Johnson-Lindenstrauss lemma is the starting point of many developments in the domain of dimension reduction and random algorithms in linear algebra. It is often referred to as the \emph{flattening lemma}, because it flattens a point cloud. The lemma given here has been borrowed from \cite{Matousek2008} and \cite{Dasgupta2003}. The demonstration is sketched in \cite{Matousek2008} and detailed in \cite[section 15.2]{Matousek2002} or \cite{Dasgupta2003}, which are similar and rely on the same observations, from which it has been borrowed and to which the reader may refer for all details omitted here. \cite{Matousek2002} gives as well some historical notes, and points to several surveys on JL lemma and its utilisation in a diversity of algorithms. Let us mention however that, for having $k < n$ with the given bound $\frac{4 \, \Log n}{\epsilon^2/2 - \epsilon^3/3}$ for small $\epsilon$, $n$  must be huge. For example, if $\epsilon = 10^{-2}$ and $n=10^6$, then $k>n$! If $n=10^7$, then $k\geq 1.28 \times 10 ^6$. For $\epsilon = 10^{-1}$, this drops to $1.2 \times 10^4$ for $n=10^7$, and $k \geq 8.6 \times 10^3$ for $n=10^5$. The threshold value of $k$ is very sensitive to $\epsilon$, as $k = \OC(\epsilon^{-2} \, \Log n)$. However, luckily, the projection is of very good quality for much lower dimensions $n$.

\newpage
\bibliographystyle{alpha}
\bibliography{diodon}

\end{document}